\documentclass[11pt,a4paper]{article}
\usepackage{graphicx}
\usepackage{amsmath}
\usepackage{amssymb}
\usepackage{amsthm}
\usepackage{geometry}
\usepackage{fancyhdr}
\usepackage{color} 
\usepackage{overpic}
\usepackage{mathabx}

\usepackage[all]{xy}
\newcommand{\al}{\alpha}

\newcommand{\be}{\beta}
\newcommand{\ga}{\gamma}
\newcommand{\Ga}{\Gamma}
\newcommand{\ra}{\rightarrow}
\newcommand{\xra}{\xrightarrow}
\newcommand{\rgl}{\rangle}
\newcommand{\lgl}{\langle}

\newcommand{\real}{\mathbb{R}} 
\newcommand{\intg}{\mathbb{Z}} 

\newcommand{\bpf}{\begin{proof}}
\newcommand{\epf}{\end{proof}}
\newcommand{\bthm}{\begin{thm}}
\newcommand{\ethm}{\end{thm}}
\newcommand{\bprop}{\begin{prop}}
\newcommand{\eprop}{\end{prop}}
\newcommand{\bcor}{\begin{cor}}
\newcommand{\ecor}{\end{cor}}
\newcommand{\blem}{\begin{lem}}
\newcommand{\elem}{\end{lem}}
\newcommand{\bdefn}{\begin{defn}}
\newcommand{\edefn}{\end{defn}}
\newcommand{\bexmp}{\begin{exmp}}
\newcommand{\eexmp}{\end{exmp}}
\newcommand{\brem}{\begin{rem}}
\newcommand{\erem}{\end{rem}}
\newcommand{\beq}{\begin{equation*}\begin{aligned}}
\newcommand{\eeq}{\end{aligned}\end{equation*}}

\newcommand{\shm}{\underline{\rm SHM}}
\newcommand{\shi}{\underline{\rm SHI}}
\newcommand{\khm}{\underline{\rm KHM}^-}
\newcommand{\khi}{\underline{\rm KHI}^-}

\newtheorem{thm}{\textbf {Theorem}}[section]
\newtheorem{cor}[thm]{\textbf{Corollary}}
\newtheorem{prop}[thm]{\textbf{Proposition}}
\newtheorem{lem}[thm]{\textbf{Lemma}}
\newtheorem{conj}[thm]{Conjecture}

\newtheorem{quest}[thm]{Question}

\theoremstyle{definition}
\newtheorem{defn}[thm]{\textbf{Definition}}

\newtheorem{exmp}[thm]{Example}

\theoremstyle{remark}
\newtheorem{rem}[thm]{Remark}

\title{Decomposing sutured monopole and instanton Floer homologies}

\author{Sudipta Ghosh and Zhenkun Li}
\date{}

\begin{document}
\bibliographystyle{plain}

\maketitle

\abstract{In this paper, we generalize the work of the second author in \cite{li2019direct} and prove a grading shifting property, in sutured monopole and instanton Floer theories, for general balanced sutured manifolds. This result has a few consequences. First, we offer an algorithm that computes the Floer homologies of a family of sutured handlebodies. Second, we obtain a canonical decomposition of sutured monopole and instanton Floer homologies and build polytopes for these two theories, which was initially achieved by Juh\'asz \cite{juhasz2010polytope} for sutured (Heegaard) Floer homology. Third, we establish a Thurston-norm detection result for monopole and instanton knot Floer homologies, which were introduced by Kronheimer and Mrowka in \cite{kronheimer2010knots}. The same result was originally proved by Ozsv\'ath and Szab\'o for link Floer homology in \cite{ozsvath2008thurston}. Last, we generalize the construction of minus versions of monopole and instanton knot Floer homology, which was initially done for knots by the second author in \cite{li2019direct}, to the case of links. Along with the construction of polytopes, we also proved that, for a balanced sutured manifold with vanishing second homology, the rank of the sutured monopole or instanton Floer homology bounds the depth of the balanced sutured manifold. As a corollary, we obtain an independent proof that monopole and instanton knot Floer homologies, as mentioned above, both detect fibred knots in $S^3$. This result was originally achieved by Kronheimer and Mrowka in \cite{kronheimer2010knots}.}


\section{Introduction}
Sutured manifold theory and Floer theory are two powerful tools in the study of $3$-dimensional topology. Sutured manifolds were first introduced by Gabai in \cite{gabai1983foliations} and in subsequent papers. The core of the sutured manifold theory is the sutured manifold hierarchy. This enables one to decompose any taut sutured manifold, in finitely many steps, into product sutured manifolds, which are the simplest possible ones. Gabai used sutured manifolds and sutured manifold hierarchies to prove some important results about $3$-manifolds, including the remarkable property R conjecture. 

Sutured (Heegaard) Floer homology was first introduced by Juh\'asz in \cite{juhasz2006holomorphic}, while some ad hoc versions were studied by Ghiggini in \cite{ghiggini2008knot} and Ni in \cite{ni2007knot}. In particular, Ni proved the celebrated result that the knot Floer homology, which was introduced by Ozsv\'ath and Sz\'abo in \cite{ozsvath2004holomorphicknot}, detects fibred knots. Ni's result is equivalent to the fact that sutured Floer homology detects product balanced sutured manifolds among homology products, which was then generalized by Juh\'asz in \cite{juhasz2010polytope}, where he proved that the rank of sutured Floer homology bounds the depth of a balanced sutured manifold with vanishing second homology. 

The combination of sutured manifold theory with gauge theory was done by Kronheimer and Mrowka in \cite{kronheimer2010knots}, where they defined sutured monopole and instanton Floer homologies. These new Floer homologies have many significant applications in the study of knots and $3$-manifolds, including a new and simpler proof of the famous property P conjecture. In \cite{kronheimer2010knots}, Kronheimer and Mrowka proved the following, in correspondence to Ni's result.

\bthm[Kronheimer and Mrowka \cite{kronheimer2010knots}]\label{thm: product detection}
Suppose $(M,\ga)$ is a balanced sutured manifold and is a homology product. Suppose further that
$${\rm rk}(\shm(M,\ga))=1~{\rm or}~{\rm rk}(\shi(M,\ga))=1.$$
Then, $(M,\ga)$ is a product sutured manifold.
\ethm

Theorem \ref{thm: product detection} has many important applications. For instance, this theorem is crucial in the proof that Khovanov homology detects unknots, by Kronheimer and Mrowka \cite{kronheimer2011khovanov}, and that Khovanov homology detects trefoils, by Baldwin and Sivek \cite{baldwin2018khovanov}. In this paper, we generalize Theorem \ref{thm: product detection} and prove the following:

\begin{thm}\label{thm: SHM and SHI bounds depth}
    Suppose $(M,\ga)$ is a taut balanced sutured manifold, $H_2(M)=0$, and 
    $${\rm rk}(\shm(M,\ga))<2^{k+1}~{\rm or}~{\rm rk}(\shi(M,\ga))<2^{k+1}.$$
    Then,
    $$d(M,\ga)\leq 2k.$$
Here $d(M,\ga)$ is the depth of a balanced sutured manifold, i.e., the minimal number of taut sutured manifold decompositions making $(M,\ga)$ into a product sutured manifold.
\end{thm}

\brem
If a balanced sutured manifold is a homology product, then $H_2(M)=0$. The converse is not necessarily true. Also, $d(M,\ga)=0$ if and only if $(M,\ga)$ is a product sutured manifold.
\erem

As a direct corollary to Theorem \ref{thm: SHM and SHI bounds depth}, we offer a new proof to the following well-known fact.
\bthm[Kronheimer and Mrowka \cite{kronheimer2010knots}]
The monopole and instanton knot Floer homologies, $KHM$ and $KHI$, as defined in \cite{kronheimer2010knots}, both detect fibred knots.
\ethm

For a balanced sutured manifold $(M,\ga)$ and a tangle $T$ inside $(M,\ga)$, Xie and Zhang \cite{xie2019tangles} constructed a version of sutured instanton Floer homology on balanced sutured manifolds with tangles, which they denote $SHI(M,\ga,T)$. In \cite{xie2019classification}, they used their construction as a tool to fully classify links whose Khovanov homologies have minimal possible ranks. One crucial step in their proofs is to show that the sutured instanton Floer homology they constructed detects product tangles inside product sutured manifolds. In this paper, with Theorem \ref{thm: SHM and SHI bounds depth}, we can prove a slightly more general result than their product-tangle-detection theorem.

\bcor\label{cor: product tangle detection}
Suppose $(M,\ga)$ is a balanced sutured manifold equipped with a vertical tangle $T$. Suppose further that $H_2(M\backslash T)=0$ and 
$SHI(M,\ga,T)\cong\mathbb{C}$. Then, $(M,\ga,T)$ is diffeomorphic to a product sutured manifold equipped with a product tangle, i.e.,
$$(M,\ga,T)\cong([-1,1]\times F,\{0\}\times\partial{F}, [-1,1]\times\{p_1,...,p_n\}).$$
Here, $F$ is a compact oriented surface-with-boundary, and $p_1,...,p_n$ are distinct points on $F$.
\ecor

\brem
Note we only need the requirement that $H_2(M\backslash T)=0$, instead of the assumption that $H_2(M)=0$ in the hypothesis of Theorem \ref{thm: SHM and SHI bounds depth}. The latter condition is stronger when tangles do exist. Our assumption is also weaker than the original assumption that $(M,\ga)$ must be a homology product in Xie and Zhang \cite{xie2019tangles}.
\erem

Despite the applications mentioned above, there are many basic aspects of sutured monopole and instanton Floer theories that remain in mystery. The usual monopole and instanton theories were defined on closed oriented $3$-manifolds, while balanced sutured manifolds are compact oriented manifolds with non-trivial boundaries. So, to define the sutured monopole and instanton Floer homologies, Kronheimer and Mrowka constructed a special class of closed oriented $3$-manifolds, called closures, out of the sutured data. However, the choice of closure is not unique, which lead to the following two questions:

\begin{quest}\label{quest: 1}
In \cite{kronheimer2010instanton}, Kronheimer and Mrowka proved that different closures give rise to isomorphic sutured monopole and instanton Floer homologies. Then, to what extent can we say that all of the essential information of sutured monopole and instanton Floer homologies is contained in the original balanced sutured manifold rather than the full closure?
\end{quest}

\begin{quest}\label{quest: 2}
The monopole Floer homology on a closed $3$-manifold decomposes along spin${}^c$ structures (see \cite{kronheimer2007monopoles}). Correspondingly, the instanton Floer homology decomposes along eigenvalue functions (see \cite{kronheimer2010knots}). Then, do we have a spin${}^c$-type decomposition for sutured monopole or instanton Floer homology?
\end{quest}

\begin{quest}\label{quest: 3}
How do sutured monopole and	instanton Floer homologies tell us information about the Thurston norm on a balanced sutured manifold?
\end{quest}

Towards answering the first question, the second author proved in \cite{li2019direct} the following proposition:
\bprop[Li \cite{li2019direct}]\label{prop: 1 1}
Suppose $(M,\ga)$ is a balanced sutured manifold with a toroidal boundary, and $\ga$ consists of two components. Suppose further that $Y$ is a closure of $(M,\ga)$, and $\mathfrak{s}_1$ and $\mathfrak{s}_2$ are two supporting spin${}^c$ structures on $Y$, then there is a $1$-cycle $x$ in $M$ so that
$$c_1(\mathfrak{s}_1)-c_1(\mathfrak{s}_2)=P.D.[x]\in H^2(Y).$$
Similar statements hold in the instanton settings.
\eprop

This theorem is central to the second author's proof of a grading shifting property for gradings associated to properly embedded surfaces inside those balanced sutured manifolds that are described in the hypothesis of Proposition \ref{prop: 1 1}. The grading shifting property has two consequences in that paper. The first is to compute the sutured monopole and instanton Floer homologies of any sutured solid tori. The second is to construct an Alexander grading on the minus versions of monopole and instanton knot Floer homologies, as well as proving many fundamental properties of them. 

However, in the hypothesis of Proposition \ref{prop: 1 1}, it is required that $M$ has a toroidal boundary and that the suture has only two components. These requirements are very restrictive. For instance, one cannot use Proposition \ref{prop: 1 1} to construct minus versions of monopole and instanton knot Floer homologies for links. This is because the sutured manifolds arising from links may have more than one boundary component, and, thus, Proposition \ref{prop: 1 1} does not apply. In this paper, we prove the same result, as in Proposition \ref{prop: 1 1}, for any balanced suture manifolds.
\bthm\label{thm: difference of supporting spin c structures}
Suppose $(M,\ga)$ is a balanced sutured manifold. Suppose further that $Y$ is a closure of $(M,\ga)$, and $\mathfrak{s}_1$ and $\mathfrak{s}_2$ are two supporting spin${}^c$ structures on $Y$. Then, there is a $1$-cycle $x$ in $M$ so that
$$c_1(\mathfrak{s}_1)-c_1(\mathfrak{s}_2)=P.D.[x]\in H^2(Y).$$

A similar result holds in the instanton settings.
\ethm

Thus, to answer Question \ref{quest: 1}, we could say that, in any closure, the difference of any two supporting spin${}^c$ structures, in terms of the Poincar\'e dual of their first Chern classes, is contained in the original balanced sutured manifold instead of the whole closure. Also, Theorem \ref{thm: difference of supporting spin c structures} leads to a generalization of the grading shifting property, which was initially discussed in Li \cite{li2019direct}, as follows.

\bthm\label{thm: general grading shifting property}
Suppose $(M,\ga)$ is a balanced sutured manifold and $\al\in H_2(M,\partial M)$ is a non-trivial homology class. Pick two surfaces $S_1$ and $S_2$ so that 
$$[S_1,\partial S_1]=[S_2,\partial{S}_2]=\al\in H_2(M,\partial{M}),$$
and they are both admissible (see Definition \ref{defn: admissibility of surfaces}) in $(M,\ga)$. Then, there exist constants $l_M,l_I\in \intg$, so that, for any $j\in\intg$, we have:
$$\shm(-M,-\ga,S_1,j)=\shm(-M,-\ga,S_2,j-l_M),$$
and
$$\shi(-M,-\ga,S_1,j)=\shi(-M,-\ga,S_2,j-l_I).$$
Note the value $l_I$ and $l_M$ depends on the surfaces $S_1$ and $S_2$.
\ethm

The general grading shifting property given by Theorem \ref{thm: general grading shifting property} helps compute the sutured monopole and instanton Floer homology of some families of sutures on a general handlebody. In section \ref{sec: general grading shifting property}, we use a concrete example to present the algorithm. Theorem \ref{thm: general grading shifting property} also leads to a generalization of the minus version of monopole and instanton knot Floer homologies for links:
\bthm
Suppose $L\subset Y$ is a link so that each component of $L$ is null-homologous. Suppose further that $L$ has $r$ components and $\textbf{p}$ is an $r$-tuple, consisting of one point on each component of $L$. Then, associated to the triple $(-Y,L,\textbf{p})$, we can construct an infinite-rank module $\khm(-Y,L,\textbf{p})$ over the rings $\mathcal{R}$, the mod $2$ Novikov ring (c.f. \cite[Remark 2.19]{baldwin2016contact}). Moreover, $\khm(-Y,L,\textbf{p})$ is well defined only up to multiplication by a unit in $\mathcal{R}$ and has the following properties.

(1) Suppose $\{S_1,...,S_r\}$ is a collection of $r$ Seifert surfaces, one for each component of $L$, then they together induce a $\intg^r$ grading on $\khm(-Y,L,\textbf{p})$ and $\khi(-Y,L,\textbf{p})$.

(2) For each $i\in\{1,...,r\}$, there is a morphism
$$U_i:\khm(-Y,L,\textbf{p})\ra \khm(-Y,L,\textbf{p}),$$
which drops the grading associated to $S_i$ by $1$ and preserves all other gradings. This makes $\khm(-Y,L,\textbf{p})$ a module over $\mathcal{R}[U_1,...,U_r]$.

(3) There exists $N\in \intg$ so that if $\textbf{j}=(j_1,...,j_r)\in\intg^r$ is a multi-grading and $j_i<N$ for some $i\in\{1,...,r\}$, then the morphism $U_i$ restricts to an isomorphism
$$U_i:\khm(-Y,L,\textbf{p},\textbf{j})\xra{\cong}\khm(-Y,L,\textbf{p},\textbf{j}^{\dag}),$$
where $\textbf{j}^{\dag}$ is obtained from $\textbf{j}$ by replacing $j_i$ with $j_i-1$.

Furthermore, using instanton theory, we can construct $\khi(-Y,L,\textbf{p})$, which is well defined up to a multiplication by a non-zero element in the field of complex numbers and properties (1), (2), and (3) all hold.
\ethm

To answer Question \ref{quest: 2}, we construct a canonical decomposition of sutured monopole and instanton Floer homologies, independent of the choices of closures. To ensure this decomposition is canonical, we need to pre-fix an element of a special type inside the sutured monopole or instanton Floer homology, which we call a homogenous element.
\bprop\label{prop: canonical decomposition}
Suppose $(M,\ga)$ is a balanced sutured manifold and $a\in\shm$ is a homogenous element. Then there is a canonical decomposition
$$\shm(M,\ga)=\bigoplus_{\rho\in H^2(M,\partial{M};\real)}\shm_a(M,\ga,\rho).$$

A similar statement holds in the instanton settings.
\eprop

\brem
See Definition \ref{defn: homogenous element} for details about the definition of homogeneous elements. Also, from the discussion in Section \ref{subsec: polytope}, as long as the Floer homology group is non-trivial, homogeneous elements always exist.
\erem

Thus, we could define polytopes for sutured monopole and instanton Floer theories, to be the convex hulls of sets of $\rho$ so that $\shm_a(M,\ga,\rho)\neq0$ or $\shi_a(M,\ga,\rho)\neq0$. The first definition of such a polytope was introduced by Juh\'asz in the context of sutured Floer theory. In this paper, we also proved the following:
\bcor\label{cor: 1_8}
Suppose $(M,\ga)$ is a taut balanced sutured manifold with $H_2(M)=0$. Suppose further that $(M,\ga)$ is reduced, horizontally prime, and free of essential product disks. Then, the polytopes must both have maximal possible dimensions. In particular, we conclude that
$${\rm rk}_{\mathcal{R}}(\shm(M,\ga))\geq b^1(M)+1,~{\rm and}~{\rm dim}_{\mathbb{C}}(\shi(M,\ga))\geq b^1(M)+1.$$
\ecor

The proofs of Corollary \ref{cor: 1_8} and Theorem \ref{thm: SHM and SHI bounds depth} both rely on a technical result proved in section \ref{sec: polytopes}, which describes in detail how sutured monopole and instanton Floer homologies behave under sutured manifold decompositions. It is closely related to the decomposition theorem, Proposition 6.9, in Kronheimer and Mrowka \cite{kronheimer2010knots}.

The polytopes we defined for sutured monopole and instanton theories are closely related to the Thurston norms on the original balanced sutured manifold as well as on the closures. In particular, the canonical decomposition in Proposition \ref{prop: canonical decomposition} and the grading shifting property in Theorem \ref{thm: general grading shifting property} enable us to prove a Thurston norm detection result for monopole and instanton knot Floer homologies. The same result was previously achieved by Ozsv\'ath and Sz\'abo in \cite{ozsvath2008thurston}, for the link Floer homology in Heegaard Floer theory. 

Suppose $Y$ is a closed oriented $3$-manifold and $L\subset Y$ is an oriented link. Let
$$L=L_1\cup...\cup L_r$$
be the components of $L$. We require the following two conditions to hold for $L\subset Y$.

(1) The link complement, $Y(L)=Y\backslash N(L)$, is irreducible.

(2) The link complement $Y(L)$ is boundary-incompressible.

Let $\Ga_{\mu}\subset \partial Y(L)$ be the suture consisting of a pair of oppositely oriented meridians on each boundary component of $Y(L)$. Then, by Proposition \ref{prop: canonical decomposition}, there is a decomposition
$$\underline{\rm KHM}(Y,L)=\shm(Y(L),\Ga_{\mu})=\bigoplus_{\rho\in H^2(Y(L),\partial{Y(L)};\mathbb{Q})}\shm_a(Y(L),\Ga_{\mu},\rho).$$
We make the following definition.
\bdefn\label{defn: the function y}
For a class $\al\in H_2(Y(L),\partial Y(L))$, define
$$y(\al)=\max_{\substack{\rho\in H^2(Y(L),\partial{Y(L)};\mathbb{Q})\\\shm_a(Y(L),\Ga_{\mu},\rho)\neq 0}}\{\rho(\al)\}-\min_{\substack{\rho\in H^2(Y(L),\partial{Y(L)};\mathbb{Q})\\\shm_a(Y(L),\Ga_{\mu},\rho)\neq 0}}\{\rho(\al)\}$$
\edefn

\bthm\label{thm: Thurston norm detection}
Under the above settings,
\begin{equation}\label{eq: Thurston norm detection}
x(\al)+\sum_{i=1}^r|\lgl\al,\mu_i\rgl|=y(\al).
\end{equation}
Here $x(\cdot)$ is the Thurston-norm defined in Definition \ref{defn: thurston norm}. $\lgl,\rgl$ is to take the algebraic intersection number of a class $\al\in H_2(Y(L),\partial Y(L))$ with a class $[\mu_i]\in H_1(Y(L))$, where $\mu_i$ is a meridian of the link component $L_i$.
\ethm

\brem
Theorem \ref{thm: Thurston norm detection} offers a complete answer to Question \ref{quest: 3} in the case when the boundary of the balanced sutured manifold consists of tori. Suppose $M$ is a connected compact oriented $3$-manifold so that $M$ is irreducible and boundary-incompressible, and its boundary consists of tori, then we can perform Dehn fillings on each boundary component of $M$. The cores of the Dehn fillings give us a Link inside the resulting closed $3$-manifold $Y$, which satisfies the hypothesis of Theorem \ref{thm: Thurston norm detection}. Since the Dehn surgery can be performed along any non-separating simple closed curves, Theorem \ref{thm: Thurston norm detection} covers all the cases when the suture $\ga$ on $\partial{M}$ consists of a pair of non-separating simple closed curves on each boundary component of $M$. For a more general suture, when it may have more than two-component on some boundary component of $M$, we can modify the coefficients of the terms $\lgl h,\mu_i\rgl$ according to the number of components, and the proof of Theorem \ref{thm: Thurston norm detection} still applies verbatim.
\erem

{\bf Organization.} In Section \ref{sec: preliminaries}, we include the basic definitions and known results that support the proofs in this paper. In Section \ref{sec: difference of supporting spin c structures}, we study the set of supporting spin${}^c$ structures on any closure of a balanced sutured manifold. This will be the basis for a generalized grading shifting result proven in Section \ref{sec: general grading shifting property}. In Section \ref{sec: general grading shifting property}, we also offer an algorithm that could potentially compute the sutured monopole and instanton Floer homologies of a family of sutured handlebodies by carrying it out on a particular example. In Section \ref{sec: polytopes}, we use the generalized grading shifting property to construct a canonical decomposition of sutured monopole or instanton Floer homology and further construct polytopes in these two theories. Also, we prove some basic properties of the polytopes as well as the results regarding the depth of a sutured manifold. In Section \ref{sec: application to links}, we present some applications to knots and links: The first is to prove the Thurston-norm detection result for link complements, and the second is to construct minus versions for links.

{\bf Acknowledgement.} The authors would like to thank their advisors David Shea Vela-Vick and Tomasz Mrowka, for their enormous help. The authors would like to thank Andr\'as Juh\'asz and Yi Ni for helpful comments or conversations. The first author was supported by his advisor David Shea Vela-Vick's NSF Grant 1907654 and Simons Foundation Grant 524876. The second author was supported by his advisor Tom Mrowka’s NSF Grant 1808794.

\section{Preliminaries}\label{sec: preliminaries}
\subsection{Basic definitions of Balanced sutured manifolds}
\bdefn\label{defn: balanced sutured manifold}
A {\it balanced sutured manifold} is a pair $(M,\ga)$ consisting of a compact oriented $3$-manifold $M$ and a closed oriented $1$-submanifold $\ga\subset\partial M$. On $\partial{M}$, let $A(\ga)=[-1,1]\times\ga$ be an annular neighborhood of $\ga\subset \partial{M}$, and let 
$$R(\ga)=\partial{M}\backslash{\rm int}(A(\ga)).$$
They satisfy the following requirements.

(1) Neither $M$ nor $R(\ga)$ has closed components.

(2) If we orient $\partial{R}(\ga)=\partial{A}(\ga)=\{\pm1\}\times\ga$ in the same way as $\ga$, then we require that the orientation on $\partial{R}(\ga)$ induces a unique orientation on $R(\ga)$. This orientation is called the {\it canonical orientation} on $R(\ga)$. We use $R_+(\ga)$ to denote the part of $R(\ga)$ whose canonical orientation coincides with the boundary orientation of $\partial{M}$, and $R_-(\ga)$ the rest. 

(3) We require that
$$\chi(R_+(\ga))=\chi(R_-(\ga)).$$
\edefn

\begin{defn}
A balanced sutured manifold $(M, \ga)$ is called a {\it product sutured manifold} if $M = [-1,1]\times R$, $A(\ga)= [-1,1] \times \partial R$, $R_+(\ga) = \{1\}\times R$, and $R_-(\ga) = \{-1\} \times R$. Here, $R$ is compact oriented surface with no closed components.    
\end{defn}

\bdefn
Suppose $M$ is a compact oriented $3$-manifold. $M$ is called {\it irreducible} if every embedded $2$-sphere $S^2\subset M$ bounds an embedded $3$-ball inside $M$.
\edefn

\bdefn
Suppose $M$ is a compact $3$-manifold and $R\subset M$ is an embedded surface. $R$ is called {\it compressible} if there is a simple closed curve $\al\subset R$ so that $\al$ does not bound a disk on $R$ but bounds an embedded disk $D\subset M$ with $D\cap R=\al$. $R$ is called {\it incompressible} if it is not compressible. A $3$-manifold is called {\it boundary-incompressible} if its boundary is incompressible.
\edefn

\bdefn[Thurston norm]\label{defn: thurston norm}
Suppose $M$ is a compact $3$-manifold, and suppose $U\subset \partial{M}$ is a submanifold of $\partial{M}$. Suppose further that $S$ is a properly embedded surface inside $M$ so that $\partial{S}\subset U$. If $S$ is connected, then define the {\it norm} of $S$ to be
$$x(S)=\max\{-\chi(S),0\}.$$
In general, suppose the components of $S$ are
$$S=S_1\cup...\cup S_n,$$
where each $S_i$ is connected, then define the {\it norm} of $S$ to be
$$x(S)=x(S_1)+...+x(S_n).$$
Moreover, suppose $\al\in H_2(M,U)$ is a non-trivial second relative homology class, then define the {\it norm} of $\al$ to be
$$x(\al)=\min\{x(S)~|~(S,\partial{S})\subset (M,U),~[S,\partial{S}]=\al\in H_2(M,U)\}.$$
\edefn

\bdefn
Suppose $M$ is a compact $3$-manifold, and $S\subset M$ is a properly embedded surface. $S$ is called {\it norm-minimizing} if
$$x(S)=x(\al),$$
where $\al=[S,\partial{S}]\in H_2(M,N(\partial{S}))$. Here, $N(\partial{S})$ is a neighborhood of $\partial{S}\subset \partial{M}$.
\edefn

\bdefn[Gabai \cite{gabai1983foliations}]
A balanced sutured manifold $(M,\ga)$ is called {\it taut} if the following is true.

(1) $M$ is irreducible.

(2) $R_{+}(\ga)$ and $R_{-}(\ga)$ are both incompressible.

(3) $R_{+}(\ga)$ and $R_{-}(\ga)$ are both norm-minimizing.
\edefn

\bdefn[Gabai \cite{gabai1983foliations}]
Let $(M, \ga)$ be a balanced sutured manifold. A {\it product annulus} $A$ in $(M, \ga)$ is an annulus properly embedded in $M$ such that $\partial A \subset R(\ga)$ and $\partial A \cap R_{\pm}(\ga) \neq \emptyset$. A product disk is a disk $D$ properly embedded in $M$ such that $\partial D \cap A(\ga)$ consists of two essential arcs in $A(\ga)$.
\edefn

Product annuli and product disks can detect where $(M,\ga)$ is locally a product. We have the following definition following Juh\'asz \cite{juhasz2010polytope}.

\bdefn[Juh\'asz \cite{juhasz2010polytope}]
A balanced sutured manifold $(M,\ga)$ is called reduced if any product annulus $A\subset M$ either bounds a cylinder $D^2\times I$ so that $\partial{D^2}\times I=A$, or is isotopic to a component of $A(\ga)$ inside $M$.
\edefn

\bdefn[Gabai \cite{gabai1983foliations}]
Let $(M, \ga)$ be a taut balanced sutured manifold. A properly embedded surface $S \subset M$ is called {\it horizontal} if the following four properties hold.

(1) $S$ has no closed components and is incompressible.

(2) $\partial S \subset A(\ga)$, and $\partial S$ is parallel to $\partial{R}_+(\ga)$ inside $A(\ga)$.

(3) $[S]=[R_+(\ga)]$ in $H_2(M, A(\ga))$.

(4) $\chi(S) = \chi (R_+(\ga))$.

We say that $(M, \ga)$ is {\it horizontally prime} if every horizontal surface in $(M, \ga)$  is parallel to either $R_+(\ga)$ or $R_-(\ga)$.
\edefn

\bdefn[Gabai \cite{gabai1983foliations}]
Suppose $(M,\ga)$ is a taut balanced sutured manifold. The {\it depth} of $(M,\ga)$, which we write $d(M,\ga)$, is the minimal integer $n$ so that there exists a sequence of sutured manifold decompositions (for definitions, see Gabai \cite{gabai1983foliations} or Scharlemann \cite{scharlemann1989suture})
$$(M_0,\ga_0)\stackrel{S_0}{\leadsto}(M_1,\ga_1)\stackrel{S_1}{\leadsto}...\stackrel{S_n}{\leadsto}(M_{n+1},\ga_{n+1}),$$
so that each $(M_i,\ga_i)$ is taut, $(M_0,\ga_0)=(M,\ga)$ and $(M_{n+1},\ga_{n+1})$ is a product sutured manifold.
\edefn

\bthm[Gabai \cite{gabai1983foliations}]
For any taut balanced sutured manifold, its depth is finite.
\ethm

\subsection{Monopole and instanton Floer homologies on balanced sutured manifolds}\label{subsec: SHM and SHI}

To define sutured monopole and instanton Floer homologies, one needs to construct a closed $3$-manifold, together with a distinguishing surface, out of a balanced sutured manifold $(M,\ga)$. To do this, pick $T$ to be a connected oriented surface so that the following is true.

(1) There is an orientation reversing diffeomorphism
$$f:\partial{T}\ra \ga.$$

(2) $T$ has genus at least $2$.

(3) There is a fixed base point $p\in T$.

Then, we can use $f$ to glue $T\times[-1,1]$ to $M$, along the annuli $A(\ga)$, and let
$$\widetilde{M}=M\mathop{\cup}_{id\times f}[-1,1]\times T.$$
The manifold $\widetilde{M}$ has two boundary components:
$$\partial\widetilde{M}=R_+\cup R_-,$$
where
$$R_{\pm}=R_{\pm}(\ga)\mathop{\cup}_{f} \{\pm 1\}\times T.$$
Let $h:R_+\ra R_-$ be an orientation preserving diffeomorphism so that $h(\{1\}\times\{p\})=\{-1\}\times\{p\}$, then we can form a closed $3$-manifold as follows:
$$Y=\widetilde{M}\mathop{\cup}_{id\cup h}[-1,1]\times R_{+}.$$
Here, we use $h$ to glue $\{1\}\times R_{+}$ to $R_-\subset \partial{\widetilde{M}}$ and use the identity map to glue $\{-1\}\times R_+$ to $R_+\subset \partial\widetilde{M}$. Let $R=\{0\}\times R_+\subset Y$, and we make the following definition. 

\bdefn\label{defn: closure}
The manifold $\widetilde{M}$ is called a {\it pre-closure} of $(M,\ga)$. The pair $(Y,R)$ is called a {\it closure} of the balanced sutured manifold $(M,\ga)$. The choices $T,f$, and $h$ are called the auxiliary data. In particular, the surface $T$ is called an {\it auxiliary surface} and $h$ a {\it gluing diffeomorphism}.
\edefn

\bdefn\label{defn: SHM}
The {\it sutured monopole Floer homology} of $(M,\ga)$ is defined as
$$SHM(M,\ga)=HM(Y|R)=\bigoplus_{\mathfrak{s}\in\mathfrak{S}^*(Y|R)}\widecheck{HM}_{\bullet}(Y,\mathfrak{s};\Ga_{\eta}).$$
Here, $\eta\subset R$ is a non-separating simple closed curve and
$$\mathfrak{S}^*(Y|R)=\{\mathfrak{s}~{\rm spin}^c~{\rm structures~on~Y},~c_1(\mathfrak{s})[R]=2g(R)-2,~\widecheck{HM}_{\bullet}(Y,\mathfrak{s})\neq 0\}$$
is called {\it the set of supporting spin${}^c$ structures}. We use the mod 2 Novikov ring $\mathcal{R}$ for coefficients. The notation $\Ga_{\eta}$ denotes the local system in \cite[Section 2.2]{kronheimer2007monopoles}.
\edefn

\begin{rem}
Note from \cite{kronheimer2010knots}, we know that
$$\widecheck{HM}_{\bullet}(Y,\mathfrak{s};\Ga_{\eta})\cong \widecheck{HM}_{\bullet}(Y,\mathfrak{s})\otimes\mathcal{R}$$
for any spin${}^c$ structure $\mathfrak{s}$ on $Y$ so that
$$c_1(\mathfrak{s})[R]=2g(R)-2.$$	
\end{rem}

\bdefn\label{defn: supporting eigenvalue function}
The {\it sutured instanton Floer homology} of $(M,\ga)$ is defined as
$$SHI(M,\ga)=I(Y|R)=\bigoplus_{\lambda\in \mathfrak{H}^*(Y|R)} I^{\omega}(Y)_{\lambda}.$$
Here, $\omega\subset Y$ is a simple closed curve having a unique transverse intersection with $R$ (c.f. \cite[Section 7]{kronheimer2010knots}), and the notation $I^{\omega}(Y)_{\lambda}$ follows from \cite[Section 6]{li2018gluing}. The set
$$\mathfrak{H}^*(Y|R)=\{\lambda\in H^2(Y;\mathbb{Q}),~\lambda[R]=2g(R)-2,~I^{\omega}(Y)_{\lambda}\neq0\}$$
is called {\it the set of supporting eigenvalue functions}. We use the field of complex numbers $\mathbb{C}$ for coefficients.
\edefn

\bthm\label{thm: isomorphism class of SHM and SHI}[Kronheimer and Mrowka \cite{kronheimer2010knots}]
The isomorphism classes of $SHM(M,\ga)$ and $SHI(M,\ga)$ are invariants for a fixed balanced sutured manifold $(M,\ga)$.
\ethm

Only knowing that the isomorphism class is invariant is sometimes not enough. In \cite{baldwin2015naturality}, Baldwin and Sivek refined the definition of closures and constructed canonical maps between different closures. In particular, they proved the following.

\bthm
Suppose $(M,\ga)$ is a balanced sutured manifold. Suppose further that $(Y_1,R_1)$ and $(Y_2,R_2)$ are two closures of $(M,\ga)$. Pick non-separating curves $\eta_1\subset R_1$ and $\eta_2\subset R_2$ to determine local coefficients. Then, there exists a map
$$\Phi_{1,2}:HM(Y_1|R_1;\Ga_{\eta_1})\xra{\cong}HM(Y_2|R_2;\Ga_{\eta_2}),$$
which is well-defined up to multiplication by a unit in the base ring. Moreover, it satisfies the following properties.

(1) If $(Y_1,R_1)=(Y_2,R_2)$, then $\Phi_{1,2}\doteq id$. Here, $\doteq$ means equal up to multiplication by a unit.

(2) If there is a third closure $(Y_3,R_3)$ of $(M,\ga)$, then
$$\Phi_{1,3}\doteq\Phi_{2,3}\circ\Phi_{1,2}.$$

Similar results hold in the instanton settings.
\ethm

Hence, for a balanced sutured manifold $(M,\ga)$, the sutured monopole or instanton Floer homologies of closures, together with the canonical maps, form a projective transitive system, and we can derive a canonical module which we denote by 
$$\shm(M,\ga)~{\rm or}~\shi(M,\ga).$$
They are well defined up to a unit in the corresponding base ring (field). For more details, readers are referred to Baldwin and Sivek \cite{baldwin2015naturality}.

\begin{figure}[h]
\centering
\begin{overpic}[width=5.0in]{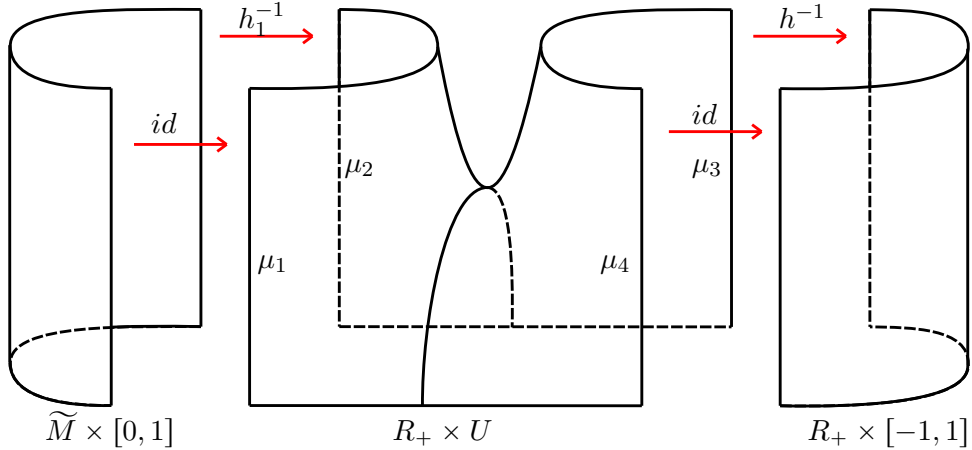}
    \put(4,-3){$\widetilde{M}\times[0,1]$}
    \put(40,-3){$R_{+}\times U$}
    \put(83,-3){$R_+\times[-1,1]$}
    \put(15,29){$id$}
    \put(24,40){$h_1^{-1}$}
    \put(71,30){$id$}
    \put(80,40){$h^{-1}$}
    \put(35,25){$\mu_2$}
    \put(26,15){$\mu_1$}
    \put(71,25){$\mu_3$}
    \put(61.5,15){$\mu_4$}
\end{overpic}
\vspace{0.05in}
\caption{Gluing three parts together to get $W$. The middle part is $U\times R_{+}$, while the $R_{+}$ directions shrink to a point in the figure.}\label{fig: excision cobordism}
\end{figure}

Floer excision is a very useful tool introduced by Kronheimer and Mrowka \cite{kronheimer2010knots} into the context of sutured monopole and instanton Floer homologies. It has several different versions, but we only present the version that is useful in later sections. Suppose $(M,\ga)$ is a balanced sutured manifold and $\widetilde{M}$ is a pre-closure of $(M,\ga)$. Recall that
$$\partial\widetilde{M}=R_+\cup R_-.$$
Suppose we use two gluing diffeomorphisms, $h_1$ and $h_2$, to obtain two closures $(Y_1,R_+)$ and $(Y_2,R_+)$, respectively. Let $h=h_1^{-1}\circ h_2$, and let $Y^h$ be the mapping torus of $h:R_+\ra R_+$ Then, we can form a cobordism $W$, which is from $Y_1\sqcup Y^h$ to $Y_2$, as follows. Let $U$ be a disk as depicted in Figure \ref{fig: excision cobordism}, and four parts of its boundary, $\mu_1,\mu_2,\mu_3,$ and $\mu_4$, are each identified with the interval $[0,1]$. Glue three pieces $[0,1]\times\widetilde{M}$, $ U\times R_+$ and $[0,1]\times[-1,1]\times R_+$ together, and let
$$W=([0,1]\times \widetilde{M})\mathop{\cup}_{id\cup h_1^{-1}}(U\times R_+) \mathop{\cup}_{id\cup h^{-1}}([0,1]\times[-1,1]\times R_+).$$
Here, we use $id$ to glue $[0,1]\times R_+$ to $\mu_1\times R_{+}$, use $h_1^{-1}$ to glue $[0,1]\times R_-$ to $\mu_2\times R_{+}$, use $h^{-1}$ to glue $\mu_3\times R_{+}$ to $[0,1]\times \{1\}\times R_-$ and use $id$ to glue $\mu_4\times R_+$ to $[0,1]\times \{-1\}\times R_{+}$. Note all the gluing maps are identity on the $[0,1]$ direction.

Pick a non-separating curve $\eta\subset T$, and suppose the diffeomorphisms $h_1$ and $h_2$ we choose at the beginning both preserve $\eta$: 
$$h_1(\{1\}\times\eta)=\{-1\}\times\eta,~{\rm and}~h_2(\{1\}\times\eta)=\{-1\}\times\eta.$$
Then, we can use $\eta$ to support local coefficients. 

\bthm[Kronheimer and Mrowka, \cite{kronheimer2010knots}]\label{thm: floer excision in monopoles}
The cobordism $W$ induces an isomorphism
$$\widecheck{HM}(W):HM(Y_1|R_+;\Ga_{\eta})\otimes HM(Y^h|R_+;\intg_2)\ra HM(Y_2|R_+;\Ga_{\eta}).$$
\ethm

There are three basic lemmas that are useful in later sections. Here we only present them in the monopole settings, but all of them have correspondences in the instanton settings.

\blem[Kronheimer and Mrowka \cite{kronheimer2010knots}]\label{lem: unique spin c structure}
The set of supporting spin${}^c$ structures $\mathfrak{S}^*(Y^h|R_+)$ consists of a unique element $\mathfrak{s}^h$. Moreover, with $\intg$ coefficients,
$$HM(Y^h|R_+)=\widecheck{HM}_{\bullet}(Y^h,\mathfrak{s}^h)\cong\intg.$$

When using local coefficients, pick a non-separating curve $\eta\subset R_+$ and suppose $\mathcal{R}$ be any suitable base ring for local coefficients, then
$$HM(Y^h|R_+;\Ga_{\eta})=\widecheck{HM}_{\bullet}(Y^h,\mathfrak{s}^h;\Ga_{\eta})\cong\mathcal{R}.$$
\elem

\blem[Kronheimer and Mrowka \cite{kronheimer2007monopoles}]\label{lem: adjunction inequality in monopoles}
Suppose $Y$ is a closed oriented $3$-manifold and $\mathfrak{s}$ is a spin${}^c$ structure on $Y$ so that there is a embedded oriented surface $R\subset Y$ so that $g(R)\geq1$, and $|c_1(\mathfrak{s})[R]|>2g(R)-2$. Then, we have
$$\widecheck{HM}_{\bullet}(Y,\mathfrak{s})=0.$$
Similarly, for any local coefficients that could possibly be used,
$$\widecheck{HM}_{\bullet}(Y,\mathfrak{s};\Ga_{\eta})=0.$$
\elem

\blem\label{lem: spin c structure shall extend}
Suppose $(W,\nu)$ is a cobordism from $Y$ to $Y'$. Suppose $\mathfrak{s}$ is a spin${}^c$ structure on $Y$ and $\mathfrak{s}'$ is a spin${}^c$ structure on $Y'$, so that
$$\widecheck{HM}(W)(\widecheck{HM}_{\bullet}(Y,\mathfrak{s}))\cap \widecheck{HM}_{\bullet}(Y',\mathfrak{s}')\neq \{0\},$$
then we know that
$$i_*(P.D.c_1(\mathfrak{s}))=(i')_*(P.D.c_1(\mathfrak{s'}))\in H_1(W).$$
Here $i:Y\ra W$ and $i':Y'\ra W'$ are the inclusions.
\elem

\subsection{A grading on sutured monopole and instanton Floer homology}\label{subsec: construction of gradings}
In \cite{li2019direct}, the second author constructed a grading on sutured monopole or instanton Floer homology, associated to a properly embedded surface with a connected boundary. We will present the construction in this subsection while dropping the condition that $\partial{S}$ is connected, using some new inputs from Kavi \cite{kavi2019pairing}.

\bdefn
 Suppose $(M, \gamma)$ is a balanced sutured manifold, and $S$ is a
properly embedded oriented surface. A {\it stabilization} of S is an isotopy of S to a surface $S'$ so that the isotopy creates a new pair of intersection points:
$$ \partial S' \cap \gamma = ( \partial S \cap \gamma) \cup \{p_+, p_{-}\}.$$
We require that there are arcs $\al \subset \partial S'$ and $\be\subset \gamma$, which are oriented in the same way as $\partial S'$ and $\gamma$, respectively, so that the following is true.

(1) We have $\partial \al= \partial \be= \{p_+, p_{-}\}$.

(2) The curves $\al$ and $\be$ cobound a disk $D$ so that
 ${\rm int}(D) \cap (\gamma \cup \partial S') = \emptyset$. 
 The stabilization is called {\it negative} if $D$ can be oriented so that
 $\partial D= \al \cup \be$ as oriented curves. It is called {\it positive} if $\partial D= (-\al) \cup \be$.
 
Denote by $S^{\pm k}$ the result of performing $k$ many positive or negative stabilizations of $S$.
 
\begin{figure}[h]
\centering
\begin{overpic}[width=4.0in]{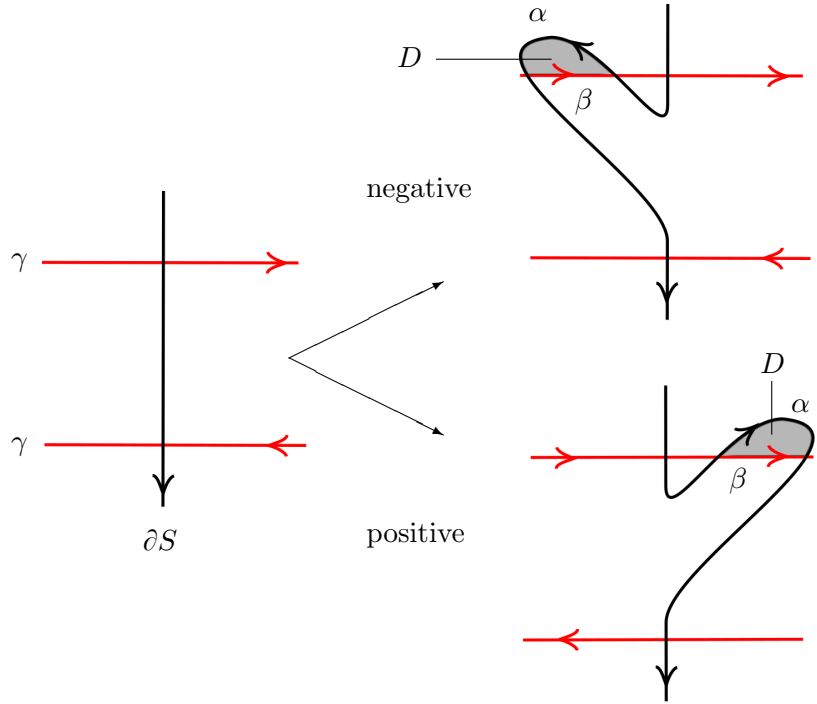}
    \put(13,20){$\partial{S}$}
    \put(-4,57){$\ga$}
    \put(-4,33){$\ga$}
    \put(32,45){\vector(2,1){20}}
    \put(32,45){\vector(2,-1){20}}
    \put(42,21){positive}
    \put(42,66){negative}
    \put(46,83){$D$}
    \put(51,84){\line(1,0){15}}
    \put(63,89){$\al$}
    \put(69,77.5){$\be$}
    \put(93,43){$D$}
    \put(94.5,42){\line(0,-1){7}}
    \put(97,38){$\al$}
    \put(89,28){$\be$}
\end{overpic}
\vspace{0.05in}
\caption{The positive and negative stabilizations of $S$.}\label{fig: pm stabilization of surfaces}
\end{figure}
\edefn

\blem\label{lem: decomposition and pm stabilization}
 Suppose $(M, \gamma)$ is a balanced sutured manifold, and S is a properly embedded oriented surface. Suppose further that $S^+$ and $S^-$ are obtained from $S$ by performing a positive and negative stabilization, respectively. Then, we have the following.
 
 (1) If we decompose $(M, \gamma)$ along $S$ or $S^+$, then the resulting two balanced sutured manifolds are diffeomorphic.
 
 (2) If we decompose $(M, \gamma)$ along $S^-$
, then the resulting balanced sutured manifold, $(M', \gamma')$, is not taut, because $R_{\pm}(\gamma')$ are both compressible.
\elem

\brem
The positive and negative stabilization on $S$ will be switched if we reverse the orientation of $S$, that is, $-(S^+)$ is the same as $(-S)^-$ and $-(S^-)$ is the same as $(-S)^+$. Accordingly, when changing the orientation of the suture, positive and negative stabilizations are also switched.
\erem

Now we present the construction of the grading. This was originally written down by Baldwin and Sivek in \cite{baldwin2018khovanov} and was then generalized one step further by the second author in \cite{li2019direct}, to fit his needs of constructing a $\intg$ grading in the minus version of monopole and instanton knot Floer homologies. In this subsection, we make a further generalization and introduce the most general setups of constructing such a grading.

\bdefn\label{defn: admissibility of surfaces}
Suppose $(M,\ga)$ is a balanced sutured manifold, and $S$ is a properly embedded surface inside $M$. Suppose further that $S$ intersects with $\ga$ transversely. $S$ is called {\it admissible} inside $(M,\ga)$ if every component of $\partial S$ intersects $\ga$ and the value $(\frac{1}{2}|S\cap \ga|-\chi(S))$ is an even integer.
\edefn

Suppose $(M,\ga)$ is a balanced sutured manifold, and $S\subset M$ is an oriented admissible properly embedded surface. Let $n=\frac{1}{2}|S\cap \ga|$. We fix an arbitrary ordering of the boundary components of $S$ and label all the intersection points of $S$ with $\ga$ as follows. On each boundary component of $S$, index them according to the orientation of $\partial{S}$. For different components of $\partial{S}$, we first index points on the boundary component that comes first in the fixed ordering. Also, the first point to be indexed on each boundary component of $S$ is chosen to be a positive intersection of $S$ with $\ga$ (on $\partial{M}$). In this way, we can assume that
$$S\cap \ga=\{p_1,...,p_{2n}\}.$$

In \cite{li2019direct}, when $S$ has a connected boundary, the second author introduced the notion of balanced pairings to help construct the grading. In Kavi \cite{kavi2019pairing}, the notion of balanced pairings was generalized to accommodate a general $S$. In this paper, we omit the detailed definitions of balanced pairings and will use the generalized definition from \cite{kavi2019pairing}.

Suppose 
$$\mathcal{P}=\{(i_1,j_1),...,(i_n,j_n)\}$$
is a balanced pairing of size $n$. Then, we can pick an auxiliary surface $T$ for $(M,\ga)$ so that the following is true.

(1) The genus of $T$ is large enough.

(2) The boundary of $T$ is identified with the suture $\ga$.

(3) There are properly embedded arcs $\al_1,...,\al_n$ inside $T$ so that the following two properties hold.

\quad (a) The classes $[\al_1],...,[\al_n]$ are linearly independent in $H_1(T,\partial T)$.

\quad (b) For $k=1,...,n$, we have
$$\partial {\al_k}=\{p_{i_k},p_{j_k}\}.$$

Then, we can form a pre-closure $\widetilde{M}$ of $(M,\ga)$:
$$\widetilde{M}=M\cup [-1,1]\times T.$$
The manifold $\widetilde{M}$ has two boundary components:
$$\partial{\widetilde{M}}=R_+\cup R_-.$$
The surface $S$ extends to a properly embedded surface $\widetilde{S}$ inside $\widetilde{M}$:
$$\widetilde{S}=S\cup [-1,1]\times\al_1\cup...\cup[-1,1]\times \al_n.$$

The definition of the balanced pairing makes sure that $\widetilde{S}\cap R_+$ and $\widetilde{S}\cap R_-$ have the same number of components, and the requirement (a) for $\al_i$ makes sure that components of $\widetilde{S}\cap R_{\pm}$ represent linearly independent classes in $H_1(R_{\pm})$. Thus, there exists an orientation preserving diffeomorphism $h: R_+\ra R_-$ so that
$$h(\widetilde{S}\cap R_+)=\widetilde{S}\cap R_-.$$
We can use $\widetilde{M}$ and $h$ to obtain a closure $(Y,R_+)$ of $(M,\ga)$, and, inside $Y$, the surface $S$ extends to a closed surface $\bar{S}$.

\bdefn\label{defn: grading}
Define
$$SHM(M,\ga, S, i)=\bigoplus_{\substack{\mathfrak{s}\in\mathfrak{S}^*(Y|R)\\c_1(\mathfrak{s})[\bar{S}]=2i}}\widecheck{HM}_{\bullet}(Y,\mathfrak{s}; \Ga_\eta).$$
We say that this grading is associated to the surface $S \subset M$. The grading defined on separate closures also induces a grading on the canonical module $\shm(M, \ga)$, as stated in Theorem \ref{thm: grading is well defined}. We write this grading on the canonical module as
$$\shm(M, \ga, S, i).$$
\edefn

\bthm[Kavi \cite{kavi2019pairing} and Li \cite{li2019direct}]\label{thm: grading is well defined}
Suppose $(M,\ga)$ is a balanced sutured manifold, and $S\subset M$ is a fixed oriented admissible properly embedded surface. Then, the grading $\shm(M,\ga,S,i)$ is independent of all the choices made in the construction and, thus, is well-defined.
\ethm

Using the grading, we can re-formulate Kronheimer and Mrowka's decomposition theorem, Proposition 6.9, in \cite{kronheimer2010knots}, as follows.
\blem\label{lem: decomposition theorem reformulated}
Suppose $(M,\ga)$ is a balanced sutured manifold and $S\subset M$ is an oriented admissible properly embedded surface. Suppose further that $S$ satisfies the hypothesis of Proposition 6.9 in \cite{kronheimer2011khovanov}, and $(M',\ga')$ is obtained from $(M,\ga)$ by a sutured manifold decomposition along $S$. Let
$$g_c=\frac{1}{4}|S\cap \ga|-\frac{1}{2}\chi(S).$$
Then, we have
$$\shm(M,\ga,S,g_c)\cong\shm(M',\ga').$$

Furthermore, the same thing holds for $\shi$.
\elem

\subsection{Bypasses}\label{subsec: by passes}
\begin{figure}[h]
\centering
\begin{overpic}[width=4.5in]{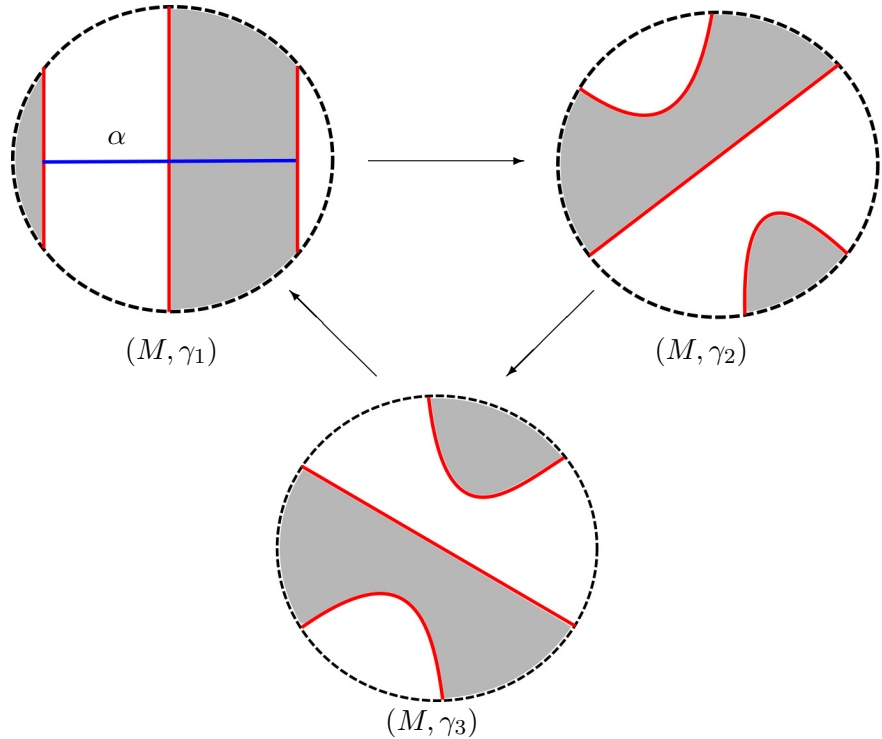}
	\put(43,-3){$(M,\ga_3)$}
	\put(13,40){$(M,\ga_1)$}
	\put(74,40){$(M,\ga_2)$}
	\put(11,65){$\al$}
	\put(41,63){\vector(1,0){18}}
	\put(67,48){\vector(-1,-1){10}}
	\put(42,38){\vector(-1,1){10}}
\end{overpic}
\vspace{0.05in}
\caption{The bypass exact triangle.}\label{fig: by pass}
\end{figure}

Suppose we have three balanced sutured manifold $(M, \ga_1)$, $(M, \ga_2)$ and $(M, \ga_3)$ so that the underlining $3$-manifolds are the same, but the sutures are different. Suppose further that $\ga_1$, $\ga_2$, and $\ga_3$ are only different with in a disk $D \subset \partial M$, and, within the disk D, they are depicted as in Figure \ref{fig: by pass}.

\bthm[Baldwin and Sivek \cite{baldwin2016contact,baldwin2018khovanov}]
\label{thm_by_pass_attachment}There are exact triangles relating the sutured monopole and instanton Floer homologies of three balanced sutured manifolds as follows.

\begin{equation*}
\xymatrix{
\shm(-M,-\ga_1)\ar[rr]^{\psi_{12}}&&\shm(-M,-\ga_2)\ar[dl]^{\psi_{23}}\\
&\shm(-M,-\ga_3)\ar[ul]^{\psi_{31}}&
}
\end{equation*}

\begin{equation*}
\xymatrix{
\shi(-M,-\ga_1)\ar[rr]^{\psi_{12}}&&\shi(-M,-\ga_2)\ar[dl]^{\psi_{23}}\\
&\shi(-M,-\ga_3)\ar[ul]^{\psi_{31}}&
}
\end{equation*}
\ethm
In contact geometry, a bypass is a half-disk, which carries some special contact structure, attached along a Legendrian arc to a convex surface. For more details, readers are referred to Honda \cite{honda2000classification}. There is a description of the maps in the above bypass exact triangle as follows. We explain how to obtain the map $\psi_{12}$, and the other two are the same. Let $Z=\partial{M}\times[0,1]$, and we can pick the suture $\ga_1$ on $\partial{M}\times\{0\}$ as well as the suture $\ga_2$ on $\partial{M}\times\{1\}$. Then, there is a special contact structure $\xi_{12}$ on $Z$, which corresponds to the bypass attachment. Hence, we can attach $Z$ to $M$ by the identification $\partial{M}\times\{0\}=\partial{M}\subset M$. The result, $(M\cup Z,\ga_2)$, is diffeomorphic to $(M,\ga_2)$, and we have
$$\psi_{12}=\Phi_{\xi_{12}}.$$
Here, $\Phi_{\xi_{12}}$ is the gluing map associated to $\xi_{12}$, as constructed by the second author in \cite{li2018gluing}.

There is a second way to interpret the maps $\psi_{\pm}$ associated to bypass attachments based on \cite{ozbagci2011contact}. In \cite{ozbagci2011contact},  Ozbagci proved that a bypass attachment could be realized by attaching a contact $1$-handle followed by a contact $2$-handle. In sutured monopole and instanton Floer homologies, there are maps associated to the contact handle attachments, due to Baldwin and Sivek \cite{baldwin2016contact,baldwin2016instanton}. So, we can compose those contact handle attaching maps to obtain $\psi_{\pm}$. This was the original way Baldwin and Sivek constructed the bypass maps (when they define bypass maps, there was no construction of gluing maps) and proved the existence of the exact triangle. The two interpretations are the same because of the functoriality of the gluing maps. For details, readers are referred to Li \cite{li2018gluing}. Both descriptions of the bypass maps are useful in later sections.

\section{Difference of supporting spin${}^c$ structures}\label{sec: difference of supporting spin c structures}
\subsection{A basic calculation}
Suppose $(M,\ga)$ is a balanced sutured manifold, and $\widetilde{M}$ is a pre-closure of $M$ with $\partial\widetilde{M}=R_+\cup R_-$. Suppose further that we pick two gluing diffeomorphisms $h_1,h_2:R_+\ra R_-$ and obtain two closures $(Y_1,R_+)$ and $(Y_2,R_+)$ of $(M,\ga)$, respectively. Let $h=h_1^{-1}\circ h_2$, and let $Y^h$ be the mapping torus of $h$. As in Subsection \ref{subsec: SHM and SHI}, we can construct a Floer excision cobordism $W$ from $Y_1\sqcup Y^h$ to $Y_2.$ Suppose $i:Y_2\ra W$ is the inclusion map. In section 4 of Li \cite{li2018gluing}, the map
$$i_*:H_1(Y_2)\ra H_1(W),$$
which is induced by the inclusion $i:Y_2\hookrightarrow W$, has played a very important role in proving Proposition \ref{prop: 1 1}. In this subsection we compute the kernel of $i_*$. A first observation is that we could just work with $\mathbb{Q}$ coefficients, since for grading purpose, torsion parts have no contributions.

From the description in Subsection \ref{subsec: SHM and SHI}, we know that $W$ is obtained by gluing three pieces together. We can compute its first homology by applying Mayer-Vietoris sequences twice and get the following result.
\begin{equation}\label{eq: first homology of W}
    H_1(W;\mathbb{Q})=[H_1(\widetilde{M};\mathbb{Q})\oplus\lgl [s_1],[s^h]\rgl]\slash [{\rm im}(h_{1,*}-1)+{\rm im}(h_*-1)].
\end{equation}
Here are a few things to be explained. First, recall in section \ref{sec: preliminaries}, we require that the gluing diffeomorphism $h_1$ and $h_2$ to fix the same base point $p\in T$. This means that there are circles $s_1\subset Y_1$ and $s^h\subset Y^h$ of the form $\{p\}\times S^1$, respectively. Then, $s_1$ and $s^h$ naturally embed into $W$, and the class $[s_1],[s^h]\in H_1(W)$ are represented by these two circles. 

Second, recall we have a map $h_1:R_+\ra R_-$, and, thus, there is a map
\begin{equation*}
\left(
\begin{array}{cc}
    -1&0\\
    (h_1)_*&0
\end{array}
\right):H_1(R_+)\oplus H_1(R_-)\ra H_1(R_+)\oplus H_1(R_-).
\end{equation*}
Note $H_1(R_+)\oplus H_1(R_-)$ can be viewed as $H_1(R_+\sqcup R_-)$ and there is an inclusion $j:R_+\sqcup R_-\ra \widetilde{M}$. So, in equation (\ref{eq: first homology of W}), we use ${\rm im}(h_{1,*}-1)$ to denote the subspace
\begin{equation*}
{\rm im}\left[j_*\circ\left(
\begin{array}{cc}
    -1&0\\
    (h_1)_*&0
\end{array}
\right)
\right]    \subset H_1(\widetilde{M}).
\end{equation*}
Note, from the maps $h:R_+\ra R_+$ and $j:R_+\sqcup R_-\ra \widetilde{M}$, we have a subspace $j_*[{\rm im}(h_*-1)]\subset H_1(\widetilde{M})$. Abusing the notation, in (\ref{eq: first homology of W}), we use ${\rm im}(h_*-1)$ to denote $j_*[{\rm im}(h_*-1)]$ and omit $j_*$ from the notation. The sum ${\rm im}(h_{1,*}-1)+{\rm im}(h_*-1)$ is the sum of the two subspaces, as described above, in $H_1(\widetilde{M})$.

In a similar way, we can compute the first homology of $Y_2$.
\begin{equation}\label{eq: first homology of Y 2}
H_1(Y_2;\mathbb{Q})=[H_1(\widetilde{M};\mathbb{Q})\oplus\lgl [s_2]\rgl]\slash[{\rm im}(h_{2,*}-1)].
\end{equation}
Here, the term ${\rm im}(h_{2,*}-1)$ is similar to the term ${\rm im}(h_{1,*}-1)$ in (\ref{eq: first homology of W}), and $s_2$ is the circle $S^1\times \{p\}\subset Y_2$, similar to $s_1$ and $s^h$ in (\ref{eq: first homology of W}). Hence, we can deduce the following. 

\blem\label{lem: kernel of i star}
Let $i:Y_2\ra W$ be the inclusion. Then, 
$$ker(i_*)\subset [{\rm im}(h_{1,*}-1)]\slash [{\rm im}(h_{2,*}-1)]\subset H_1(Y_2;\mathbb{Q}).$$
Here, the term ${\rm im}(h_{1,*}-1)$ is the same as the one appeared in (\ref{eq: first homology of W}), and the term ${\rm im}(h_{2,*}-1)$ is the same as the one appeared in (\ref{eq: first homology of Y 2}).
\elem

\bpf
It is straight forward to check that
$$i_*([s_2])=[s_1]+[s^h]\in H_1(W;\mathbb{Q}),$$
and it is clear that $[s_1]+[s^h]$ is non-zero in $H_1(W;\mathbb{Q})$ from (\ref{eq: first homology of W}).
Hence, from (\ref{eq: first homology of Y 2}), we know that the kernel must come from the quotient of $H_1(\widetilde{M};\mathbb{Q})$:
$$ker(i_*)\subset H_1(\widetilde{M};\mathbb{Q})\slash[{\rm im}(h_{2,*}-1)]\subset H_1(Y_2;\mathbb{Q}).$$

Let $\al\in ker(i_*)$ be any element in the kernel, then, from (\ref{eq: first homology of Y 2}) and (\ref{eq: first homology of W}), we can find a lift $\tilde{\al}\in H_1(\widetilde{M};\mathbb{Q})$ so that $\tilde{\al}\in {\rm im}(h_{1,*}-1)+{\rm im}(h_*-1)$. Equivalently, there are classes $\be,\ga\in H_1(R_+;\mathbb{Q})$ so that (recall $j:R_+\sqcup R_-\ra \widetilde{M}$ is the inclusion)
$$\widetilde{\al}=j_*[(h_{1})_*(\be)-\be]+j_*[h_*(\ga)-\ga].$$
Write 
$$\tilde{\al}'=h_*(\ga)-\ga\in H_1(R_+;\mathbb{Q})\subset H_1(R_+\sqcup R_-;\mathbb{Q}),$$ 
then we know that (recall $h_2=h_1\circ h$)
\beq
\tilde{\al}'&=(h_1)_*(\tilde{\al}')+[\tilde{\al}'-(h_1)_*(\tilde{\al}')]\\
&=(h_2)_*(\ga)-(h_1)_*(\ga)+[\tilde{\al}'-(h_1)_*(\tilde{\al}')]\\
&=(h_2)_*(\ga)-\ga+[\ga-(h_1)_*(\ga)]+[\tilde{\al}'-(h_1)_*(\tilde{\al}')].
\eeq

Hence, we have
\beq
\tilde{\al}=&j_*[(h_2)_*(\ga)-\ga]\\
&+j_*[\ga-(h_1)_*(\ga)]+j_*[\tilde{\al}'-(h_1)_*(\tilde{\al}')]+j_*[(h_{1})_*(\be)-\be].
\eeq
The first term is in ${\rm im}(h_{2,*}-1)$ and the rest are in ${\rm im}(h_{1,*}-1)$. So, we know that
$$\al\in[{\rm im}(h_{1,*}-1)]\slash [{\rm im}(h_{2,*}-1)]\subset H_1(\widetilde{M},\mathbb{Q})\slash [{\rm im}(h_{2,*}-1)]=H_1(Y_2;\mathbb{Q}),$$
and this concludes the proof of Lemma \ref{lem: kernel of i star}.
\epf

\subsection{Adding $1$-handles}\label{subsec: adding one handles}
\bdefn\label{defn: one handle}
Suppose $(M,\ga)$ is a balanced sutured manifold. A {\it product (or contact) $1$-handle} is a tuple $(\phi,S,D^3,\delta)$, where $S\subset \partial{D}^3$ is the disjoint union of two embedded disks on $\partial{D}^3$, $\delta$ is a simple closed curve on $\partial{D}^3$, which intersects each component of $S$ in an arc, and $\phi:S\ra \partial{M}$ is an embedding so that $\phi(\delta\cap S)= \ga\cap \phi(S)\subset\partial{M}$. Then, we can form a new balanced sutured manifold
$$(M',\ga')=(M\mathop{\cup}_{\phi}D^3,\ga'=\ga\backslash \phi(S)\cup (\delta\backslash S)).$$ 
\edefn

\brem
In Kronheimer and Mrowka \cite{kronheimer2010knots}, this process is called attaching a product $1$-handle, while in Baldwin and Sivek \cite{baldwin2016contact}, the same process is called attaching a contact $1$-handle.
\erem

\blem[Kronheimer, Mrowka \cite{kronheimer2010knots} or Baldwin, Sivek \cite{baldwin2016contact}]\label{lem: one handle does not influence pre closures}
When using auxiliary surfaces of large enough genus, any pre-closure of $(M,\ga)$ is a pre-closure of $(M',\ga')$, and vice-versa. 
\elem

Hence, by Lemma \ref{lem: one handle does not influence pre closures}, we can freely add $1$-handles to the original $(M,\ga)$ without changing its closure. A straightforward observation is the following.

\blem\label{lem: set of one handles}
For any balanced sutured manifold $(M,\ga)$, there exists a set of $1$-handles, $\{h_1,..., h_n\}$ so that the following is true.

(1) If $(M',\ga')$ is the resulting balanced sutured manifold after attaching all $1$-handles $h_1$,..., $h_n$, then $\ga'$ is connected.

(2) For $l=1,...,n$, if $(M_l,\ga_l)$ is the resulting balanced sutured manifold after attaching all $1$-handles $h_1$, ..., $h_n$ except $h_l$, then $R_{\pm}(\ga_l)$ are both connected.
\elem

\brem
The first condition is used in the proof of Lemma \ref{lem: difference of spin c structures is in M prime}, and the second is used in the proof of Lemma \ref{lem: difference of spin c structures is in M l}.
\erem

Thus, from any balanced sutured manifold $(M,\ga)$, we can find a set of $1$-handles, $\{h_1,..., h_n\}$, according to Lemma \ref{lem: set of one handles}. Let $(M',\ga')$ be the resulting balanced sutured manifold after attaching all those $1$-handles, then we have the following lemma. 
\blem\label{lem: difference of spin c structures is in M prime}
Suppose $(Y,R)$ is a closure of $(M,\ga)$ with $g(R)$ large enough. Then, $(Y,R)$ can also be regarded as a closure of $(M',\ga')$ by Lemma \ref{lem: one handle does not influence pre closures}. Suppose further that $\mathfrak{s}_1$ and $\mathfrak{s}_2$ are two supporting spin${}^c$ structures on $Y$, then there is a $1$-cycle $x$ in $M'$ so that
$$P.D.c_1(\mathfrak{s}_1)-P.D.c_1(\mathfrak{s}_2)=[x]\in H_1(Y;\mathbb{Q}).$$

Similar results hold in the instanton settings.
\elem

\bpf
We start by constructing a special reference closure of $(M',\ga')$. Pick an auxiliary surface $T'$, and let 
$$\widetilde{M}=M'\cup [-1,1]\times T'$$
be a pre-closure of $(M',\ga')$ that is also a pre-closure of $(M,\ga)$. We have
$$\partial{\widetilde{M}}=R_+\cup R_-,~R_{\pm}=R_{\pm}(\ga')\cup \{\pm 1\}\times T'.$$
We can pick a special gluing diffeomorphism $h^r:R_+\ra R_-$ so that $h^r|_{\{1\}\times T'}=id_{T'}$. Let $(Y^r,R_+)$ be the closure of $(M',\ga')$ arising from $h^r$ and $\widetilde{M}$. We know that $h^r$ can be split into two parts, $f^r=h^r|_{R_+(\ga')}$ and $id_{T'}$.

\begin{figure}[h]
\centering
    \begin{overpic}[width=5in]{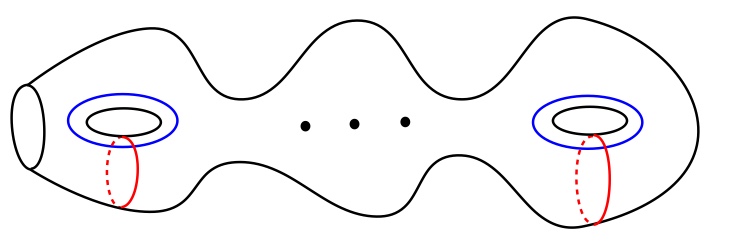}
    \put(19,9){$a_1$}
    \put(82,9){$a_k$}
    \put(24.5,16){$b_1$}
    \put(86.5,16){$b_k$}
\end{overpic}
\vspace{0.05in}
\caption{A basis for $H_1(T')$.}\label{fig: basis}
\end{figure}
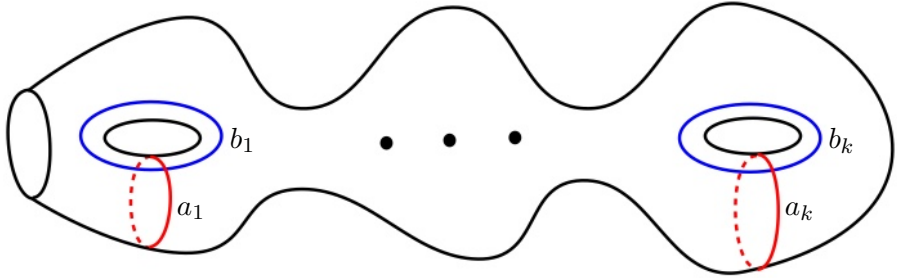

Thus, we have an alternative interpretation for $Y^r$. First, we can use $f^r:R_+(\ga')\ra R_-(\ga')$ to glue $R_+(\ga')\subset \partial{M}'$ to $R_-(\ga')\subset \partial{M}'$, and $M'$ becomes a manifold $M'_1$ with a toroidal boundary. Note $f^r|_{\partial{R_+(\ga')}}=id$, so we have a natural framing $s^r$ and $\ga'$ on $\partial{M}_1'$. Here, $s^r$ is obtained as follows. If $q\in\ga'$ is a point, then we have an arc $[-1,1]\times \{q\}\subset A(\ga)$. Note $f^r$ identifies $\{1\}\times\{q\}$ with $\{-1\}\times\{q\}$, so $[-1,1]\times\{q\}$ becomes a circle $s^r$ inside $\partial{M}_1'$. Second, we can glue $M_1'$ and $S^1\times T'$ together to form $Y^r$:
$$Y^r=M_1'\mathop{\cup}_{\phi} S^1\times T',$$
where $\phi:\partial{M}_1'\ra \partial(S^1\times T')$ maps $s^1$ to the $S^1$ direction and maps $\ga'$ to $\partial{T}'$ direction. Let $g(T)'=k$, and let $\{a_1,b_1,...,a_k,b_k\}$ be a set of generators of $H_1(T')$ as in Figure \ref{fig: basis}. Then, we can use Mayer–Vietoris sequence to conclude the following:
\begin{equation}\label{eq: first homology of Y r}
    H_1(Y^r;\mathbb{Q})= \Big(H_1(M';\mathbb{Q})\slash[{\rm im}(f^r_*-1)]\Big)\oplus\lgl [s^r],[a_1],...,[b_k]\rgl.
\end{equation}

Here, the term ${\rm im}(f^r_*-1)$ is similar to the term ${\rm im}(h_{1,*}-1)$ in (\ref{eq: first homology of W}). Suppose $\mathfrak{s}$ is a supporting spin${}^c$ structure on $Y^r$, then we can write $P.D.c_1(\mathfrak{s})$ in terms of the description of $H_1(Y^r)$ in \ref{eq: first homology of Y r}. The coefficient of the class $[s^r]$ can be understood to be $2g(R)-2$ by looking at the pairing 
$$c_1(\mathfrak{s})[R]=2g(R)-2.$$
The coefficients of $[a_1],..., [b_k]$ are all zero, since we can apply the adjunction inequality in Lemma \ref{lem: adjunction inequality in monopoles} to the tori $b_1\times S^1,..., a_k\times S^1\subset Y$. Thus, we conclude the following.
\blem\label{lem: supporting spin c structures on Y r}
Suppose $\mathfrak{s}$ is a supporting spin${}^c$ structure on $Y^r$, then there is a $1$-cycle $[x]$ in $M'$ so that
$$P.D.c_1(\mathfrak{s})=[x]+(2g-2)[s^r]\in H_1(Y^r;\mathbb{Q}).$$

Similar statement (Replacing $c_1(\mathfrak{s})$ by a supporting eigenvalue function $\lambda$, c.f. Definition \ref{defn: supporting eigenvalue function}) holds in the instanton settings.
\elem

Now suppose $h:R_+\ra R_-$ is any gluing diffeomorphism, and $(Y,R)$ is the resulting closure of $(M',\ga')$. Let $\psi=(h^r)^{-1}\circ h: R_+\ra R_+$ be a diffeomorphism, and let $Y^{\psi}$ be the mapping torus of $\psi$. We can form an excision cobordism $W$ from $Y^r\sqcup Y^{\psi}$ to $Y$, as in Subsection \ref{subsec: SHM and SHI}. Then, we can compute via Mayer–Vietoris sequences that 
\begin{equation}\label{eq: first homology of mapping torus}
    H_1(Y^{\psi};\mathbb{Q})= \Big(H_1(R_+;\mathbb{Q})\slash[{\rm im}(\psi_*-1)]\Big)\oplus\lgl [s^{\psi}]\rgl.
\end{equation}

Here, $s^{\psi}$ is the same as the term $s^h$ in (\ref{eq: first homology of W}). Let $i:Y\ra W$, $i^r:Y^r\ra W$, and $i^{\psi}:Y^{\psi}\ra W$ be the inclusions. Suppose $\mathfrak{s}$ is a supporting spin${}^c$ structure on $Y$, then, from Lemma \ref{lem: spin c structure shall extend}, we know that there is a supporting spin${}^c$ structure $\mathfrak{s}^r$ on $Y^r$ so that
$$i_*[P.D.c_1(\mathfrak{s})]=i^r_*[P.D.c_1(\mathfrak{s^r})]+i^{\psi}_*[P.D.c_1(\mathfrak{s^{\psi}})]\in H_1(W;\mathbb{Q}).$$
Here, $\mathfrak{s}^{\psi}$ is the unique spin${}^c$ structure on $Y^{\psi}$, as in Lemma \ref{lem: unique spin c structure}. To find all possible values of $P.D.c_1(\mathfrak{s})$, we first find a class $[z]\in H_1(Y;\mathbb{Q})$ so that
$$i_*([z])=i^r_*[P.D.c_1(\mathfrak{s^r})]+i^{\psi}_*[P.D.c_1(\mathfrak{s^{\psi}})]\in H_1(W;\mathbb{Q}),$$
and then
$$P.D.c_1(\mathfrak{s})\in [z]+{\rm ker}(i_*)\subset H_1(Y;\mathbb{Q}).$$
Note ${\rm ker}(i_*)$ has been understood by Lemma \ref{lem: kernel of i star}.

To find the class $[z]$, by Lemma \ref{lem: supporting spin c structures on Y r}, we know that there exists a $1$-cycle $x$ in $M'$ so that
$$P.D.c_1(\mathfrak{s}^r)=[x]+(2g-2)[s^r]\in H_1(Y^r;\mathbb{Q}).$$
For $\mathfrak{s}^{\psi}$, we can check similarly (or see Subsection 4.2 in Li \cite{li2019direct}) that there is a $1$-cycle $y$ in $R_+$ so that
$$P.D.c_1(\mathfrak{s}^{\psi})=[y]+(2g-2)[s^{\psi}]\in H_1(Y^{\psi};\mathbb{Q}).$$
Inside $W$, there are annuli $x\times[0,1]$ and $y\times[0,1]$ and a pair of pants from $s^r\sqcup s^{\psi}\subset Y^r\sqcup Y^{\psi}$ to $s\subset Y$. Hence, we can take $z$ to be the $1$-cycle in $Y$, of the form
$$z=x+y+(2g-2)s.$$
Note $x$ is in $M'$ and $y$ is in $R_+$, so there are natural ways to regard them as $1$-cycles in $Y$. Thus, we know that
$$i_*([z])=i^r_*[P.D.c_1(\mathfrak{s^r})]+i^{\psi}_*[P.D.c_1(\mathfrak{s^{\psi}})]\in H_1(W;\mathbb{Q})$$
and, hence, 
\begin{equation}\label{eq: pd of c 1 of a supporting spin c structure on Y}
P.D.c_1(\mathfrak{s})\in [x]+[y]+(2g-2)[s]+{\rm ker}(i_*)\subset H_1(Y;\mathbb{Q}).
\end{equation}

From Lemma \ref{lem: kernel of i star}, we know that
$${\rm ker}(i_*)\subset [{\rm im}(h^r_*-1)]\slash [{\rm im}(h_*-1)].$$
Note that $h^r_*-1=f^r_*-1$ since $h|_{T'\times\{1\}}=id_{T'}$. By construction, $[{\rm im}(f^r_*-1)]$ is contained in $H_1(M')$, and, thus, we know that ${\rm ker}(i_*)$ can only contribute to the part $[x]$ in (\ref{eq: pd of c 1 of a supporting spin c structure on Y}). As a result, we conclude that there is a $1$-cycle $x'$ in $M'$ so that
$$P.D.c_1(\mathfrak{s})=[x']+[y]+(2g-2)[s]\in H_1(Y;\mathbb{Q}).$$
Note $[y]+(2g-2)[s]$ is independent of the choice of the supporting spin${}^c$ structure on $Y$, so we conclude the proof of Lemma \ref{lem: difference of spin c structures is in M prime}.
\epf

\subsection{Dropping $1$-handles}
Recall we have a balanced sutured manifold $(M,\ga)$, and we have a set of $1$-handles, $\{h_1,...,h_n\}$, as in Lemma \ref{lem: set of one handles}. Recall further that $(M',\ga')$ is the resulting balanced sutured manifold after attaching all of those $1$-handles. In Lemma \ref{lem: difference of spin c structures is in M prime}, we prove that, in terms of the Poincar\'e dual of the first Chern classes, the difference of two supporting spin${}^c$ structures on $Y$ is contained in $M'$. In this subsection, we sharpen the result and prove that the difference must lie in $M$ instead of the whole $M'$, which is exactly the statement of Theorem \ref{thm: difference of supporting spin c structures}.

Suppose, for $l=1,...,n$, $(M_l,\ga_l)$ is the resulting balanced sutured manifold after attaching all $1$-handles $h_1,..., h_n$ except $h_l$. Then, $(M',\ga')$ is obtained from $(M_l,\ga_l)$ by attaching $h_l$. Recall, from Subsection \ref{subsec: adding one handles}, we have an auxiliary surface $T'$ for $M'$ and a pre-closure $\widetilde{M}$. Then, from Lemma \ref{lem: one handle does not influence pre closures}, $\widetilde{M}$ is also a pre-closure for $(M_l,\ga_l)$. Thus, any closure (Y,R) arising from $\widetilde{M}$ is also a closure for $(M_l,\ga_l)$. 
\blem\label{lem: difference of spin c structures is in M l}
For any fixed $l$, suppose $(Y,R)$ is a closure of $(M,\ga)$ with $g(R)$ large enough, then $(Y,R)$ is also a closure for $(M_l,\ga_l)$ by Lemma \ref{lem: one handle does not influence pre closures}. Suppose further that $\mathfrak{s}_1, \mathfrak{s}_2$ are two supporting spin${}^c$ structures on $Y$, then there is a $1$-cycle $x$ in $M_l$ so that
$$P.D.c_1(\mathfrak{s}_1)-P.D.c_1(\mathfrak{s}_2)=[x]\in H_1(Y;\mathbb{Q}).$$

A similar result holds in the instanton settings.
\elem

\bpf[Proof of Theorem \ref{thm: difference of supporting spin c structures}.]
Note $(M',\ga')$ and $(M_l,\ga_l)$ are obtained from $(M,\ga)$ by attaching $1$-handles. So, there are injections
$$H_1(M;\mathbb{Q})\hookrightarrow H_1(M_l;\mathbb{Q})\hookrightarrow H_1(M';\mathbb{Q}).$$
Also, inside $H_1(M';\mathbb{Q})$, we have
$$H_1(M)=\bigcap_{l=1}^n H_1(M_l).$$
So, Lemma \ref{lem: difference of spin c structures is in M l}, together with Lemma \ref{lem: one handle does not influence pre closures}, implies Theorem \ref{thm: difference of supporting spin c structures}.
\epf

\bpf[Proof of Lemma \ref{lem: difference of spin c structures is in M l}.]
By Lemma \ref{lem: difference of spin c structures is in M prime}, we know that there is a $1$-cycle $x$ is in $M'$ which satisfies the statement of the lemma, but our goal is to show that $x$ can be chosen inside $M_{l}$. Recall, in the proof of Lemma \ref{lem: difference of spin c structures is in M prime}, we pick a special reference closure $Y^r$ by requiring $h^r|_{T'\times\{1\}}=id_{T'}$ (and $f^r=h^r|_{R_+(\ga')}$). To prove the current lemma, we will need an even more special closure, by making further restrictions on $f^r$. Though we will keep using the notations $Y^r$, $h^r$ and $f^r$, etc.

\begin{figure}[h]
\centering
    \begin{overpic}[width=4.5in]{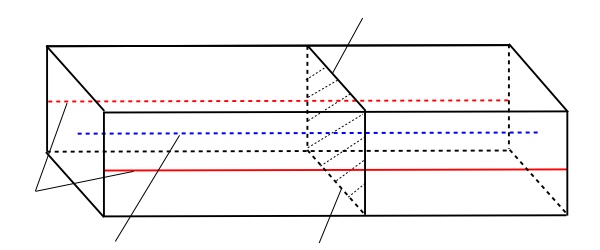}
    \put(51,-1){$\be_-$}
    \put(60,40){$\be_+$}
    \put(62,8){$D_l$}
    \put(17,-1){$\al_l$}
    \put(4,6.5){$\ga'$}
\end{overpic}
\vspace{0.05in}
\caption{The one handle $h_l$.}\label{fig: one handle}
\end{figure}
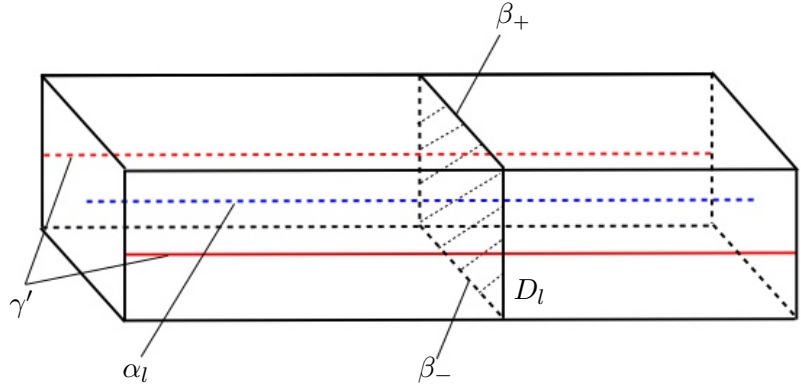

To explain the further restriction on $f^r$, recall that $(M',\ga')$ is obtained from $(M_l,\ga_l)$ by attaching the $1$-handle $h_l$. We can write
$$H_1(M';\mathbb{Q})=H_1(M_l;\mathbb{Q})\oplus\lgl \al_l\rgl,$$
where $\al_l$ consists of the core of the $1$-handle $h_l$ together with an arc inside $M_l$. Let $D_l$ be the co-core of the $1$-handle $h_l$. It is a properly embedded disk $D_l\subset M'$ that intersects $\al_l$ transversely at one point, and $\partial{D}_l$ intersects $R_{\pm}(\ga')$ in arcs $\be_{l,\pm}$. See Figure \ref{fig: one handle}. Note, by the condition (2) in Lemma \ref{lem: set of one handles}, the arcs $\be_{l,\pm}$ are non-separating inside $R_{\pm}(\ga')$. Thus, from Lemma 3.6 in Li \cite{li2018gluing}, we can find a map $f^r$ which sends $\be_{l,+}$ to $\be_{,-}$. 

Now we can use the new $h^r$ to obtain the closure $(Y^r,R_+)$. It not only satisfies all requirements needed to conclude the proof of Lemma \ref{lem: difference of spin c structures is in M prime}, but also has some new features. Suppose the co-core $D_l$ intersects $\ga'$ at two points $p_l,q_l$ and let $\eta\subset T'$ be a simple arc with end points $p_l,q_l$. Then, we know that 
$$\widetilde{D}_{l}=D_l\cup(\eta\times[-1,1])$$
is a properly embedded surface in $\widetilde{M}$ so that
$$\partial{\widetilde{D}}_l\cap R_{\pm}=C_{\pm}=\be_{l,\pm}\cup\eta\times\{\pm1\}.$$
From the above construction, we know that $h^r(C_+)=C_-$. Hence, $\widetilde{D}_l$ becomes a torus $T_l\subset Y^r$. (It is straightforward to compute the Euler characteristic to see that it is indeed a torus.) Note this torus $T_l$ intersects $\al_l$ transversely at one point, but is disjoint from all other generators of
$$H_1(M';\,\mathbb{Q})=H_1(M_l;\,\mathbb{Q})\oplus \lgl[\al_l]\rgl.$$

As a result, we can use the adjunction inequality in Lemma \ref{lem: adjunction inequality in monopoles} to make refinement of Lemma \ref{lem: supporting spin c structures on Y r} as follows.

\blem\label{lem: supporting spin c structures on Y r refined}
Suppose $\mathfrak{s}$ is a supporting spin${}^c$ structure on $Y^r$, then there is a $1$-cycle $x$ in $M_l$ so that
$$P.D.c_1(\mathfrak{s})=[x]+(2g-2)[s^r]\in H_1(Y^r;\mathbb{Q}).$$

Similar results hold in the instanton settings.
\elem

Suppose $(Y,R_+)$ is an arbitrary closure of $(M,\ga)$, arising from the pre-closure $\widetilde{M}$ and a gluing diffeomorphism $h$. We can form the diffeomorphism $\psi=(h^r)^{-1}\circ h$, the mapping torus $Y^{\psi}$, and the excision cobordism $W$, as in the proof of Lemma \ref{lem: difference of spin c structures is in M prime}. Suppose further that $\mathfrak{s}$ is a supporting spin${}^c$ structure on $Y$, then we know from (\ref{eq: pd of c 1 of a supporting spin c structure on Y}) that
$$P.D.c_1(\mathfrak{s})\in [x]+[y]+(2g-2)[s]+{\rm ker}(i_*)\subset H_1(Y;\mathbb{Q}),$$
where $i:Y\ra W$ is the inclusion, $y$ is a $1$-cycle on $R_+$, $s$ is the curve $\{p\}\times S^1\subset Y$, and $x$ is a $1$-cycle in $M_l$ guaranteed by Lemma \ref{lem: supporting spin c structures on Y r refined}. 

From Lemma \ref{lem: kernel of i star}, we know that
$${\rm ker}(i_*)\subset [{\rm im}(h_{*}^r-1)]\slash [{\rm im}(h_{2,*}-1)].$$
By construction, we know that
$${\rm im}(h_{*}^r-1)={\rm im}(f_{*}^r-1)\subset H_1(M';\mathbb{Q}).$$ 
To conclude the proof of Lemma \ref{lem: difference of spin c structures is in M l}, we need to show that 
$${\rm im}(h_{*}^r-1)={\rm im}(f_{*}^r-1)\subset H_1(M_l;\mathbb{Q}),$$
and this is equivalent to show that, under the decomposition
$$H_1(M';\mathbb{Q})=H_1(M_l;\mathbb{Q})\oplus\lgl[\al_l]\rgl,$$
any element $[z]\in {\rm im}(f_{*}^r-1)$ can not have a non-zero $[\al_l]$ component. 

To prove this final statement, suppose $[z]\in {\rm im}(f_{*}^r-1)$ is of the form 
$$[z]=a\cdot [\al_l]+H_1(M_l;\mathbb{Q})$$
for some $a\in\mathbb{Q}$, then we need to show that $a=0$. Note, from (\ref{eq: first homology of Y r}), we know that $[z]=0\in H_1(Y^r)$. Also, inside $Y^r$, $[z]\cdot [T^2_l]=a$, since, by construction, $[\al_l]\cdot [T^2_l]=1$, and $T^2_l\cap M_l=\emptyset$. Thus, we know that $a=0$, and this concludes the proof of lemma \ref{lem: difference of spin c structures is in M l}. As mentioned above, this also concludes the proof of Theorem \ref{thm: difference of supporting spin c structures}. 
\epf

\section{General grading shifting formula}\label{sec: general grading shifting property}
In this section, we prove the generalized grading shifting formula, Theorem \ref{thm: general grading shifting property}, as stated in the introduction. We also use it to compute the sutured monopole and instanton Floer homologies of some particular sutured handlebodies.

\bprop\label{prop: general grading shifting property}
Suppose $(M,\ga)$ is a balanced sutured manifold, and $S\subset M$ is a properly embedded surface. Pick $i\in\intg$ so that the surface $S^i$, which is obtained from $S$ by performing $i$ stabilizations, is admissible (see Definition \ref{defn: admissibility of surfaces}). Pick any $k\in\intg$. Then, there exist constants $l_M,l_I\in \intg$ so that, for any $j\in\intg$, we have:
$$\shm(-M,-\ga,S^{i+2k},j)=\shm(-M,-\ga,S^{i},j-l_M),$$
and
$$\shi(-M,-\ga,S^{i+2k},j)=\shi(-M,-\ga,S^{i},j-l_I).$$
Moreover, if the sutured manifold decompositions of $(M,\ga)$ along $S$ and $-S$ are both taut, then $l_M=l_I=k$. 
\eprop

\begin{figure}[h]
\centering
    \begin{overpic}[width=5.5in]{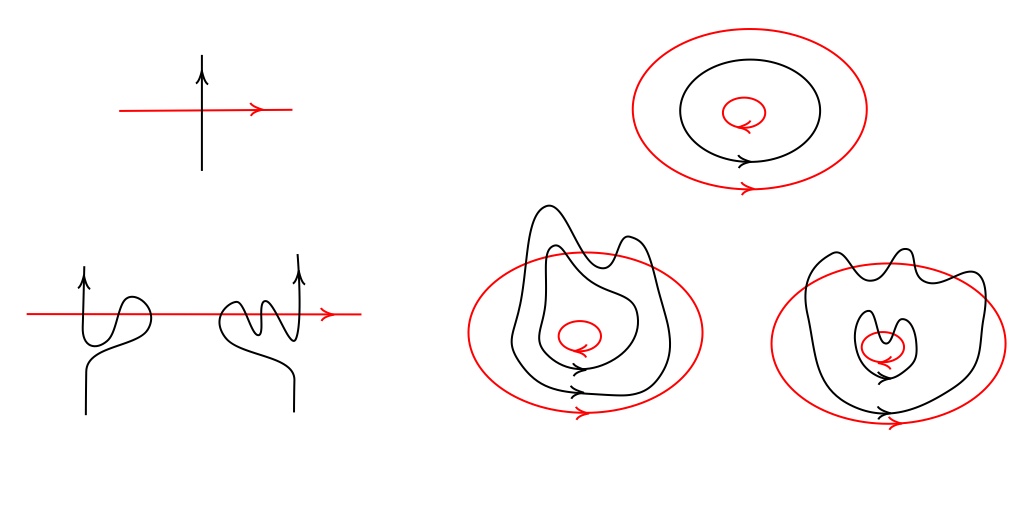}
    \put(19,46){$S$}
    \put(29,39){$\ga$}
    \put(7,6){$S^i$}
    \put(26,6){$S^{i+2k}$}
    \put(72,38.5){$\ga$}
    \put(82,38.5){$\ga$}
    \put(77.5,38.5){$S$}
    \put(59,17){$S^i$}
    \put(59,27.5){$S^{i+2k}$}
    \put(85,20){$S^{i+2k}$}
    \put(87,26.5){$S^{i}$}
\end{overpic}
\vspace{-0.3in}
\caption{Pushing off $S^{i+2k}$.}\label{fig: push off}
\end{figure}
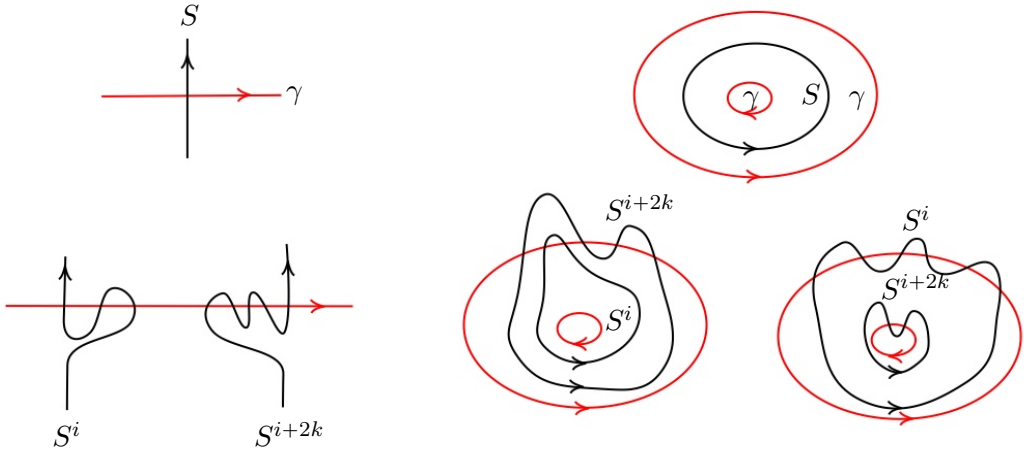

\bpf
We now prove the first part of the proposition, namely the existence of $l$ and the fact that it is independent of $j$. As usual, we argue this in the monopole settings, but the same is true for the instanton settings. Suppose $(M,\ga)$ and $S$ are defined the same as in the hypothesis of the proposition. We follow the idea in Subsection \ref{subsec: construction of gradings}, to construct a closure $(Y,R)$ of $(M,\ga)$ so that both $S^i$ and $S^{i+2k}$ extend to closed surfaces. To do this, put both $S^i$ and $S^{i+2k}$ in $M$ so that they are transverse and $\partial{S}^i\cap \partial S^{i+2k}=\emptyset$ in $Y$. Note we can always achieve this. If $\partial{S}$ intersects $\ga$, then the positive and negative stabilizations can be performed in an arbitrarily small neighborhood of $S\cap\ga$, so we can simply start with two parallel copies of $S$ and perform $i$ and $i+2k$ stabilizations respectively. If $\partial{S}\cap\ga=\emptyset$, then positive and negative stabilizations happen on different sides of $S$, and we can always perturb them to be distinct. See figure \ref{fig: push off}.

Now we pick a connected auxiliary surface $T$ for $(M,\ga)$, which is of large enough genus. For both $S^i$ and $S^{i+2k}$, we can apply the construction of  gradings as in Subsection \ref{subsec: construction of gradings}. If the genus of $T$ is chosen to be large enough, then we could arrange the arcs, which come from both $\partial{S}^i$ and $\partial{S}^{i+2k}$, all represent linearly independent classes in $H_1(T,\partial{T})$. We can then form the pre-closure $\widetilde{M}=M\cup [-1,1]\times T$ and it has two boundary components $\partial{\widetilde{M}}=R_+\cup R_-$. As in Subsection \ref{subsec: construction of gradings}, we know that $S^i$ and $S^{i+2k}$ both extend to  properly embedded surfaces $\widetilde{S}^{i}, \widetilde{S}^{i+2k}\subset \widetilde{M}$, and there are equal number of boundary components on $R_+$ and on $R_-$.

Thus, we can pick an orientation preserving diffeomorphism $h:R_+\ra R_-$ so that
$$h(\widetilde{S}^i\cap R_+)=\widetilde{S}^{i}\cap R_-,~{\rm and}~h(\widetilde{S}^{i+2k}\cap R_+)=\widetilde{S}^{i+2k}\cap R_-.$$
Using $h$ and $\widetilde{M}$, we obtain a closure $(Y,R_+)$ of $(-M,-\ga)$ so that there are closed surfaces $\widebar{S}^{i}\subset Y$ and $\widebar{S}^{i+2k}\subset Y$. Pick a non-separating simple closed curve $\eta$ which is disjoint from $\widebar{S}^i\cup \widebar{S}^{i+2k}$. Then, as in Subsection \ref{subsec: construction of gradings}, we know that
$$SHM(-M,-\ga,S^{i},j)=\bigoplus_{\substack{\mathfrak{s}\in\mathfrak{S}^*(Y|R_+)\\c_{1}(\mathfrak{s})[\widebar{S}^{i}]=2j}}\widecheck{HM}_{\bullet}(Y,\mathfrak{s};\Ga_{\eta}),$$
and
$$SHM(-M,-\ga,S^{i+2k},j)=\bigoplus_{\substack{\mathfrak{s}\in\mathfrak{S}^*(Y|R_+)\\c_{1}(\mathfrak{s})[\widebar{S}^{i+2k}]=2j}}\widecheck{HM}_{\bullet}(Y,\mathfrak{s};\Ga_{\eta}).$$

Now, suppose $\mathfrak{s}_1,\mathfrak{s}_2\in\mathfrak{S}^*(Y|R_+)$ are two supporting spin${}^c$ structures, then, from Theorem \ref{thm: difference of supporting spin c structures}, we know that there is a $1$-cycle $x\subset M$ so that
$$[c_1(\mathfrak{s}_1)-c_1(\mathfrak{s}_2)]=P.D.[x]\in H^2(Y).$$
Since $x\subset M$ and $S^i$ is isotopic to $S^{i+2k}$ in $M$, we know that
$$[c_1(\mathfrak{s}_1)-c_1(\mathfrak{s}_2)][\widebar{S}^{i}]=[c_1(\mathfrak{s}_1)-c_1(\mathfrak{s}_2)][\widebar{S}^{i+2k}].$$
Thus, the number
$$l_M=\frac{1}{2}(c_1(\mathfrak{s}_1)[\widebar{S}^{i}]-c_1(\mathfrak{s}_1)[\widebar{S}^{i+2k}])=\frac{1}{2}(c_1(\mathfrak{s}_2)[\widebar{S}^{i}]-c_1(\mathfrak{s}_2)[\widebar{S}^{i+2k}])$$
is the desired constant in the statement of the proposition.

When the decomposition of $(M,\ga)$ along $S$ and $-S$ are both taut, then we can settle down the value of $l$ by looking at the top or bottom non-vanishing grading and conclude $l=k$. This part of the proof is exactly the same as the proof of Proposition 4.3 in Li \cite{li2019direct}.
\epf

\begin{conj}
In general, we always have $l=k$.
\end{conj}

\bcor\label{cor: disk with four intersection with the suture}
Suppose $(M,\ga)$ is a taut balanced sutured manifold and $D\subset M$ is a disk so that $D$ intersects $\ga$ transversely four times. Suppose further that the decompositions of $(M,\ga)$ along $D$ and $-D$ are $(M',\ga')$ and $(M'',\ga'')$, respectively. If at least one of the two decompositions is taut, then
$$\shm(M,\ga)\cong\shm(M',\ga')\oplus\shm(M'',\ga'').$$
Similar results holds for the instanton settings.
\ecor

\bpf
Without loss of generality, we can assume that $(M',\ga')$ is taut. We can perform a positive stabilization on $D$ to make it admissible. Then, it induces a grading
$$\shm(M,\ga,D^+,i).$$
$D^+$ is a disk intersecting the suture six times. From the construction of grading in Section \ref{subsec: construction of gradings}, it becomes a genus-two closed surface inside a suitable closure of $(M,\ga)$. The adjunction inequality in Lemma \ref{lem: adjunction inequality in monopoles} then tells us that there are only three potentially non-trivial gradings, being $i=-1,0,1$. Lemma \ref{lem: decomposition and pm stabilization} and Lemma \ref{lem: decomposition theorem reformulated} then imply that
$$\shm(M,\ga,D^+,1)\cong\shm(M',\ga'),$$ 
$$\shm(M,\ga,D^+,-1)=0.$$
Applying Proposition \ref{prop: general grading shifting property} and Lemma \ref{lem: decomposition theorem reformulated}, we know that
\beq
\shm(M,\ga,D^+,0)&=\shm(M,\ga,D^-,-1)\\
&=\shm(M,\ga,(-D)^+,1)\\
&\cong\shm(M'',\ga'')
\eeq
Hence, we are done.
\epf

This gives an affirmative answer to Conjecture 4.3 in Li \cite{li2018contact}.

\bcor
Suppose $V$ is a solid torus, and $\ga^4$ consists of four longitudes. When using $\intg$ coefficients, we have
$$SHM(V,\ga^4)\cong \intg^2.$$
\ecor
\bpf
We apply Corollary \ref{cor: disk with four intersection with the suture}: There is a meridian disk $D$ intersecting the suture $\ga^4$ four times. The decomposition of $(V,\ga)$ along $D$ and $-D$ both result in a $3$-ball with a connected suture on the boundary, whose Floer homology is simply $\intg$.
\epf

With the help of Proposition \ref{prop: general grading shifting property}, we can prove the general grading shifting property, Theorem \ref{thm: general grading shifting property}, as stated in the introduction.

\bpf[Proof of Theorem \ref{thm: general grading shifting property}]
We first sketch the proof of the theorem as follows: we first find a closure $(Y,R)$ of $(M,\ga)$ so that, after some suitable isotopies, $S_1$ and $S_2$ both extend to some closed surfaces. Next, when we want to compare the gradings associated to surfaces $S_1$ and $S_2$, we need only to compare the difference of the first Chern classes of different supporting spin${}^c$ structures evaluated on the fundamental classes of the extensions of the two surfaces $S_1$ and $S_2$. By Proposition \ref{prop: 1 1}, the difference of supporting spin${}^c$ structures are due to a $1$-cycle inside the sutured manifold $M$. Since $S_1$ and $S_2$ represent the same relative homology class in $M$, we could conclude the result of an overall grading shift.

Now assume that $S_1$ and $S_2$ are transverse to each other. We need to isotope $S_1$ and $S_2$ into $S_1^{i_1,i_1'}$ and $S_2^{i_2,i_2'}$, respectively. Here, $i_1, i_2\geq0$ indicate the number of positive stabilizations on $S_1$ and $S_2$, respectively. Similarly, $i_1', i_2'\leq0$ correspond to the negative stabilizations. We require the following six conditions to hold.

(1) Both $S_1^{i_1,i_1'}$ and $S_2^{i_2,i_2'}$ are admissible, and no more intersection points are created during the stabilizations.

(2) Any positive intersection of $S_1^{i_1,i_1'}$ with $S_2^{i_2,i_2'}$ is contained in $R_+(\ga)$, and any negative intersection is contained in $R_-(\ga)$.

(3) If $\theta_1$ is a component of $\partial S_1^{i_1,i_1'}$ and $\theta_1\cap S_2^{i_2,i_2'}\neq\emptyset$, then $\theta_1\cap\ga\neq\emptyset$.

(4) If $\theta_2$ is a component of $\partial S_2^{i_2,i_2'}$ and $\theta_2\cap S_1^{i_1,i_1'}\neq\emptyset$, then $\theta_2\cap\ga\neq\emptyset$.

(5). If $\delta_1$ is a component of $\partial S_1^{i_1,i_1'}\cap R(\ga)$, then $\delta$ intersects $S_2^{i_2,i_2'}$ at most once.

(6). If $\delta_2$ is a component of $\partial S_2^{i_2,i_2'}\cap R(\ga)$, then $\delta$ intersects $S_1^{i_1,i_1'}$ at most once.

Pick a connected auxiliary surface $T$ of large enough genus, and form the pre-closure $\widetilde{M}=M\cup [-1,1]\times T$. It has two boundary components: 
$$\partial{\widetilde{M}}=R_+\cup R_-.$$ 
We can carry out the construction of gradings in Subsection \ref{subsec: construction of gradings} on both $S_1^{i_1,i_1'}$ and $S_2^{i_2,i_2'}$. Suppose
$$n_1=\frac{1}{2}|S_1^{i_1,i_1'}\cap \ga|,~{\rm and}~n_2=\frac{1}{2}|S_2^{i_2,i_2'}\cap\ga|.$$
Pick two balanced pairings $\mathcal P_1$ and $\mathcal{P}_2$, as introduced in Subsection \ref{subsec: construction of gradings}, for the two surfaces $S_1^{i_1,i_1'}$ and $S_2^{i_2,i_2'}$, respectively. Inside $T$, pick a set of disjoint properly embedded arcs $\{\al_1,...,\al_{n_1},\be_1,...,\be_{n_2}\}$ so that the end points of $\al_k$ are identified with the intersection points of $S_1^{i_1,i_1'}\cap\ga$, according to the balanced pairing $\mathcal P_1$, and the end points of $\be_j$ are identified with the intersection points of $S_2^{i_2,i_2'}\cap\ga$, according to the balanced pairing $\mathcal P_2$. There is one special requirement for $\mathcal{P}_1$: 

(7) If $\delta_1$ is a component of $\partial S_1^{i_1,i_1'}\cap R(\ga)$, which intersects $S_2^{i_2,i_2'}$ non-trivially, then there is an arc $\al_{k_0}$ so that $\partial{\al}_{k_0}$ is identified with $\partial{\delta}_1$ 

Strictly speaking, $\partial{\delta}\subset\partial{R(\ga)}$, but we could regard $\partial{\delta}_1$ to be on $\ga$, since $\partial{R}(\ga)$ is parallel to $\ga$. Similarly, we require the following for $\mathcal{P}_2$.

(8) If $\delta_2$ is a component of $\partial S_2^{i_2,i_2'}\cap R(\ga)$, which intersects $S_1^{i_1,i_1'}$ non-trivially, then there is an arc $\be_{j_0}$ so that $\partial{\be}_{j_0}$ is identified with $\partial{\delta}_2$.

Note when we perform enough positive and negative stabilizations, the balanced pairing satisfying the constraints (7) and (8) always exist. When the genus of $T$ is large enough, we can choose the arcs, $\al_1,...,\al_{n_1},\be_1,...,\be{n_2}$, to represent linearly independent classes in $H_1(T,\partial{T})$. Then, inside $\widetilde{M}$, $S_1^{i_1,i_1'}$ and $S_2^{i_2,i_2'}$ extend to properly embedded surfaces $\widetilde{S}_1^{i_1,i_1'}$ and $\widetilde{S}_2^{i_2,i_2'}$, respectively. From the construction, we know that $\partial \widetilde{S}_1^{i_1,i_1'}\cap R_+$ and $\partial \widetilde{S}_1^{i_1,i_1'}\cap R_-$ have equal number of boundary components. Thus, let
$$\partial \widetilde{S}_1^{i_1,i_1'}\cap R_+=C_{+,1}\cup...\cup C_{+,s},~{\rm and}~ \partial \widetilde{S}_1^{i_1,i_1'}\cap R_-=C_{-,1}\cup...\cup C_{-,s}.$$
Similarly, we can assume
$$\partial \widetilde{S}_2^{i_2,i_2'}\cap R_+=D_{+,1}\cup...\cup D_{+,t},~{\rm and}~ \partial \widetilde{S}_2^{i_2,i_2'}\cap R_-=D_{-,1}\cup...\cup D_{-,t}.$$

Note the intersection points of $\widetilde{S}_1^{i_1,i_1'}$ and $\widetilde{S}_2^{i_2,i_2'}$ are in one-to-one-correspondence to the intersection points of $S_1$ and $S_2$ by requirement (1). We claim that $\partial S_1\cap \partial S_2$ consists of an even number of positive and negative points. Indeed, $\partial{S}_i\cap \partial{S}_j=\partial{S}_i\cap S_j$, and it is clear that the algebraic intersection number of $\partial{S}_i$ and $S_j$ is zero. Hence, on $R_{\pm}$, we have a collection of circles $C_{\pm,1},...,C_{\pm,s},D_{\pm,1}...D_{\pm,t}$. They represent linearly independent classes in $H_1(R_+)$. There might be intersections between $C_{\pm,k}$ with $D_{\pm,j}$, but, by requirement (5), (6), (7), and (8), each $C_{\pm,k}$ intersects with at most one $D_{\pm,j}$,  and each $D_{\pm,j}$ intersects with at most one $C_{\pm,k}$. Hence, the pattern of $C_{+,1},...,C_{+,s},D_{+,1},...,D_{+,t}$ on $R_+$ is exactly the same as the pattern of $C_{-,1},...,C_{-,s},D_{-,1},...,D_{-,t}$ on $R_-$. As a result, there exists an orientation preserving diffeomorphism $h:R_+\ra R_-$ so that
$$h(\partial \widetilde{S}_1^{i_1,i_1'}\cap R_+)=\partial \widetilde{S}_1^{i_1,i_1'}\cap R_-,~{\rm and}~h(\partial \widetilde{S}_2^{i_2,i_2'}\cap R_+)=\partial \widetilde{S}_2^{i_2,i_2'}\cap R_-.$$

Hence, we can obtain a closure $(Y,R_+)$ of $(M,\ga)$ from $\widetilde{M}$ and $h$. Inside $Y$ there are two closed surfaces $\bar{S}_1^{i_1,i_1'}$ and $\bar{S}_n^{i_2,i_2'}$, and they induce gradings on $\shm(M,\ga)$ that are associated to $S_1^{i_1,i_1'}$ and $S_2^{i_2,i_2'}$, respectively. The rest of the proof is then exactly the same as the proof of Proposition \ref{prop: general grading shifting property}.
\epf

\begin{figure}[h]
\centering
\begin{overpic}[width=6.0in]{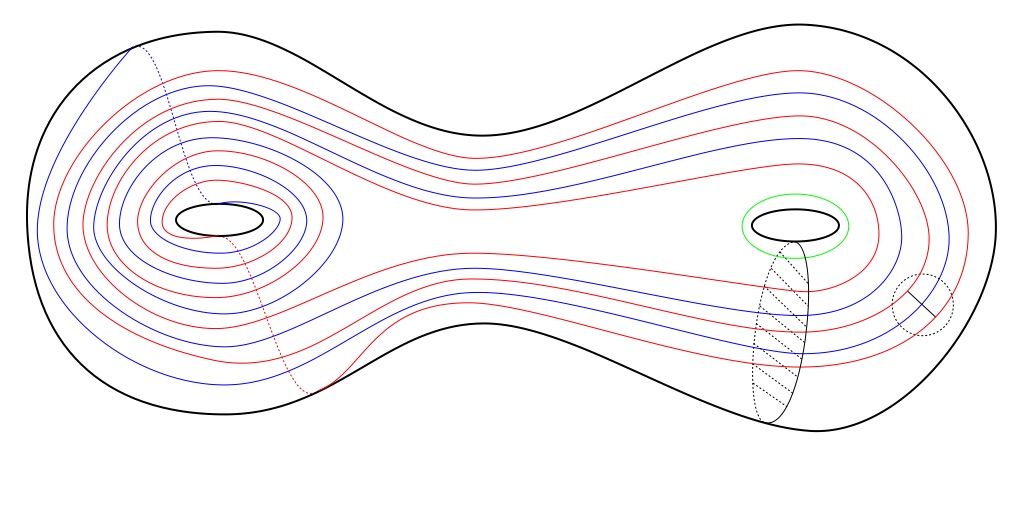}
    \put(78,10){$D$}
    \put(90,19){\line(1,-1){7}}
    \put(93,9){By-pass}
\end{overpic}
\vspace{0.05in}
\caption{The sutured manifold $(M,\ga)$. The three curves (red, blue and green) are the sutures. The disk $D$ is used to construct a grading.}\label{fig: genus two handle body}
\end{figure}

Next, we present one example to demonstrate the usage of the techniques developed above. The simplest sutured manifolds having non-trivial Floer homology are sutured handlebodies. All genus-one sutured handlebodies have been dealt with in the second author's previous paper \cite{li2019direct}. So we will work with a genus-two handle body. In particular, we want to compute the sutured monopole and instanton Floer homologies of the following balanced sutured manifold $(M,\ga)$, as depicted in Figure \ref{fig: genus two handle body}. (Strictly speaking, they are the sutured monopole and instanton Floer homologies of $(-M,-\ga)$.)

The idea is to apply the by-pass exact triangles repeatedly. There is a graded version of by-pass exact triangles, as in Li \cite{li2019direct}, generalizing the by-pass exact triangle introduced in Subsection \ref{subsec: by passes}.

\begin{equation}\label{eq: graded by-pass exact triangle}
    \xymatrix{
    \shm(-M,-\gamma,D,i)\ar[r]&\shm(-M,-\ga_1,D_1^{-2},i)\ar[d]\\
    &\shm(-M,-\ga_2,D_2^{+1},i)\ar[ul].
    }
\end{equation}
Here, the surface $D$ is chosen as in Figure \ref{fig: genus two handle body} so that it has six transverse intersections with the suture $\ga$. The surfaces $D_1, D_2\subset M$ are isotopic to $D$, but having minimal possible transverse intersection with $\ga_1$ (two intersections) and $\ga_2$ (four intersections) respectively. The superscripts in $D_2^{+1}$ and $D_1^{-2}$ imply the number of positive or negative stabilizations performed on the surfaces $D_1\subset M$ and $D_2\subset M$, as introduced in Subsection \ref{subsec: construction of gradings}. A direct check shows that $(M,\ga_1)$ is a product sutured manifold, and the suture $\ga_2$ is depicted as in Figure \ref{fig: genus two handle body, 2}.

\begin{figure}[h]
\centering
\begin{overpic}[width=6.0in]{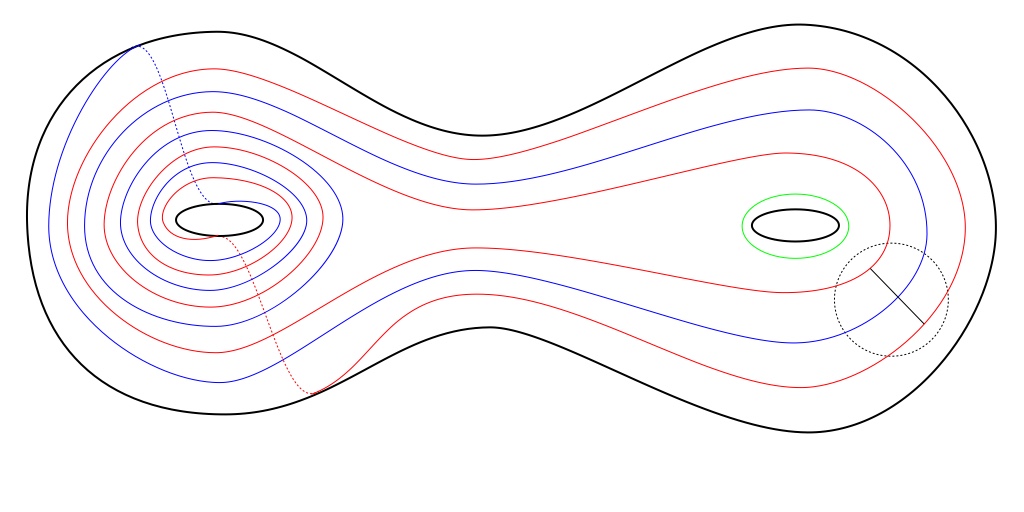}
    \put(86,18){\line(1,-1){9}}
    \put(91,6){By-pass}
\end{overpic}
\vspace{0.05in}
\caption{The sutured manifold $(M,\ga_2)$. The three curves (red, blue and green) are the sutures.}\label{fig: genus two handle body, 2}
\end{figure}

Now we know that
$$\shm(-M,-\ga_1)=\shm(-M,-\ga_1,D_1,0)=\mathcal{R}.$$
As a result, by Lemma \ref{lem: decomposition and pm stabilization} and Lemma \ref{lem: decomposition theorem reformulated}, we know that
\begin{equation}\label{eq: SHM of (M,ga_1)}
    \shm(-M,-\ga_1,D_1^{-2},i)=\left\{
    \begin{array}{cc}
        \mathcal{R}&i=1\\
        0& {\rm others}
    \end{array}
    \right.
\end{equation}
Note $D_1$ is a disk intersecting the suture six times. So it becomes a genus-two closed surface inside some suitable closure as in Subsection \ref{subsec: construction of gradings}. From the adjunction inequality in Lemma \ref{lem: adjunction inequality in monopoles}, we know that
$$\shm(-M,-\ga_2,D_2^{+1},i)=0$$
for $i>1$ or $i<-1$. From Lemma \ref{lem: decomposition and pm stabilization} and Lemma \ref{lem: decomposition theorem reformulated}, we know that
$$\shm(-M,-\ga_2,D_2^{+1},1)=0.$$
Note from above discussions, we know that for each $i$, either
$$\shm(-M,-\ga_1,D_1^{-2},i)=0$$
or
$$\shm(-M,-\ga_2,D_2^{+1},i).$$
So from the graded exact triangle in (\ref{eq: graded by-pass exact triangle}), we conclude that
\begin{equation}\label{eq: SHM of (M,ga)}
    \shm(-M,-\ga)=\mathcal{R}\oplus\shm(-M,-\ga_2).
\end{equation}

To compute the sutured monopole Floer homology of $(M,\ga_2)$, we perform the same trick once more (or equivalently apply Corollary \ref{cor: disk with four intersection with the suture}) and conclude that
\begin{equation}\label{eq: SHM of (M,ga_2)}
    \shm(-M,-\ga_2)=\mathcal{R}\oplus\shm(-M,-\ga_3),
\end{equation}
where $(M,\ga_3)$ is the balanced sutured manifold as depicted in Figure \ref{fig: genus two handle body, 4}.

\begin{figure}[h]
\centering
\begin{overpic}[width=6.0in]{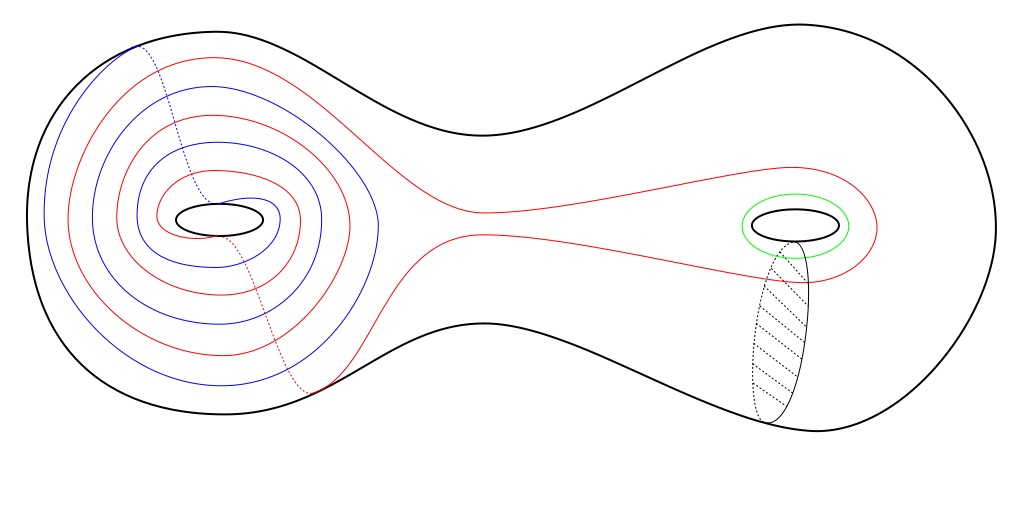}
    \put(78,11){$D$}
\end{overpic}
\vspace{0.05in}
\caption{The sutured manifold $(M,\ga_3)$. The three curves (red, blue and green) are the sutures.}\label{fig: genus two handle body, 4}
\end{figure}

To compute the sutured monopole Floer homology of $(M,\ga_3)$, we could perform a sutured manifold decomposition along the disk $D$ as depicted in Figure \ref{fig: genus two handle body, 4}. Suppose $(M_4,\ga_4)$ is the resulting balanced sutured manifold, then, from Kronheimer and Mrowka \cite{kronheimer2010knots}, we know that 
$$\shm(-M,-\ga_3)\cong\shm(-M_4,-\ga_4).$$
Furthermore, the balanced sutured manifold $(M_4,\ga_4)$ is a solid torus equipped with two curves of slope $\frac{1}{3}$ as the suture, as depicted in Figure \ref{fig: genus two handle body, 5}, so from Li \cite[Proposition 1.3]{li2019direct}, we know that
$$\shm(-M,-\ga_3)\cong\shm(-M_4,-\ga_4)\cong\mathcal{R}^3.$$
Finally, we conclude that
\bprop
If $(M,\ga)$ is the sutured manifold as shown in Figure \ref{fig: genus two handle body}, then
$$\shm(-M,-\ga)\cong\mathcal{R}^5.$$
\eprop

\begin{figure}[h]
\centering
\begin{overpic}[width=6.0in]{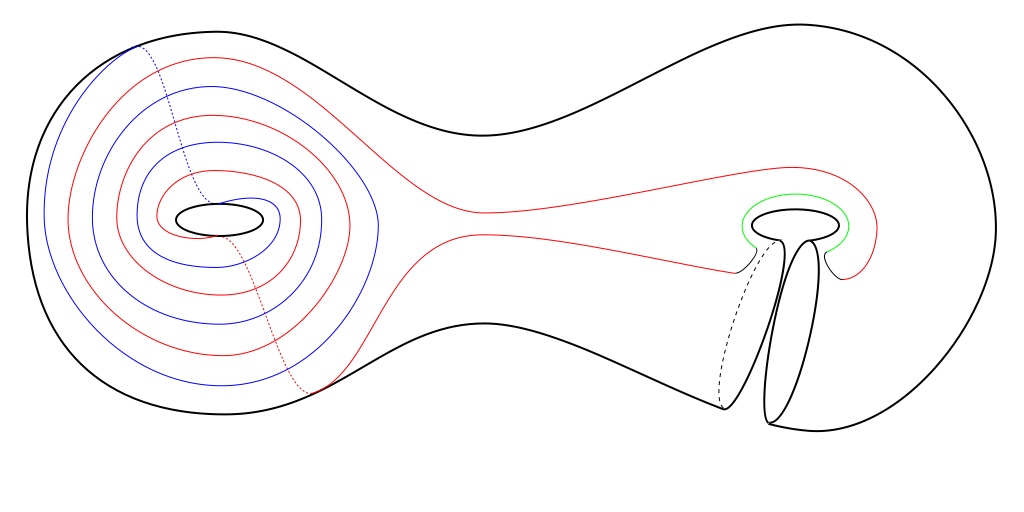}
\end{overpic}
\vspace{0.05in}
\caption{The sutured manifold $(M,\ga_4)$. The three curves (red, blue and green) are the sutures.}\label{fig: genus two handle body, 5}
\end{figure}

The same type of arguments in instanton theory yields the following.
\bcor
If $(M,\ga)$ is the sutured manifold as shown in Figure \ref{fig: genus two handle body}, then
$$\shi(-M,-\ga)\cong\mathbb{C}^5.$$
\ecor

\section{Polytopes}\label{sec: polytopes}
\subsection{Constructing the polytope}\label{subsec: polytope}
Suppose $(M,\ga)$ is a balanced sutured manifold, and $\al\in H_2(M,\partial M)$ is a second relative homology class. From Theorem \ref{thm: general grading shifting property}, there is a grading on $\shm(M,\ga)$ or $\shi(M,\ga)$ associated to $\al$, which is well-defined up to an overall grading shifting. Hence, we can define the following.
\bdefn\label{defn: homogenous element}
An element $a$ in $\shm(M,\ga)$ or $\shi(M,\ga)$ is called {\it homogenous} if for any homology class $\al\in H_2(M,\partial M)$, the element $a$ is homogenous with respect to the grading associated to $\al$.
\edefn

\blem\label{lem: homogenous element exists}
For any taut balanced manifold $(M,\ga)$, there exist non-zero homogenous elements in $\shm(M,\ga)$ and $\shi(M,\ga)$.
\elem
\bpf
Pick a basis $\al_1,...,\al_{n}$ for $H_2(M,\partial M)$. Pick admissible properly embedded surfaces $S_1,...,S_{n}$ to represent $\al_1,...,\al_n$, respectively. In the proof of Theorem \ref{thm: general grading shifting property}, for two surfaces $S_1$ and $S_2$ inside $(M,\ga)$, we constructed a closure $(Y,R_+)$ of $(M,\ga)$ where (suitable isotopies of) $S_1$ and $S_2$ both extend to closed surfaces inside $Y$. we can carry out a similar construction here, though this time we have $n$ surfaces instead of two. As a result, we obtain a special closure $(Y,R)$ of the balanced sutured manifold $(M,\ga)$ so that, inside $Y$, the surfaces $S_1,...,S_n$ extends to closed surfaces $\bar{S}_1,...,\bar{S_n}$, respectively, and the surfaces $\bar{S}_1,...,\bar{S_n}$ are the ones used to define gradings on $\shm(M,\ga)$ associated to $S_1,...,S_n$, respectively. Also, pick some suitable non-separating simple closed curves $\eta$ on $R$ to support local coefficients. On $Y$, the surfaces $\bar{S}_1,...,\bar{S}_n$ induces a $\intg^r$-grading on $HM(Y|R;\Ga_{\eta})$ by looking at the evaluation of the first Chern classes of spin${}^c$ structures on those closed surfaces. Suppose $\mathfrak{s}_0$ is a supporting spin${}^c$ structures (see Definition \ref{defn: SHM}) on $Y$, and $a\in \widecheck{HM}_{\bullet}(Y,\mathfrak{s}_0;\Ga_{\eta})$ is a non-zero element, then we know that the element $a$ is homogenous with respect to all the gradings induced by $\bar{S}_1,...,\bar{S}_n$. We claim that it is a homogenous element as defined in Definition \ref{defn: homogenous element}.

To prove the claim, suppose $\al\in H_2(M,\ga)$ is any homology class. Since $\al_1,...,\al_n$ form a basis of $H_2(M,\ga)$, $\al$ is a linear combination of $\al_1,...,\al_n$. Thus, we can perform a sequence of double curve surgeries (for definition, see Scharlemann \cite{scharlemann1989suture}) on a few parallel copies of  $S_1,...,S_n$, with the number of copies depending on coefficients of $\al$, to obtain a properly embedded surface $S$ that represents the class $\al\in H_2(M,\ga)$. Correspondingly, we can perform the same set of double curve surgeries on $\bar{S}_1,...,\bar{S_n}$ to obtain a closed surface $\bar{S}\subset Y$, which extends $S$ and which induces the grading associated to $S$. Then, we know that the element $a$ is a homogenous element with respect to the grading associated to $S$, and this concludes the proof of lemma \ref{lem: homogenous element exists}.
\epf

\bdefn\label{lem: difference of homogenous element}
Suppose $(M,\ga)$ is a balanced sutured manifold, and $a,b\in\shm(M,\ga)$ are two homogenous elements. Then, we define an element $\rho(a,b)\in H^2(M,\partial{M};\mathbb{Q})$ associated to the (ordered) pair $(a,b)$ as follows:
We first construct the map
$$\rho(a,b):H_2(M,\partial M)\ra \intg.$$
For any class $\al\in H_2(M,\partial{M})$, we pick a surface $S$ that represents the class $\al$ and is admissible. Define
$$\rho(a,b)(\al)={\rm difference~between~}a~{\rm and}~b~{\rm under~the~grading~associated~to~}S.$$
This is well defined by Theorem \ref{thm: general grading shifting property}. This map is linear by essentially the same type of argument as in the proof of Lemma \ref{lem: homogenous element exists}. Then, we can regard $\rho(a,b)$ as an element in $H^2(M,\partial{M};\mathbb{Q})$.

We can carry out a similar construction in the instanton setups.
\edefn

\bdefn\label{defn: polytope}
Suppose $(M,\ga)$ is a balanced sutured manifold, and $a\in\shm(M,\ga)$ is a homogenous element. For an element $\rho\in H^2(M,\partial{M};\real)$, define
$$\shm_a(M,\ga,\rho)=\{b\in\shm(M,\ga),~\rho(a,b)=\rho\in H^2(M,\partial{M};\real)\},$$
and
$$\shi_a(M,\ga,\rho)=\{b\in\shi(M,\ga),~\rho(a,b)=\rho\in H^2(M,\partial{M};\real)\}.$$
Let 
$$SM_a(M,\ga)=\{\rho\in H^2(M,\partial M;\real),~\shm_a(M,\ga,\rho)\neq0\}$$
and
$$SI_a(M,\ga)=\{\rho\in H^2(M,\partial M;\real),~\shi_a(M,\ga,\rho)\neq0\}.$$
Define the {\it polytopes} $PM_a(M,\ga)$ and $PI_a(M,\ga)$ to be the convex hull of $SM_a(M,\ga)$ and $SI_a(M,\ga)$, respectively. 
\edefn

\blem
Suppose $a$ and $b$ are two homogenous elements in $\shm(M,\ga)$, then the polytopes $PM_a(M,\ga)\subset H^2(M,\partial M;\real)$ is a translate of $PM_b(M,\ga)\subset H^2(M,\partial M;\real)$. The same result holds in the instanton setups.
\elem
\bpf
It is straightforward from the construction.
\epf

\subsection{Dimension formula}
\blem\label{lem: decompose along product annuli}
Suppose $(M,\ga)$ is a taut balanced sutured manifold with $H_2(M)=0$, and $A$ is an incompressible product annulus. Then, we can pick an orientation of $A$ so that the sutured manifold decomposition of $(M,\ga)$ along the oriented $A$ yields a taut balanced sutured manifold $(M',\ga')$, that $SHM(M',\ga')$ is a direct summand of $SHM(M,\ga)$, and that $SHI(M',\ga')$ is a direct summand of $SHI(M,\ga)$.
\elem
\bpf
Since $A$ is incompressible, we know that no components of $\partial{A}\subset R(\gamma)$ bound a disk. Note that \cite[Lemma 4.2]{scharlemann2006three} makes sure that no matter which orientation of $A$ we choose, the balanced sutured manifold after the decomposition is taut. There are three cases:

{\bf Case 1}. Both components of $\partial{A}$ are homologically essential on $R(\ga)$. 

{\bf Case 2}. Both components of $\partial{A}$ are homologically trivial on $R(\ga)$. Then, there are $V_+\subset R_+(\ga)$ and $V-\subset R_-(\ga)$ so that $\partial{V}_+\cup \partial{V}_-=\partial{A}$ (as unoriented curves). Thus, we have a closed surface $V_+\cup A\cup V_-$. The fact that $H_2(M)=0$ implies that this closed surface is separating, or equivalently, $A$ separates $M$ into two parts, of which one has boundary $V_+\cup A\cup V_-$. Thus, this part is disjoint from $\ga$.

For the above two cases, the lemma follows from Proposition 6.7 of Kronheimer and Mrowka \cite{kronheimer2010knots}. Actually in these two cases $SHI(M',\ga')$ is isomorphic to $SHI(M,\ga)$.

{\bf Case 3}. One component of $\partial{A}$ is homologically essential, and the other is inessential. Then we can choose a suitable orientation of $A$ to make $\partial{A}$ being boundary coherent, in the sense of Kronheimer and Mrowka \cite{kronheimer2010knots}, so that Proposition 6.9 in that paper applies, and, thus, we conclude the proof of Lemma \ref{lem: decompose along product annuli}.
\epf

\blem\label{lem: adding product one handle does not change SHM}
Suppose $(M,\ga)$ is a taut balanced sutured manifold and $S\subset M$ is a properly embedded decomposing surface. Suppose ${p,q}\subset S\cap \ga$ are two points of different signs. Then, we can attach a product $1$-handle, in the sense of Definition \ref{defn: one handle}, to obtain a new taut balanced sutured manifold $(M_1,\ga_1)$ and a new properly embedded surface $S_1\subset M_1$ so that the decomposition
$$(M,\ga)\stackrel{S}{\leadsto}(M',\ga')$$
is taut if and only if the decomposition
$$(M_1,\ga_1)\stackrel{S_1}{\leadsto}(M'_1,\ga'_1)$$
is taut. Furthermore, there is a commutative diagram
\begin{equation*}
\xymatrix{
SHM(M',\ga')\ar@{^{(}->}[r]\ar[d]^{=}& SHM(M,\ga)\ar[d]^{=}\\
SHM(M_1',\ga_1')\ar@{^{(}->}[r]&SHM(M_1,\ga_1).
&
}    
\end{equation*}

\begin{figure}[h]
\centering
    \begin{overpic}[width=5.5in]{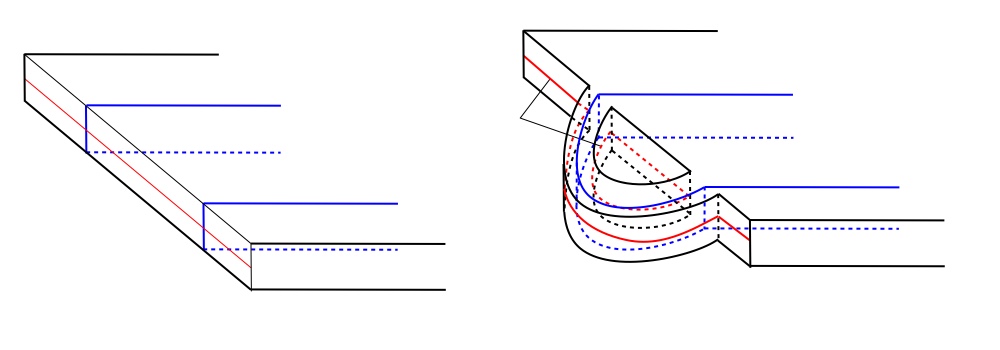}
    \put(5,25){\line(0,-1){7}}
    \put(4,16){$\ga$}
    \put(7,22){\line(0,-1){10}}
    \put(3,10){$A(\ga)$}
    \put(25,25){$S$}
    \put(35,15){$S$}
    \put(50,21){$\ga_1$}
    \put(75,26){$S$}
    \put(85,16.7){$S$}
\end{overpic}
\vspace{0.05in}
\caption{Adding a product $1$-handle.}\label{fig: product one handle}
\end{figure}
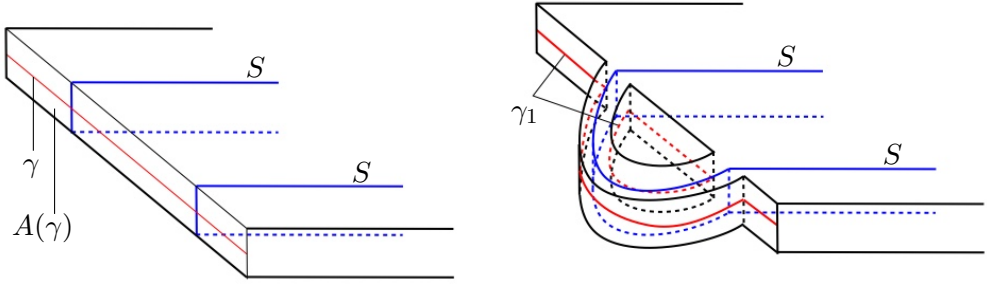

A similar statement holds in the instanton settings.
\elem
\bpf
This is how Kronheimer and Mrowka proved Proposition 6.9 in \cite{kronheimer2010knots}. $S_1$ is obtained from $S$ by attaching a $2$-dimensional $1$-handle inside the $3$-dimensional product $1$-handle. See Figure \ref{fig: product one handle}.

\epf

\bprop\label{prop: distinct direct summand}
Suppose $(M,\ga)$ is a balanced sutured manifold with $H_2(M)=0$ and is taut, reduced, horizontally prime, and free of non-separating essential product disks. Suppose $\al\in H_2(M,\partial{M})$ is a non-zero class. Then, we can find properly embedded surfaces $S$ and $S'$ in $M$ so that

(1) $[S]=-[S']=\al\in H_2(M,\partial{M})$.

(2) The sutured manifold decompositions
$$(M,\ga)\stackrel{S}{\leadsto}(M',\ga')~{\rm and}~(M,\ga)\stackrel{S'}{\leadsto}(M'',\ga'')$$
are both taut.

(3) $SHM(M',\ga')$ and $SHM(M'',\ga'')$ are direct summands of $SHM(M,\ga)$, (by \cite[Proposition 6.9]{kronheimer2010knots}) and
$$SHM(M',\ga')\cap SHM(M'',\ga'')={0}$$
in $SHM(M,\ga)$.

(4) The same result holds for $SHI$.
\eprop

\bpf
This will be a very long proof, so we first sketch the proof as follows. First, The proof is parallel to the proof of \cite[Theorem 6.1]{juhasz2010polytope}, though we make some major modifications to adapt to the monopole and instanton setups. We start with a non-trivial class $\al\in H_2(M,\partial M)$. Gabai \cite{gabai1983foliations} introduced a way to find well-groomed surfaces $S$ and $S'$, representing the class $\al$ and $-\al$, respectively, so that the decompositions of $(M,\ga)$ along $S$ and $S'$ are both taut. From \cite[Proposition 6.9]{kronheimer2010knots}, these two taut decompositions provides two summands of $SHM(M,\ga)$. To un-package the proof of \cite[Proposition 6.9]{kronheimer2010knots}, we will find a closure $(Y,R)$ of $(M,\ga)$ so that $S$ and $S'$ both extend to some closed surfaces $\bar{S}$ and $\bar{S}'$, respectively. In order to show that these two summands from two taut decompositions are distinct [conclusion (3) of the proposition], we need to show that there is no spin${}^c$ structure on $Y$ whose first Chern class has the desired evaluation on both $\bar{S}$ and $\bar{S}'$. This is done via the adjunction inequality in \ref{lem: adjunction inequality in monopoles} and an inequality from \cite{juhasz2010polytope}, see \ref{lem: inequality on the Euler characteristics}. In the construction of $\bar{S}$ and $\bar{S}'$, for some technical reasons, we need balance two values $s$ and $t$. The idea is to study some basic pieces $(V,\delta)$ of a solid torus with two longitudinal sutures. By adding enough copies of this standard piece, we can make $s$ equals to $t$ and complete the construction of $\bar{S}$ and $\bar{S}'$.

Now suppose $\al$ is a non-trivial class in $H_2(M,\partial M)$. From Lemma 0.7 in Gabai \cite{gabai1987foliations2}, we can pick an $S$ so that the following holds

(i) $S$ represents the class $\al\in H_2(M,\partial{M})$.

(ii) For any component $V$ of $R(\ga)$, if $S\cap V$ has a closed component, then $S\cap V$ consists of parallel, parallel oriented non-separating simple closed curves.

(iii) For any component $\delta$ of $\partial{R(\ga)}$, all the intersection points of $S$ with $\delta$ are of the same sign.

(iv) $S$ intersects $A(\ga)$ in parallel and coherently oriented essential arcs.

(v) The sutured manifold decomposition
$$(M,\ga)\stackrel{S}{\leadsto}(M',\ga')$$
is taut.

\begin{figure}[h]
\centering
    \begin{overpic}[width=5.5in]{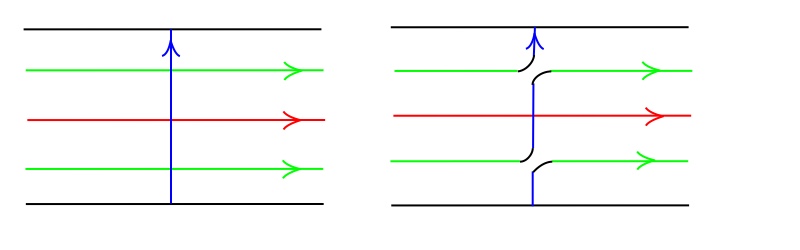}
    \put(41.5,22){$R(\ga)$}
    \put(41.5,27){$A(\ga)$}
    \put(44,16){$\ga$}
    \put(41.5,10){$R(\ga)$}
    \put(41.5,5){$A(\ga)$}
    \put(19,2){$-\partial{S}$}
\end{overpic}
\caption{The smoothing inside $A(\ga)$.}\label{fig: smoothing}
\end{figure}
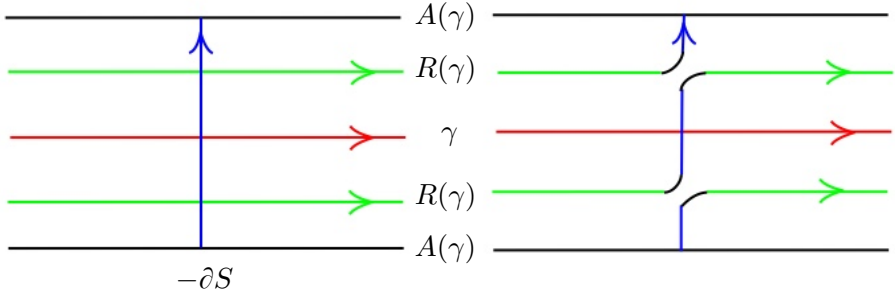

To find $S'$, we proceed as follows, according to Gabai \cite{gabai1983foliations}, Scharlemann \cite{scharlemann2006three} or Juhasz \cite{juhasz2010polytope}. Pick $U=N(\partial{S})\cup A(\ga)\subset \partial{M}$. Pick a large enough $k\in\intg_+$ so that
$$x(-\al+(k+1)[R(\ga)])=x(-\al+k[R(\ga)])+x(R(\ga)),$$
where $x(\cdot)$ is the Thurston norm for classes in $H_2(M,U)$. Pick a norm-minimizing embedded surface which represents the class $-\al+k[R(\ga)]$, and disregard all components of it, which represent the zero homology class, then the remaining surface $S'$ is the desired one. Note
$$[\partial{S}']=-[\partial{S}]+2k\cdot[\ga]\in H_1(U).$$
Here, it is $2k$ rather than $k$, because each copy of $R(\ga)$ contributes $2[\ga]$. We can arrange so that $\partial{S}'$ is obtained from $-\partial{S}$ and $2k$ copies of $\ga$ by an oriented smoothing as in Figure \ref{fig: smoothing}. From construction of $S'$, we know that the sutured manifold decomposition
$$(M,\ga)\stackrel{S'}{\leadsto}(M'',\ga'')$$
is taut.

\begin{figure}[h]
\centering
    \begin{overpic}[width=5.5in]{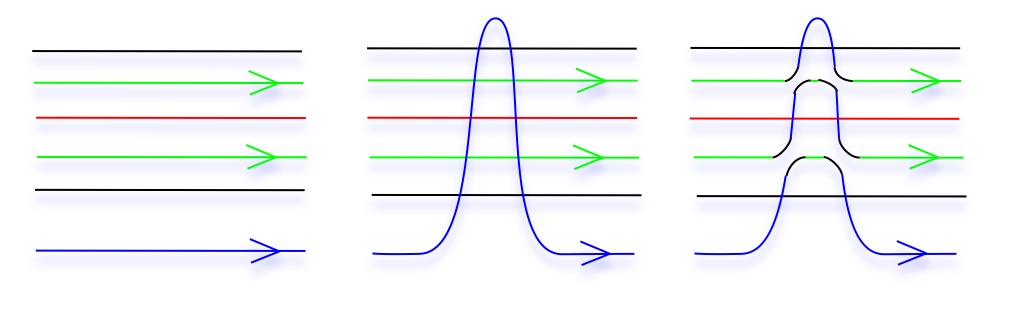}
    \put(-2,23.2){$\partial S'$}
    \put(1,19.5){$\ga$}
    \put(-2,15.5){$\partial S'$}
    \put(-0.5,6.5){$\partial S'$}
    \put(37,2){Positive stabilization}
    \put(69,2){Double curve surgery}
\end{overpic}
\caption{A positive stabilization on $S'$ followed by a double curve surgery.}\label{fig: positive stabilization and double curve surgery}
\end{figure}
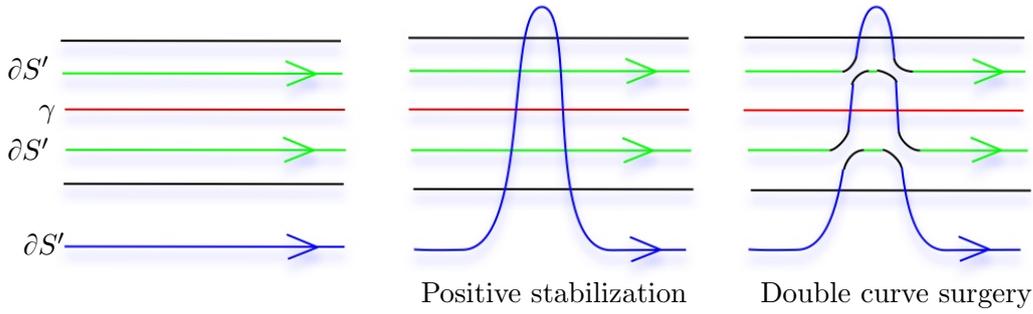

From the above construction of $S'$, we know that $S\cap\ga=S'\cap \ga$. So, assume that
$$n=\frac{1}{2}|S\cap\ga|=\frac{1}{2}|S'\cap\ga|.$$
Also, assume that $\partial{S}\cap R(\ga)$ (and thus $\partial{S}'\cap R(\ga)$) has $m$ closed components. Write them as
$$B_1\cup...\cup B_m\subset \partial{S}\cap R(\ga),$$
and orient $B_i$ by the boundary orientation of $S$. Write $S_0$ the surface obtained from $S$ by performing one negative stabilization on each $B_i$, and let $S^{\dag}$ be the surface obtained by performing one positive stabilization on each $B_i$. By \cite[Lemma 4.5]{juhasz2008floer}, both $S_0$ and $S^{\dag}$ exist (i.e., the negative and positive stabilizations do exist). We also want to modify $S'$ correspondingly. Note part of $\partial{S}'$ is coming from $-\partial S$, so positive stabilizations near this part of $\partial S'$ corresponds to negative stabilizations on $S$. Thus, we perform positive stabilizations on $S'$, in correspondence to the negative stabilizations performed on $S$. Since $\partial{S'}$ contains many copies of $\partial{R}(\ga)\subset A(\ga)$, when performing the positive stabilizations  to create intersections of $B_i$ with $\ga$, the isotopy also creates self intersections of $S'$. We then perform double curve surgeries, in the sense of Scharlemann \cite{scharlemann2006three}, to resolve all self intersections created by the positive stabilizations, and let $S_0'$ be the resulting surface. See Figure \ref{fig: positive stabilization and double curve surgery}. On $A(\ga)$, this double curve surgery behaves exactly in the same way as the oriented smoothing that is depicted in Figure \ref{fig: smoothing}. It is important that we choose positive stabilizations to perform on $S'$, so the decomposition of $(M,\ga)$ along $S_0'$ gives the same result as decomposing along $S'$:
$$(M,\ga)\stackrel{S'_0}{\leadsto}(M'',\ga'').$$

We can attach $m+n$ copies of product $1$-handle, as in Lemma \ref{lem: adding product one handle does not change SHM}, along the intersections of $S_0$ with $\ga$. We require that the pair of intersection points created by a negative stabilization on $S$ are paired together by a product $1$-handle. Let $(M_1,\ga_1)$ be the resulting balanced sutured manifold. The surface $S_0$ extends to a properly embedded surface $S_1\subset M_1$ as in Lemma \ref{lem: adding product one handle does not change SHM}.

\begin{figure}[h]
\centering
    \begin{overpic}[width=5.5in]{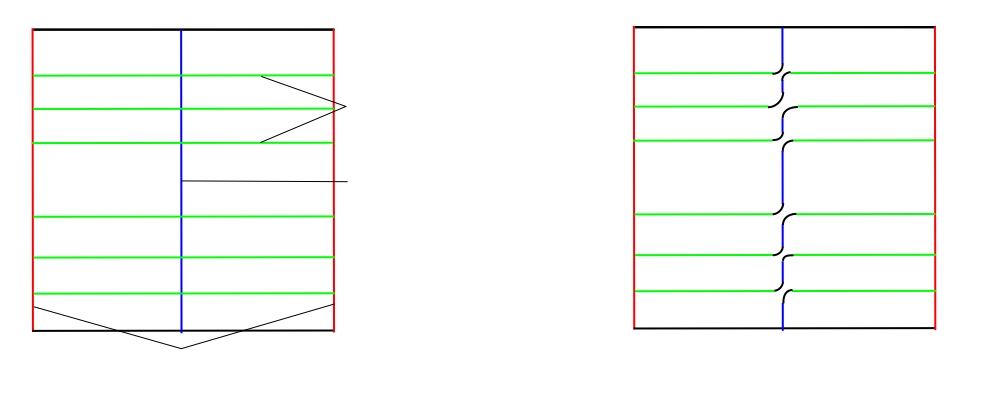}
    \put(15,2){$A(\ga)$}
    \put(36,21){$-S_1$}
    \put(36,29){$R(\ga_1)$}
    \put(38,18){\vector(1,0){23}}
    \put(44,19){smoothing}
\end{overpic}
\vspace{0.05in}
\caption{Double curve surgery on a cross section of the product $1$-handle. On the right: the $2$-dimensional $1$-handles after the double curve surgery will cut the original product $1$-handle into $2k+2$ small ones.}\label{fig: double curve surgery in product one handle}
\end{figure}
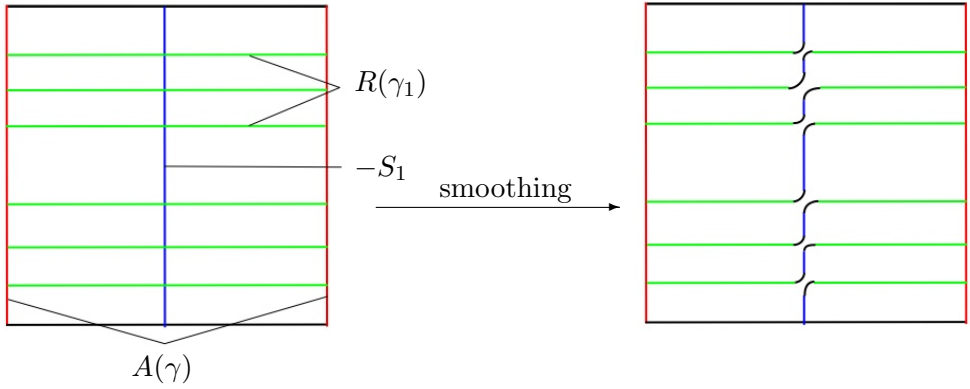

The surface $S'_0$ can also extend, though in a slightly complicated way. As in Figure \ref{fig: double curve surgery in product one handle}, in each product $1$-handle, there is one (vertical) $2$-dimensional $1$-handle to be glued to the part of $\partial{S}'_0$ that comes from $-\partial{S}$, and $2k$ copies of (horizontal) $2$-dimensional $1$-handles to be glued to the part of $\partial{S}'_0$ that corresponds to the $k$ copies of $\partial R(\ga)$. We can perform a double curve surgery on those two collections of 2-dimensional $1$-handles, as in Figure \ref{fig: double curve surgery in product one handle}, and then glue the resulting surface to $S'_0$, when gluing the product $1$-handle to $(M,\ga)$. In this way, $S'_0$ extends to a properly embedded surface $S_1'$ in $(M_1,\ga_1)$. Let $(M_1',\ga_1')$ and $(M_1'',\ga_1'')$ be obtained from $(M_1,\ga_1)$ by decomposing along $S_1$ and $S_1'$, respectively. $(M_1'',\ga_1'')$ is taut by Lemma \ref{lem: adding product one handle does not change SHM}, since it is obtained from $(M'',\ga'')$ by attaching a few product $1$-handles, as indicated in Figure \ref{fig: double curve surgery in product one handle}. After decomposing along $S_1'$, each product $1$-handle is decomposed into $(2n+1)$ many product $1$-handles that are attached to $(M'',\ga'')$. In that figure, $n=3$, and the six arcs in the right-subfigure divides the original product $1$-handle into $7$ parts.

From the construction, we also know that
$$S_1'\cap R_{\pm}(\ga_1)=-S_1\cap R_{\pm}(\ga_1),$$
and $S_1'\cap A(\ga_1)$ consists of $2k$ parallel copies of $\ga_1$. Now let
\begin{equation}\label{eq: the value s}
	S_1\cap R_{+}(\ga_1)=B_{+,1}\cup...\cup B_{+,m}\cup C_{+,1}\cup...\cup C_{+,s}
\end{equation}

\begin{equation}\label{eq: the value t}
	S_1\cap R_{-}(\ga_1)=B_{-,1}\cup...\cup B_{-,m}\cup C_{-,1}\cup...\cup C_{-,t}
\end{equation}

Here, $B_{\pm,i}$ are the boundary components of $S_1$ that come from attaching a product $1$-handle along the pair of intersection points created by a negative stabilization on $S$ near $B_i$. Note $s$ and $t$ are not necessarily equal. Without loss of generality, we could assume that $s\geq t$. 

Note in order to obtain a closed surface from $S_1$ inside a suitable closure of $(M,\ga)$, we need to require that $s=t$ as in Section \ref{subsec: construction of gradings}. However, we might have the possibility that $s>t$. As explained at the beginning, we need to add $(s-t)$ many copies of the following standard piece $(V,\delta)$ into consideration.

\begin{figure}[h]
\centering
    \begin{overpic}[width=5.5in]{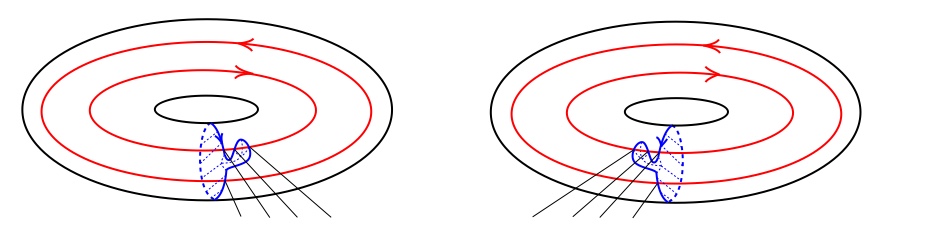}
    \put(21,2){$D$}
    \put(69,2){$D^{\dag}$}
    \put(24,2){$p_4$}
    \put(27,2){$p_1$}
    \put(30,2){$p_2$}
    \put(34,2){$p_3$}
    \put(54,2){$p_3$}
    \put(58,2){$p_2$}
    \put(62,2){$p_1$}
    \put(65,2){$p_4$}
\end{overpic}
\vspace{0.05in}
\caption{The balanced sutured manifold $(V,\delta)$ together with the surfaces $D_0$ and $D^{\dag}$.}\label{fig: add on product manifold}
\end{figure}
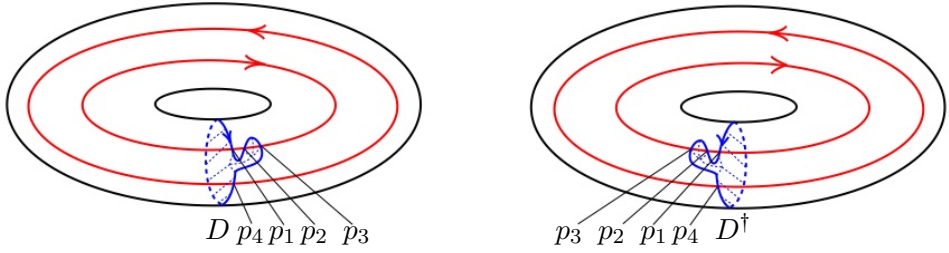

Pick $(V,\delta)$ be a balanced sutured manifold where $V=S^1\times D^2$ is a solid torus and $\delta$ consists of two longitudes. Let $D\subset V$ be a standard meridian disk in $V$, which has two intersections with the suture $\delta$. Let $D_0$ be the surface obtained from $D$ by performing a negative stabilization, as shown in Figure \ref{fig: add on product manifold}. It has four intersections with the suture $\delta$. Attach two product $1$-handle along two pairs of points $(p_1,p_4)$ and $(p_2,p_3)$, as labeled in the figure, and let the resulting balanced sutured manifold be $(V_1,\ga_1)$. $D_0$ extends to a properly embedded surface $D_1$ as in Lemma \ref{lem: adding product one handle does not change SHM}. Note $D_1\cap R_+(\delta_1)$ has one components and $D_1\cap R_-(\delta_1)$ has two components. Now we want to construct another surface $D_1'$ inside $V_1$. First, let $D'$ be the result of a double curve surgery of $-D_0$ with $k$ copies of $R(\delta)$. It is crucial that our choice of $D_0$ makes the sutured manifold decomposition
$$(V,\delta)\stackrel{D'}{\leadsto}(V'',\ga'')$$
taut. We can check directly that $V''$ is a $3$-ball and $\ga''$ is a simple closed curve on $\partial{V}''$, and, hence, it is a product sutured manifold. Next, when attaching two product $1$-handles to $(V,\delta)$, we repeat the procedure explained above, with which we construct $S'_1$ out of $S'_0$, and thus construct a properly embedded surface $D_1'$ inside $V_1$ out of $D_0'$. It is straightforward to check that the decomposition of $(V_1,\ga_1)$ along $D_1'$ yields a product sutured manifold $(V_1'',\ga_1'')$.

So, we could form the disjoint union
$$(M_2,\ga_2)=(M_1,\ga_1)\sqcup (s-t)(V_1,\delta_1),~S_2=S_1\cup (s-t)D_1,~S_2'=S_1'\cup (s-t)D_1'.$$
Since the decomposition of $(M_1,\ga_1)$ along $S_1'$ is taut (as we have explained) and the decomposition of $(V_1,\ga_1)$ along $D_1'$ is also taut (by a direct check), we know that the decomposition of $(M_2,\ga_2)$ along $S_2'$ is taut. However, the decomposition along $S_2$ is not, as we will explain later.

Pick a connected auxiliary surface $T$ for $(M_2,\ga_2)$ and form a pre-closure
$$\widetilde{M}=M_2\cup [-1,1]\times T,~\partial{\widetilde{M}}=R_+\cup R_-,~R_{\pm}=R_{\pm}(\ga_2)\cup \{\pm1\}\times T.$$

Since $(M_1,\ga_1)$ is obtained from $(M,\ga)$ by attaching product $1$-handles, and $(V_1,\delta_1)$ is itself a product sutured manifold, we know that $\widetilde{M}$ is also a pre-closure of $(M,\ga)$. By construction, we know that
$$S_2\cap R_{\pm}=-S_2'\cap R_{\pm},$$
and there are same number of components of $S_2\cap R_+$ and $S_2\cap R_-$. Moreover, since $T$ is connected, the components of $S_2\cap R_{\pm}$ represent linearly independent classes in $H_1(R_{\pm})$. Thus, we could find an orientation preserving diffeomorphism $h:R_+\ra R_-$ so that
$$h(S_2\cap R_+)=S_2\cap R_-,$$
and we can use $\widetilde{M}$ as well as $h$ to construct a closure $(Y,R_+)$ for both $(M,\ga)$ and $(M_2,\ga_2)$. Inside $Y$, $S_2$ becomes a closed surface $\bar{S}$.

Though there is no sutures on the boundary of $\widetilde{M}$, the theory of balanced sutured manifolds in Kronheimer and Mrowka \cite{kronheimer2010knots} extends to $\widetilde{M}$ effectively. In particular, we could define 
$$SHM(\widetilde{M})=HM(Y|R_+)=SHM(M,\ga)=SHM(M_2,\ga_2).$$
The surface $S_2'$ extends to a surface $\widetilde{S}_2'\subset \widetilde{M}$ as the union of $S_2'$ with $2k$ copies of $T$. Suppose $(\widetilde{M}'',\tilde{\ga}'')$ is the result of the sutured manifold decomposition of $\widetilde{M}$ along $\widetilde{S}_2'$, then $(\widetilde{M}'',\tilde{\ga}'')$ can be obtained from $(M_2'',\ga_2'')$ by attaching $2k+1$ copies of the product region $T\times[-1,1]$. Recall that $(M_2'',\ga_2'')$ is obtained from $(M_2,\ga_2)$ by decomposing along $S_2'$ and thus is the disjoint union of $(M_1'',\ga_1'')$ with $(s-t)$ copies of product sutured manifolds $(V_1'',\ga_1'')$. Furthermore, $(M_1'',\ga_1'')$ is obtained from $(M'',\ga'')$ by attaching a few product $1$-handles, so we finally conclude that
$$SHM(\widetilde{M}'',\widetilde{\ga}'')\cong SHM(M'',\ga''),$$
since attaching product regions (or product $1$-handles), disjoint union with product manifolds will never change the sutured monopole Floer homology.

Back to the point that $(\widetilde{M}'',\tilde{\ga}'')$ is the decomposition of $\widetilde{M}$ along $\widetilde{S}_2'$. The decomposition theorem, Proposition 6.9 in Kronheimer and Mrowka \cite{kronheimer2010knots}, continues to hold in this case, and we conclude that $SHM(\widetilde{M}'',\widetilde{\ga}'')$ is a direct summand of $SHM(\widetilde{M})$. More precisely, $\widetilde{S}_2'$ becomes a closed surface $\bar{S}'\subset Y$, and we have
\begin{equation}\label{eq: summand corresponds to S prime}
SHM(M'',\ga'')=SHM(\widetilde{M}'',\widetilde{\ga}'')=\bigoplus_{\substack{\mathfrak{s}\in\mathfrak{S}^*(Y|R_+)\\c_1(\mathfrak{s})[\bar{S}']=2g(\bar{S}')-2}}\widecheck{HM}_{\bullet}(Y,\mathfrak{s}).
\end{equation}

We also want to identify the summand $SHM(M',\ga')$ inside $SHM(\widetilde{M})=SHM(M,\ga)$. We cannot proceed directly as we did for $SHM(M'',\ga''),$ since the decomposition of $\widetilde{M}$ along $S_2$ is not taut. This is because the decomposition of $(M,\ga)$ along $S_0$ and $(V,\ga)$ along $D_0$ are both not taut, since, at the beginning, we picked $S_0$ and $D_0$ by performing negative stabilizations (see Lemma \ref{lem: decomposition and pm stabilization}). Let $S^{\dag}$ and $D^{\dag}$ be obtained from $S$ and $D$ performing positive stabilizations instead of negative ones. We can repeat the whole construction again with $S_0$ and $D_0$ replaced by $S^{\dag}$ and $D^{\dag}$, respectively. Attach $m+n$ product $1$-handle along the intersection points of $\partial S^{\dag}$ with $\ga$, and let $(M_1^{\dag},\ga_1^{\dag})$ be the result. There is a properly embedded surface $S_1^{\dag}\subset M_1^{\dag}$. Similarly, attach two product $1$-handles to $(V,\ga)$, and let $(V^{\dag}_1,\ga^{\dag}_1)$ be the result. Then, there is a properly embedded surface $D_1^{\dag}\subset V_1^{\dag}$. We can form the disjoint union
$$(M_2^{\dag},\ga_2^{\dag})=(M_1^{\dag},\ga^{\dag}_1)\sqcup (s-t)(V_1^{\dag},\ga_1^{\dag}),~S^{\dag}_2=S_1^{\dag}\cup(s-t)D_1^{\dag}.$$
Pick an auxiliary surface $T^{\dag}$ that has the same genus as $T$ and form a pre-closure $\widetilde{M}^{\dag}$. We can pick a suitable gluing diffeomorphism $h^{\dag}$ and obtain a closure $(Y^{\dag},R_+^{\dag})$. Inside $Y^{\dag}$, the surface $S_2^{\dag}$ becomes a closed surface $\bar{S}^{\dag}$.The decomposition of $\widetilde{M}^{\dag}$ along $S^{\dag}_2$ is taut. The argument that concludes (\ref{eq: summand corresponds to S prime}) applies again, and we have
\begin{equation}\label{eq: summand corresponds to S dag}
SHM(M',\ga')=\bigoplus_{\substack{\mathfrak{s}\in\mathfrak{S}^*(Y^{\dag}|R^{\dag}_+)\\c_1(\mathfrak{s})[\bar{S}^{\dag}]=2g(\bar{S}^{\dag})-2}}\widecheck{HM}_{\bullet}(Y^{\dag},\mathfrak{s}).
\end{equation}

Since $g(T^{\dag})=g(T)$, we know that $g(R_+^{\dag})=g(R_+)$. Thus, as in Subsection \ref{subsec: SHM and SHI}, there is an excision cobordism $W$, from $Y^{\dag}\cup Y_T$ to $Y$, which induces an isomorphism
$$\Phi:HM(Y^{\dag}|R^{\dag}_+)\ra HM(Y|R_+).$$
Here, $Y_T$ is a mapping torus of a diffeomorphism on $R_{+}^{\dag}$, arising from $h^{\dag}$, $h$, and a suitable identification $R_+^{\dag}=R_+$. The difference between the surfaces $\bar{S}^{\dag}$ and $\bar{S}$ originates from whether we performed positive or negative stabilizations. So, the proof of Proposition 4.1 in Li \cite{li2019direct} applies to the present context, and we know that inside $W$,
$$[\bar{S}^{\dag}]=[\bar{S}]+[\Sigma_1]+...[\Sigma_m]+[\Sigma'_1]+...+[\Sigma'_{s-t}].$$
Here, each $\Sigma_i$ or $\Sigma'_j$ is a connected closed oriented surface of genus $2$. $\Sigma_i$ corresponds to a positive or negative stabilization on $S$, and $\Sigma_j'$ corresponds to a positive or negative stabilization on $D$. As a result of the adjunction inequality in Lemma \ref{lem: adjunction inequality in monopoles}, we have

\begin{equation}\label{eq: summand corresponds to S}
SHM(M',\ga')\subset\bigoplus_{\substack{\mathfrak{s}\in\mathfrak{S}^*(Y|R_+)\\c_1(\mathfrak{s})[\bar{S}]\geq2g(\bar{S})-2-2(s-t)-2m}}\widecheck{HM}_{\bullet}(Y,\mathfrak{s}).
\end{equation}

Finally, we argue that
$$SHM(M',\ga')\cap SHM(M'',\ga'')=\{0\}.$$
Suppose not, then, from (\ref{eq: summand corresponds to S prime}) and (\ref{eq: summand corresponds to S}), there is a supporting spin${}^c$ structure $\mathfrak{s}\in S^{*}(Y|R_+)$ so that
\begin{equation*}
c_1(\mathfrak{s})[\bar{S}]\geq 2g(\bar{S})-2-2(s-t)-2m,~{\rm and}~c_1(\mathfrak{s})[\bar{S}']=2g(\bar{S}')-2.
\end{equation*}
From the construction, we know that 
$$[\bar{S}']=-[\bar{S}]+2k\cdot[R_+]\subset H_2(Y).$$
Hence, the above equalities and inequalities imply
$$2g(\bar{S})-2-2(s-t)-2m+2g(\bar{S}')-2\leq c_1(\mathfrak{s})[\bar{S}]+c_1(\mathfrak{s})[\bar{S}']=2k\cdot [2g(R_+)-2],$$
which is equivalent to
\begin{equation}\label{eq: inequality on the Euler characteristics}
\chi(\bar{S})+\chi(\bar{S}')+2(s-t)+2m\geq 2k\cdot \chi(R_+).    
\end{equation}
Now let us compute each of the three terms in (\ref{eq: inequality on the Euler characteristics}) regarding the Euler characteristics. First, $\chi(\bar{S})=\chi(S_2)$, and, by construction, $S_2=S_1\sqcup (s-t)D_1$. Furthermore, we know that $S_1$ is obtained from $S_0$ by attaching $m+n$ copies of $2$-dimensional $1$-handles, that $S_0$ is isotopic to $S$, and that $D_1$ is obtained from a disk $D$ by attaching two copies of $2$-dimensional $1$-handles. Thus, we conclude that
$$\chi(S_2)=\chi(S)-(m+n)+(s-t)(-1)=\chi(S)-(m+n)-(s-t).$$

Second, we know that $\chi(\bar{S}_2')=\chi(\widetilde{S}_2')$, and
$$\widetilde{S}_2'=S_1'\cup D_1'\cup (2k)\cdot T.$$
Note that there are $m+n$ product $1$-handles attached to $(M,\ga)$, and, inside each product $1$-handle, there are $(2k+1)$ copies of $2$-dimensional $1$-handles attached to $S'$. Thus, we have that
$$\chi(S_1')=\chi(S')-(2k+1)(m+n).$$
Similarly, we conclude that
$$\chi(D_1')=\chi(D')-2(2k+1).$$
Also, $D'$ is obtained by a double curve surgery on $-D$, which is a disk, with $k$ copies of $R(\delta)$, which is the disjoint union of two annuli. Thus, we conclude that
$$\chi(\bar{S}')=\chi(S')-(m+n)-(s-t)+2k\cdot\chi(T)-2k(m+n)-4k(s-t).$$

Third, we know that 
$$R_+=R_+(\ga_1)\cup (s-t)\cdot R_+(\delta_1)\cup T.$$
Here, $R_+(\ga_1)$ is obtained from $R_+(\ga)$ by attaching $m+n$ copies of $2$-dimensional $1$-handles, and $R_+(\delta_1)$ is obtained from $R_+(\delta)$ (an annulus) by attaching $2$ copies of $2$-dimensional $1$-handles. Thus, we know that
$$\chi(R_+)=\chi(R_+(\ga))+\chi(T)-(m+n)-2(s-t).$$

Putting everything together, (\ref{eq: inequality on the Euler characteristics}) is equivalent to
$$\chi(S)+\chi(S')-2n\geq2k\cdot \chi(R_+(\ga)).$$
This directly contradicts Lemma \ref{lem: inequality on the Euler characteristics}, since, by definition, we have
$$n=\frac{1}{2}|S\cap \ga|.$$
\epf

\blem\label{lem: inequality on the Euler characteristics}
Suppose $(M,\ga)$, $S$ and $S'$ are as above in Proposition \ref{prop: distinct direct summand}, then
$$\chi(S)+\chi(S')-|S\cap \ga|<2k\cdot \chi(R_+(\ga)).$$
\elem
\bpf
This is exactly the inequality
$$\chi(S)+\chi(S')+I(S)+I(S')<r(S,t)+r(S',t)$$
in the proof of \cite[Theorem 6.1]{juhasz2010polytope} by Juh\'asz. It is by definition that
$$I(S)=I(S')=-\frac{1}{2}|S\cap \ga|,$$
and Juh\'asz also proved in Theorem 6.1 that
$$r(S,t)+r(S',t)=2k\cdot \chi(R_+(\ga)).$$
Note $2k$ in this paper corresponds to $k$ in his paper.
\epf

\brem
Note the assumption in Proposition \ref{prop: distinct direct summand} is slightly stronger than the hypothesis of Theorem 6.1 in Juh\'asz \cite{juhasz2010polytope}, i.e., we require both being reduced and containing no essential product disks, while Juh\'asz only required being reduced. The difference between the two setups is some special family of balanced sutured manifolds, which are reduced but also contain essential product disks. By Lemma 2.13 in \cite{juhasz2010polytope}, there are only two such balanced sutured manifolds, namely the product sutured manifolds $([-1,1]\times F,\{0\}\times \partial{F})$, where $F$ is a sphere with two or three disks removed. It is also worth mentioning that the requirement of containing no essential product disks. Clearly, the two special product sutured manifolds described above are counterexamples to Theorem 6.1 in Juh\'asz \cite{juhasz2010polytope}. The small error made in the proof of Theorem 6.1 in his paper is that, at some point, he used the assumption of being reduced and applied Lemma 2.13 in \cite{juhasz2010polytope} to rule out essential product disks from the discussion, but he didn't exclude the two special product sutured manifolds from the hypothesis, where clearly Lemma 2.13 in \cite{juhasz2010polytope} failed.
\erem

\bcor\label{cor: maximal dimension}
Suppose $(M,\ga)$ is a balanced sutured manifold with $H_2(M)=0$. Suppose further that $(M,\ga)$ is taut, horizontally prime, reduced, and free of essential product disks. Then, the dimensions of $PM_a(M,\ga)$ and $PI_a(M,\ga)$ (see definition \ref{defn: polytope}) are both ${\rm dim}_{\mathbb{Q}} H^2(M,\partial{M};\mathbb{Q})$.
\ecor
\bpf
We prove in the monopole settings. Pick any $\al\in H_2(M,\partial{M})$ so that $\al\neq0$. As in the proof of Proposition \ref{prop: distinct direct summand}, we can find two properly embedded surfaces $S$ and $S'$ in $M$ representing $\al$ and $-\al$, a closure $Y$ of $(M,\ga)$ and suitable closed surfaces $\bar{S}$ and $\bar{S}'$ originate from $S$ and $S'$, respectively. Let $b$ be a homogenous elements in the top grading induced by $\bar{S}$ and let $c$ be a homogenous elements in the top grading induced by $\bar{S}'$. Suppose further that they are supported by spin${}^c$ structures $\mathfrak{s}_b$ and $\mathfrak{s}_c$ on $Y$. We know from the definition that
$$\rho(b,c)(\al)=\frac{1}{2}[c_1(\mathfrak{s}_b)-c_1(\mathfrak{s}_c)](\al).$$

We claim that $\rho(b,c)(\al)\neq0$. Suppose the contrary, then we know that
\begin{equation}\label{eq: contrary 1}
c_1(\mathfrak{s}_c)[\bar{S}]=c_1(\mathfrak{s}_b)[\bar{S}]\geq 2g(\bar{S})-2-2(s-t)-2m,
\end{equation}
and
\begin{equation}\label{eq: contrary 2}
c_1(\mathfrak{s}_c)[\bar{S}']=2g(\bar{S}')-2.
\end{equation}
Here, $s,t,m$ are constants in the construction of $\bar{S}$ and $\bar{S}'$, as in the proof of Proposition \ref{prop: distinct direct summand}. Then, formulae (\ref{eq: contrary 1}) and (\ref{eq: contrary 2}) lead to exactly the same contradiction as in the proof of Proposition \ref{prop: distinct direct summand}.

Since $\rho(b,c)(\al)\neq0$ and $\al\in H_2(M,\partial{M})$ is chosen arbitrarily, we conclude that the polytope $PM_a(M,\ga)$ must have maximal possible dimension.
\epf

\bcor
Suppose $(M,\ga)$ is a taut balanced sutured manifold with $H_2(M)=0$. Suppose further that $(M,\ga)$ is free of essential horizontal surfaces, and is free of non-trivial product disks or product annuli. Then
$${\rm rk}(\shm(M,\ga))\geq {\rm dim}_{\real}H^2(M,\partial{M};\real)+1$$
and
$${\rm dim}_{\mathbb{C}}(\shi(M,\ga))\geq {\rm dim}_{\real}H^2(M,\partial{M};\real)+1$$
\ecor

\bpf
To have dimension $d={\rm dim}_{\real}H^2(M,\partial{M};\real)$, there must be at least $d+1$ points inside the polytope and we are done.
\epf

Now we are ready to prove Theorem \ref{thm: SHM and SHI bounds depth} as stated in the introduction.
\bpf[Proof of Theorem \ref{thm: SHM and SHI bounds depth}]
This is essentially the proof of Proposition 7.6 in Juhasz \cite{juhasz2010polytope} but carried out in the monopole or the instanton settings. We will present the proof in the monopole settings, and the instanton case follows from a similar argument.

First suppose $k=0$. By Proposition 2.16 and Proposition 2.18 in Juh\'asz \cite{juhasz2010polytope}, we can perform a sutured manifold decomposition on $(M,\ga)$ to obtain $(M',\ga')$ that is taut, reduced, and horizontally prime. By Lemma 5.1 in Juh\'asz \cite{juhasz2010polytope}, $H_2(M')=0$. By Proposition 6.6 in Kronheimer and Mrowka \cite{kronheimer2010knots} and Lemma \ref{lem: decompose along product annuli} in this paper, 
$$1\leq{\rm rk}(SHM(M',\ga'))\leq {\rm rk}(SHM(M,\ga))<2.$$
Hence, ${\rm rk}(SHM(M',\ga'))=1$. By Proposition \ref{prop: distinct direct summand}, this implies that $H_2(M',\partial M')=0$ and consequently, $H_1(\partial{M}')=0$, which means $\partial{M}'$ is a sphere. ($H_2(M')=0$ implies that $\partial{M}'$ is connected.) Since $M'$ is irreducible, $M'$ must be a $3$-ball, and $\ga'$ must be connected due to tautness. So, $(M',\ga')$ is a product sutured manifold and so is $(M,\ga)$. This implies that $d(M,\ga)=0$.

Now we assume that the conclusion of the proposition holds for $k-1$, and next, we prove it for $k$. This part of the proof is exactly the same as in Juh\'asz \cite{juhasz2010polytope}, so we only sketch as follows: if $(M,\ga)$ is not horizontally prime, then we could decompose along non-boundary-parallel horizontal surfaces and get a disjoint union of balanced sutured manifolds. Each component has a sutured monopole Floer homology of rank at most $2^{k}$, and thus inductive hypothesis applies, and we conclude that $d(M,\ga)\leq 2k$. If $(M,\ga)$ is horizontally prime, we can perform a sutured manifold decomposition to make it reduced and applying Proposition \ref{prop: distinct direct summand} to choose a suitable decomposition surface so that a (second) sutured manifold decomposition along the chosen surface will reduce the dimension by at least a half. Then, the inductive hypothesis applies, and again we conclude that $d(M,\ga)\leq 2k$. This concludes the proof of theorem \ref{thm: SHM and SHI bounds depth}.
\epf

As a corollary to the above proposition, we offer a new proof to the fact that the monopole and instanton knot Floer homology constructed by Kronheimer and Mrowka \cite{kronheimer2010knots} detects fibred knots in $S^3$.

\bcor
Suppose $K\subset S^3$ is a knot. Then, the following three things are equivalent.

(1) ${\rm rk}(KHM(S^3,K,g(K)))=1$.

(2) ${\rm rk}(KHI(S^3,K,g(K)))=1.$

(3) $K$ is a fibred knot.
\ecor
\bpf
We only prove that (1) and (3) are equivalent. From Kronheimer and Mrowka \cite{kronheimer2010knots}, we know that
$$KHM(S^3,K)=SHM(S^3(K),\Ga_{\mu})$$
where $S^3(K)$ is the knot complement and $\Ga_{\mu}$ is the suture consisting of two meridians on $\partial S^3(K)$. Pick a minimal genus Seifert surface, $S\subset S^3$ of $K$. We know that the decomposition
$$(S^3(K),\Ga_{\mu})\stackrel{S}{\leadsto}(M,\ga)$$
is taut, and
$$KHM(S^3,K,g(K))\cong SHM(M,\ga).$$
Thus, the corollary follows from Theorem \ref{thm: SHM and SHI bounds depth}.
\epf

\bpf[Proof of Corollary \ref{cor: product tangle detection}.]
By Lemma 7.10 in Xie and Zhang \cite{xie2019tangles}, we have an isomorphism
$$SHM(M,\ga,T)\cong SHI(M_T,\ga_T),$$
where $(M_T,\ga_T)$ is some balanced sutured manifold arising from the triple $(M,\ga,T)$ as explained in Section 7 in \cite{xie2019tangles}, and $SHI(M_T,\ga_T)$ is the usual sutured instanton Floer homology defined by Kronheimer and Mrowka in \cite{kronheimer2010knots}. From the description of $M_T$ in \cite{xie2019tangles}, we know that $H_2(M_T)=H_2(M\backslash T)$ and hence Theorem \ref{thm: SHM and SHI bounds depth} applies.
\epf

\section{Application to knots and links}\label{sec: application to links}
\subsection{Thurston-norm detection}
In this subsection, we prove Theorem \ref{thm: Thurston norm detection}. We only work in the monopole settings, and the proof in the instanton settings is exactly the same. First, we need some preparations.

\blem\label{lem: choose a suitable representative}
Suppose $(M,\ga)$ is a taut balanced sutured manifold so that $M$ is boundary-incompressible and the boundary of $M$ consists of a few tori. Suppose further that $\al\in H_2(M,\partial M)$ is a non-zero second relative homology class. Then, there is a properly embedded surface $S\subset M$ with the following properties.

(1) $[S,\partial{S}]=\al\in H_2(M,\partial M).$

(2) $\chi(S)=-x(\al)$.

(3) For any component $\Sigma$ of $\partial{M}$, $S\cap \Sigma$ consists of a disjoint union of coherently oriented non-separating simple closed curves on $\Sigma$. 

(4). $S$ is incompressible.
\elem

\bpf
Pick a surface $S$ so that $[S,\partial S]=\al\in H_2(M,\partial M)$, and $x(S)=x(\al)$. Since $M$ is irreducible, we can assume that there is no spherical component of $S$. We have assumed that $M$ is boundary-incompressible, so we can also assume that there is no disk component of $S$. Thus, we have 
$x(\al)=x(S)=-\chi(S)$.

To achieve condition (3) in the conclusion of the lemma, if a component $\al$ of $\partial{S}$ bounds a disk on $\partial{M}$, then we can cap off $\al$ using the disk it bounds on $\partial{M}$. Capping off by a disk does not increase the norm, so we can assume that the surface $S$ does not have a boundary component that bounds a disk on $\partial{M}$. 

Pick a component $\Sigma$ of $\partial{M}$. By assumption, $\Sigma$ is a torus. Since no component of $S\cap \Sigma$ bounds a disk, we know that $S\cap \Sigma$ consists of a disjoint union of parallel non-separating simple closed curves. If two components of $S\cap \Sigma$ are adjacent on $\Sigma$ but are oriented reversely, we can glue the annulus, which they co-bound on $\Sigma$, to $S$. After possible compressions and throwing away any spherical or disk components arising from the compression, we still call the resulting surface $S$. Note gluing annuli, performing compressions, and throwing away spherical, and disk components do not increase the norm. Thus, we conclude the proof of Lemma \ref{lem: choose a suitable representative}.
\epf

\blem\label{lem: decomposing along S is taut}
Suppose $(M,\ga)$ is a taut balanced sutured manifold so that $M$ is boundary-incompressible, and the boundary of $M$ consists of a few tori. Suppose $\al\in H_2(M,\ga)$ is a non-zero second relative homology class, and $S$ is a properly embedded surface inside $M$ satisfying conditions (1)-(3) in Lemma \ref{lem: choose a suitable representative}. Then, the decomposition of $(M,\ga)$ along $S$ is taut.
\elem

\bpf
Suppose the sutured manifold decomposition of $(M,\ga)$ along $S$ yields $(M',\ga')$, then we can regard $M'$ as a submanifold of $M$. Suppose $\Sigma$ is a component of $\partial{M}$, then by assumption, $\Sigma$ is a torus. If $\Sigma\cap S=\emptyset$, then $\Sigma$ is also a component of $\partial{M}$. Thus, we have
$$\ga'\cap \Sigma=\ga\cap\Sigma,~{\rm and}~R_{\pm}(\ga')\cap \Sigma=R_{\pm}(\ga)\cap V.$$

If $\Sigma\cap S\neq\emptyset$, then $\partial{M'}\cap \Sigma$ consists of a disjoint union of annuli, which, regardless of the orientations, are bounded by pairs of parallel curves in $\partial{S}\cap \Sigma$. Let $A\subset \Sigma$ be a component of $\partial{M}'\cap \Sigma$. There are two cases, depending on the intersection of the suture $\ga$ with the surface $S$. In both cases, it is straightforward to check how $\ga'$ looks like.

\begin{figure}[h]
\centering
    \begin{overpic}[width=5in]{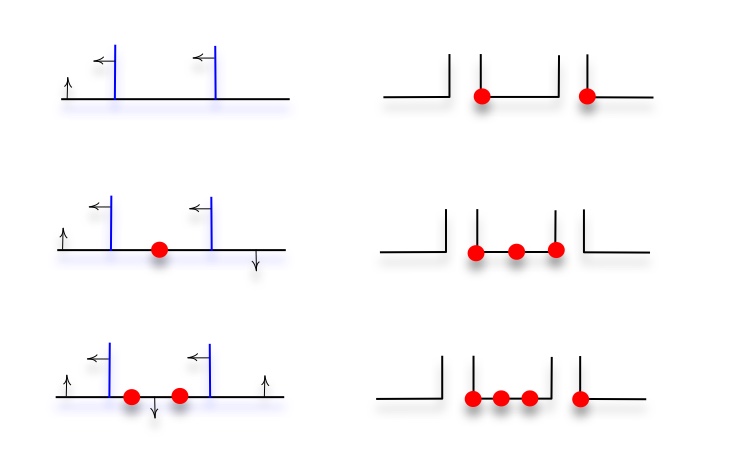}
    \put(20,44){$A$}
    \put(5,46){$\Sigma$}
    \put(14,55){$S$}
    \put(27,55){$S$}
     \put(20,30){$\ga$}
\end{overpic}
\vspace{0.05in}
\caption{Left, before the decomposition. The vertical (blue) curves represent $S$, and the horizontal curves represent $\Sigma$. The (red) dots represent the suture $\ga$. Right, after the decomposition. The (red) dots represent the suture $\ga'$.}\label{fig: decomposition 1}
\end{figure}
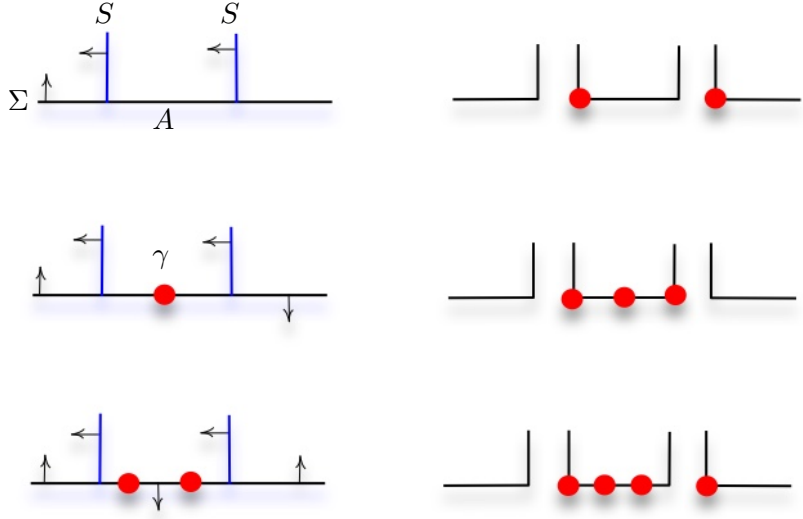

{\bf Case 1}. $(S\cap \Sigma)\cap (\ga\cap \Sigma)=\emptyset$. In this case, the annulus $A$ possibly contains multiple components of $\ga$, and they remains in $\ga'$. Thus, $A$ may contain either one or three components of $\ga'$ (note $\ga\cap \Sigma$ has two components), and each component of $\ga'$ is parallel to $\partial{A}$. See Figure \ref{fig: decomposition 1}.

\begin{figure}[h]
\centering
    \begin{overpic}[width=5in]{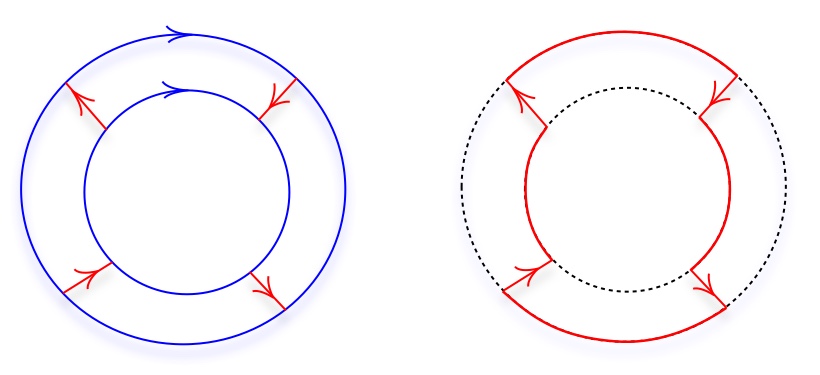}
    \put(20,0){$A$}
    \put(11,21){$\partial{A}\subset\partial S$}
    \put(12,33){$\ga$}
    \put(74,0){$A$}
\end{overpic}
\vspace{0.05in}
\caption{Left, before the decomposition. The two (blue) circles are the boundary of $A$. The (red) arcs represent the suture $\ga$. Right, after the decomposition. The two dashed circles are the boundary of $A$, and the (red) solid curve represent the suture $\ga'$.}\label{fig: decomposition 2}
\end{figure}

{\bf Case 2}. $(S\cap \Sigma)\cap (\ga\cap \Sigma)\neq\emptyset$. In this case, $\ga\cap A$ consists of an even number of essential arcs in $A$, and adjacent arcs are oriented oppositely. Then, after the decomposition, $A$ contains exactly one components of $\ga'$, and this component is parallel to $\partial A$. See Figure \ref{fig: decomposition 2}.

Suppose $A_1,...,A_n$ are all the annular components of $\partial{M}'\cap \partial{M}$ that contain three components of $\ga'$. Push the interiors of $A_i$ into the interior of $M'$ to make them properly embedded. Then, we can perform a sutured manifold decomposition on $(M',\ga')$, along the surface $A_1\cup...\cup A_n$, after the pushing off. The resulting balanced sutured manifold $(M'',\ga'')$ is a disjoint union:
$$(M'',\ga'')=(M''',\ga''')\cup (V_1,\ga_1^4)\cup...\cup (V_n,\ga_n^4).$$
Here, for $i=1,...,n$, $V_i$ is a framed solid torus and $\ga_i^4$ is the suture on $\partial V$ consisting of four longitudes. With $\mathbb{Q}$ coefficients, we know from Li \cite{li2018contact} that
$$SHM(V_i,\ga_i^4)\cong \mathbb{Q}^{2}.$$
From Proposition 6.9 in Kronheimer and Mrowka \cite{kronheimer2010knots}, we know
$$SHM(M',\ga')\cong SHM(M'',\ga'')\cong SHM(M''',\ga''')\otimes_{\mathbb{Q}}\mathbb{Q}^{2^n}.$$
Thus, $(M',\ga')$ is taut if and only if $(M''',\ga''')$ is. 

It is then suffice to prove that $(M''',\ga''')$ is taut. Note we can regard $M'''=M'$, and, thus, we can assume 
$$M'''=M\backslash {\rm int}(N(S)).$$

Let $S_{\pm}$ be parallel copies of $S$ in $\partial{N(S)}$, then $S_{\pm}$ are part of the boundary of $M'''$. Let $\Sigma_{\partial}$ be the union of components of $\partial{M}$ which are disjoint from $S$, and let 
$$F_{\pm}=\Sigma_{\partial}\cap R_{\pm}(\ga).$$
Then, we can describe $R_{\pm}(\ga''')$ as follows:
$$R_{\pm}(\ga''')=F_{\pm}\sqcup S_{\pm}.$$

By assumption, both $F_{\pm}$ and $S_{\pm}$ are incompressible and norm-minimizing in $M$, hence they are also incompressible and norm-minimizing in $M'''$. The fact that $M$ is irreducible implies that $M'''$ is the same. Thus, we conclude the proof of lemma \ref{lem: decomposing along S is taut}.
\epf

\bcor\label{cor: decomposing along S and -S are both taut}
Suppose $(M,\ga)$ and $S$ are the same as in the hypothesis of Lemma \ref{lem: decomposing along S is taut}. Then the decomposition of $(M,\ga)$ along $-S$ is also taut.
\ecor

\bpf[Proof of Theorem \ref{thm: Thurston norm detection}]
For any $\al\in H_2(M,\partial M)$, pick a surface $S$ as in Lemma \ref{lem: choose a suitable representative}. Then, by corollary \ref{cor: decomposing along S and -S are both taut}, the decomposition of $(Y(L),\Ga_{\mu})$ along $S$ and $-S$ are both taut. Suppose 
$$n=\frac{1}{2}|S\cap \ga|,$$
then we know, from condition (3) in the statement of Lemma \ref{lem: choose a suitable representative}, that
$$n=\sum_{i=1}^r|\lgl\al,\mu_i\rgl|.$$

We possibly need to perform a stabilization on $S$ to achieve admissibility. Suppose the $S^{+m}$ is obtained from $S$ by perform $m$ many positive stabilizations on $S$. Here, $m=0$ or $1$. Note $m=0$ means that the original $S$ is admissible, and we take $S^0=S$. From Lemma \ref{lem: adjunction inequality in monopoles}, Lemma \ref{lem: decomposition and pm stabilization}, and Lemma \ref{lem: decomposition theorem reformulated}, we know that
$$\shm(Y(L),\Ga_{\mu},S^{+m}, g_c)\neq 0,~{\rm and}~\shm(Y(L),\Ga_{\mu},S^{+m}, i)=0~{\it for}~i>g_c.$$
Here, we have $g_c=\frac{1}{2}(-\chi(S)+n+m).$

Similarly, we have
$$\shm(Y(L),\Ga_{\mu},(-S)^{+m}, g_c)\neq 0,~{\rm and}~\shm(Y(L),\Ga_{\mu},(-S)^{+m}, i)=0~{\it for}~i>g_c.$$
Note $(-S)^{+m}=-(S^{-m})$, and, hence,
$$\shm(Y(L),\Ga_{\mu},S^{-m}, -g_c)=\shm(Y(L),\Ga_{\mu},(-S)^{+m}, g_c)\neq 0,$$
and
$$\shm(Y(L),\Ga_{\mu},S^{-m}, i)=0~{\it for}~i<-g_c.$$

Applying proposition \ref{prop: general grading shifting property}, we know that
\beq
\shm(Y(L),\Ga_{\mu},S^{+m}, m-g_c)&=\shm(Y(L),\Ga_{\mu},(-S)^{-m}, g_c-m)\\&
=\shm(Y(L),\Ga_{\mu},(-S)^{+m}, g_c)\\
&\neq 0,
\eeq
and
$$\shm(Y(L),\Ga_{\mu},S^{+m}, i)=0~{\it for}~i<m-g_c.$$

From the definition of the function $y(\cdot)$ in Definition \ref{defn: the function y} and the construction of the canonical decomposition of sutured monopole Floer homology in \ref{defn: polytope}, we know that
\beq
y(\al)&=\max\{i~|~\shm(Y(L),\Ga_{\mu},S^{+m}, i)\neq0\}-\min\{i~|~\shm(Y(L),\Ga_{\mu},S^{+m}, i)\neq0\}\\
&=g_c-(m-g_c)\\
&=2g_c-m\\
&=-\chi(S)+n+m-m\\
&=x(\al)+\sum_{i=1}^r|\lgl\al,\mu_i\rgl|.
\eeq
This concludes the proof of Theorem \ref{thm: Thurston norm detection}.
\epf

Suppose $L\subset S^3$ is a link with $r$ components. Then we know that $H_2(S^3(L),\partial{S^3(L)})\cong\intg^r$. Thus there is a $\intg^r$ grading on $KHM(S^3,K)$ and $KHI(S^3,K)$, according to the proof of lemma \ref{lem: homogenous element exists}. Here $KHM$ and $KHI$	 are the monopole and instanton knot Floer homologies introduced by Kronheimer and Mrowka in \cite{kronheimer2010knots}. This leads to the following question.
\begin{quest}
Can we recover the multi-variable Alexander polynomial using the $\intg^r$ grading on $KHM$ or $KHI$?
\end{quest}

\subsection{Minus version for links}\label{sec:Minus version for links }
Suppose $Y$ is a closed oriented $3$-manifold and $L \subset Y$ is an oriented link. Let $L_1$, $L_2$,...,$L_r$ be the components of $L$. We assume further that each component of $L$ is null-homologous in $Y$. Thus, for $i=1,...,r$, we can find (and fix) a Seifert surface $S_i\subset Y$ for $L_i$. Note $S_i$ possibly intersects with $L_j$, for $j\neq i$. Also, for $i=1,...,r$, let $p_i\in L_i$ be a fixed base point. Let $\textbf{p}=(p_1,...,p_r)$. In this subsection, we construct minus versions of monopole and instanton knot Floer homologies for the triple $(-Y,L,\textbf{p})$.

\brem
Here, we require that each component of $L$ to be null-homologous to fix a Seifert surface for each component of the link. It is possible to weaken this condition by simply requiring that the whole link $L$ represents the zero class in $H_1(Y)$ and fix a Seifert surface $S$ for it. The construction in this subsection can be easily adapted to the more general setup.
\erem

Let $Y(L)=Y\backslash N(L)$ be the knot complement and let $T_i$ be the boundary component of $Y(L)$ corresponding to the knot $L_i$. The Seifert surface $S_i$ induces a framing on $T_i$. We call the longitude $\lambda_i$ and the meridian $\mu_i$. Note $\lambda_i$ is oriented in the same way as $L_i$ and $\mu_i$ is oriented so that $\lambda_i\cdot\mu_i=1$ on $\partial Y(L)$. For $\textbf{n}=(n_1,...,n_r)\in\intg^r$, let $\Ga_{\textbf{n}}$ be the suture on $\partial{Y(L)}$ so that $\Ga_{\textbf{n}}\cap T_i$ consists of two parallel simple closed curves of class $\pm[\mu_i-n\lambda_i]$. We have the following lemma.

\brem
It seems that the choice of base points $\textbf{p}$ does not appear in the above setup. However, $\textbf{p}$ helps to resolve the ambiguity arising from the choice of the link complements. Since this issue is fully clarified in Baldwin and Sivek \cite{baldwin2015naturality} and Li \cite{li2019direct}, we won't discuss it anymore in this paper.
\erem

\brem
The construction in this subsection originates from \cite{etnyre2017sutured}. The authors work with knots and Heegaard Floer theory in their paper. A parallel construction for knots was made in monopole and instanton theory by the second author in \cite{li2019direct}.
\erem

\blem
Suppose $\textbf{n}=(n_1,...,n_r)\in(\intg_+)^r$. Let $\textbf{n}'$ be obtained from $\textbf{n}$ by replacing $n_i$ with $n_i+1$, and let $\textbf{n}''$ be obtained from $\textbf{n}$ by replacing $n_i$ with $+\infty$, then there are exact triangles: (when all signs are positive or all signs are negative)
\begin{equation*}
\xymatrix{
\shm(-Y(L),-\Ga_{\textbf{n}})\ar[rr]^{\psi_{\pm,\textbf{n},i}}&&\shm(-Y(L),-\Ga_{\mathbf{n'}})\ar[dl]^{\psi_{\pm,\textbf{n}',i}}\\
&\shm(-Y(L),-\Ga_{\mathbf{n''}})\ar[ul]^{\psi_{\pm,\textbf{n}'',i}}&
}    
\end{equation*}
Here, $\psi_{\pm,\textbf{n},i}$ are the map associated to a positive or negative by-pass attached to $Y(L),\Ga_{\textbf{n}}$, which is performed on the boundary component $T_i$ of $Y(L)$. (See Subsection \ref{subsec: by passes}.) 

There are similar exact triangles in the instanton settings.
\elem
\bpf
This is a direct application of Theorem \ref{thm_by_pass_attachment}. For more details, readers are referred to Section 2 of Li \cite{li2019direct}.
\epf

To use a better notation, for $i=1,...,r$, let $\mathbf{e}^i=(0,...,1,...,0)\in\intg^r$ be the vector whose entries are all $0$ except of being $1$ on the $i$-th place.

\blem\label{lem: commutative diagram for direct limit}
For any $\textbf{n}\in(\intg_+)^r$ and $i,j\in\{1,...,r\}$, we have the following commutative diagram:
\begin{equation*}
\xymatrix{
\shm(-Y(L), -\Ga_{\textbf{n}})\ar[rr]^{\psi_{-,\textbf{n},i}}\ar[dd]^{\psi_{-,\textbf{n},j}}&&\shm(-Y(L),-\Ga_{\textbf{n}+\textbf{e}^i})\ar[dd]^{\psi_{-,\textbf{n}+\textbf{e}^i,j}}\\
&&\\
\shm(-Y(L), -\Ga_{\textbf{n}+\textbf{e}^j})\ar[rr]^{\psi_{-,\textbf{n}+\textbf{e}^j,i}}&&\shm(-Y(L),-\Ga_{\textbf{n}+\textbf{e}^i+\textbf{e}^j})
}
\end{equation*}
\elem

\begin{proof}
When $i=j$, the diagram commutes obviously. Assume $i\neq j$. As explained in Subsection \ref{subsec: by passes}, the by-pass maps, $\psi_{\pm,\textbf{n},i}$, ultimately come from contact handle attaching maps. Since the by-passes corresponding to vertical and horizontal maps happen on different boundary components of $Y(L)$, the corresponding contact handle attachments commute, and so do the by-pass maps.
\end{proof}

\bdefn \label{defn: minus version for links}
We define the minus version of monopole link Floer homology of a based link $L \subset -Y$, which is denoted by $\underline{\rm KHM}^-(-Y,L,\textbf{p})$ , to be the
direct limit of the direct system $$\{\psi_{-,\textbf{n}, i}: \shm (-Y(L), \Gamma_\textbf{n})\rightarrow \shm (-Y(L), \Gamma_{\textbf{n}+ \textbf{e}^i}), \textbf{n}\in (\mathbb{Z}_+)^r, i\in {1,...,r}\}.$$

We define $\khi(-Y,L,\textbf{p})$ in a similar manner.
\edefn

\blem \label{lem: commutative diagram for U maps}
For any $\textbf{n}\in(\intg_+)^r$ and $i,j\in\{1,...,r\}$, we have the following commutative diagram:
\begin{equation*}
\xymatrix{
\shm(-Y(L), -\Ga_{\textbf{n}})\ar[rr]^{\psi_{-,\textbf{n},i}}\ar[dd]^{\psi_{+,\textbf{n},j}}&&\shm(-Y(L),-\Ga_{\textbf{n}+\textbf{e}^i})\ar[dd]^{\psi_{+,\textbf{n}+\textbf{e}^i,j}}\\
&&\\
\shm(-Y(L), -\Ga_{\textbf{n}+\textbf{e}^j})\ar[rr]^{\psi_{-,\textbf{n}+\textbf{e}^j,i}}&&\shm(-Y(L),-\Ga_{\textbf{n}+\textbf{e}^i+\textbf{e}^j})
}
\end{equation*}

There are similar commutative diagrams in the instanton settings.
\elem

\bpf
If $i\neq j$, this follows from exactly the same argument as in the proof of Lemma \ref{lem: commutative diagram for direct limit}. If $i=j$, this follows from the proof of the same type of commutative diagram in the construction of minus versions for knots by the second author. See Corollary 2.22 in Li \cite{li2019direct}.
\epf

\begin{defn}\label{defn: U maps for links}
For any fixed $i\in\{1,...,r\}$, The set of maps
$$\{\psi_{+,\textbf{n},i}:\shm(-Y(L), -\Ga_{\textbf{n}})\ra \shm(-Y(L), -\Ga_{\textbf{n}+\textbf{e}^i})\}$$
induces a map
$$U_i: \khm(-Y,L,\textbf{p})\ra \khm(-Y,L,\textbf{p}),$$
which we call the {\it $i$-th $U$ map}.

We define
$$U_i: \khi(-Y,L,\textbf{p})\ra \khi(-Y,L,\textbf{p})$$
in a similar manner.
\end{defn}

\bprop\label{prop: U maps all comuute}
For any $i,j\in\{1,...,r\}$, the maps $U_i$ and $U_j$ commute with each other.
\eprop

\begin{proof}
The proof is exactly the same as the proof of Lemma \ref{lem: commutative diagram for direct limit}.
\end{proof}

Next, we construct a $\mathbb{Z}^r$ grading on $\khm(-Y,L,\textbf{p})$, based on the chosen Seifert surfaces $S_1,...,S_r$ of $L_1,...,L_r$. 

Recall that $L\subset Y$ has components $L_1, L_2,..., L_r$, and $L_i$ has a Seifert surface $S_i$ that could possibly intersect other components of the link. By a slight abuse of notation, let $S_i$ also denote the intersection of the original Seifert surface with the link complement $Y(L)$. Thus, the boundary of $S_i$ consists of a longitude on $T_i$ and a few (possibly none) meridians on $T_j$, for $j\neq i$. For fixed $i\in\{1,...,r\}$ and $\textbf{n}=(n_1,...,n_r)\in(\intg_+)^r$, let $S_{i,\textbf{n}}$ be the isotopy of $S_i$ so that $\partial{S_{i,\textbf{n}}}$ has the least possible intersections with the suture $\Ga_{\textbf{n}}\subset \partial{Y(L)}$. This means that the longitudinal boundary component of $S_{i,\textbf{n}}$ intersects $\Ga_{\textbf{n}}$ at $2n_i$ points, and each meridional boundary component of $S_{i,\textbf{n}}$ intersects $\Ga_{\textbf{n}}$ at two points. Applying the construction of gradings in Subsection \ref{subsec: construction of gradings}, the surface $S_{i,\textbf{n}}$, or its stabilizations, give rise to a grading on $\shm (-Y(L), -\Gamma_\textbf{n})$. We then have the following proposition.

\bprop\label{prop: grading shifting for by-pass maps}
Fix any $i\in\{1,...,r\}$ and $\textbf{n}\in(\intg_+)^r$. If $n_i$ is even, then, for any $j\in \intg$, we have
$$\psi_{\pm, \textbf{n}, i}(\shm (-Y(L), -\Gamma_{\textbf{n}}, S_{i,\textbf{n}}^{\pm}, j))\subset \shm (-Y(L), -\Gamma_{\textbf{n}+ {\textbf{e}^i}}, S_{i,\textbf{n}+\textbf{e}^i},j)$$
If $n_i$ is odd, then, for any $j\in \mathbb{Z}$,
$$\psi_{\pm, \textbf{n}, i}(\shm (-Y(L), -\Gamma_\textbf{n}, S_{i,\textbf{n}}^{\pm 2}, j))\subset \shm (-Y(L), -\Gamma_{\textbf{n}+ {\textbf{e}^i}}, S_{i,\textbf{n}+ \textbf{e}^i}^{\pm},j)$$

Furthermore, for any $k\neq i$, the maps $\psi_{\pm,\textbf{n},k}$ preserve the gradings associated to $S_{i,\textbf{n}}$ and its stabilizations. 

Similar statements hold for the instanton settings.
\eprop

\bpf
For $\psi_{\pm, \textbf{n}, i}$, the proof is exactly the same as the proof of Proposition 5.5 in Li \cite{li2019direct}. For $\psi_{\pm, \textbf{n}, k}$ with $k\neq i$, note $S_{i,\textbf{n}}$ has a few meridional components on the $T_k$, so the by-passes, which corresponds to the maps $\psi_{\pm, \textbf{n}, k}$, can actually be made disjoint from $S_{i,\textbf{n}}$.
\epf

Similar to the constructions in Section 5 of Li \cite{li2019direct}, if $n_i$ is odd, let $S^{\tau}_{i,\textbf{n}}$ be just $S_{i,\textbf{n}}$, and if $n_i$ is even, let $S^{\tau}_{i,\textbf{n}}$ be a negative stabilization of $S_{i,\textbf{n}}$ performed near $T_i$. We can use $S^{\tau}_{i,\textbf{n}}$ to define a grading on $\shm(-Y(L),\Ga_{\textbf{n}})$. We also need to perform a grading shift. Let
$$\shm(-Y(L),-\Ga_{\textbf{n}},S^{\tau}_{i,\textbf{n}},j)\{\sigma_i\}=\shm(-Y(L),-\Ga_{\textbf{n}},S^{\tau}_{i,\textbf{n}},j+\lfloor\frac{n}{2}\rfloor),$$
where $\lfloor x\rfloor$ is to take the maximal integer which is no larger than $x$.

\bprop
Using the grading $\shm(-Y(L),-\Ga_{\textbf{n}},S^{\tau}_{i,\textbf{n}},j)\{\sigma_i\}$, we can construct a $\intg$-grading on $\khm(-Y,L,\textbf{p})$. The $i$-th $U$ map, $U_i$, drops the grading by $1$, and all other $U$ maps, $U_k$ with $k\neq i$, preserve the grading.

Furthermore, all Seifert surfaces, $S_1,...,S_r$, together induce a $\intg^r$ grading on $\khm(-Y,L,\textbf{p})$, which we write as
$$\khm(-Y,L,\textbf{p},\textbf{j}).$$
Here, $\textbf{j}\in\intg^r$ denote a multi-grading.
As a result, together with the commutativity of all $U$ maps in Proposition \ref{prop: U maps all comuute}, there is an $\mathcal{R}[U_1,...,U_r]$ module structure on $\khm(-Y,L,\textbf{p})$.

Similar results hold in the instanton settings.
\eprop

\bpf
The first half of the proposition follows from Proposition \ref{prop: grading shifting for by-pass maps}, and the second half of the proposition follows from the proof of Lemma \ref{lem: homogenous element exists}.
\epf



The first computable example is the case of unlinks.

\bprop
Suppose $Y$ is a closed oriented 3-manifold and $L \subset Y$ is an unlink of $r$ components, i.e., there exists an embedded disk $S_i \cong D^2$, for each $i\in\{1,..., r\}$, so that $\partial S_i=L_i$, and all $S_i$ are disjoint from each other. Then, 
$$\khm (-Y, L, \textbf{p}) \cong \shm(-Y(r),-\delta^r) \, \otimes_\mathcal{R} \mathcal{R}[U_1,…,U_r].$$
Here, $\textbf{p}$ is a chosen set of base points, and $(Y(r), \delta^r)$ is the balanced sutured manifold obtained from $Y$ by removing $r$ disjoint 3-balls and picking one simple closed curve on each spherical boundary of $Y(r)$ as the suture.

Similar statements hold in the instanton settings.
\eprop

\bpf
For any $\textbf{n}\in(\intg_+)^r$, we know that $(Y(L),\Ga_{\textbf{n}})$ can be obtained from the disjoint union 
$$(Y(r),\delta^r)\sqcup (S^3(L_1),\Ga_{n_1})\sqcup ...\sqcup (S^3(L_r),\Ga_{n_r})$$
by attaching $r$ many contact $1$-handles (see Definition \ref{defn: one handle}). Each $1$-handle connects some $(S^3(L_i),\Ga_{n_i})$ to $(Y(r),\delta^r)$. As in Subsection \ref{subsec: by passes}, the by-pass maps $\psi_{\pm,\textbf{n},i}$ can be realized as contact handle attaching maps and those contact handles are disjoint form the contact $1$-handles just described above. Hence, under the isomorphism
$$\shm(-Y(L),\Ga_{\textbf{n}})\cong \shm(Y(r),\delta^r)\otimes\shm(S^3(L_1),\Ga_{n_1})\otimes...\otimes \shm(S^3(L_r),\Ga_{n_r}),$$
we have an identification
$$\psi_{\pm,\textbf{n},i}=id\otimes...\otimes\psi_{\pm,n_i}\otimes...\otimes id,$$
where 
$$\psi_{\pm,n_i}:\shm(-S^3(L_i),-\Ga_{n_i})\ra\shm(-S^3(L_i),-\Ga_{n_i+1}).$$
Hence, we are done.
\epf

\bprop
Under the above setups, the direct system stabilizes, that is, for any fixed $j\in \mathbb{Z}$, there exists $N\in \mathbb{Z}$ , so that for all $i\in\{1,...,r\}$ and $\textbf{n}=(n_1,...,n_r) \in(\intg_+)^r$ such that $n_i>N$, we have an isomorphism 
$$\psi_{-,\textbf{n},i}: \shm (-Y(L), -\Gamma_\textbf{n}, S_{i,\textbf{n}}^{\tau}, j)\{\sigma_i\} \cong \shm (-Y(L), -\Gamma_{\textbf{n}+ {\textbf{e}^i}},  S_{i,\textbf{n}+ {\textbf{e}^i}}^{\tau}, j)\{\sigma_i\}.$$
\eprop
\bpf
This follows from exactly the same argument as in the proof of Proposition 5.10 in Li \cite{li2019direct}.
\epf

\bprop
Under the above setups, there exists an integer $N_0$, so
that for any fixed $i\in\{1,...,r\}$ and any multi-grading $\textbf{j}=(j_1,...,j_r)\in \intg^r$ with $j_i< N_0$, the map $U_i$ restricts to an isomorphism
$$U_i:\khm(Y, L, \textbf{p}, \textbf{j}) \cong \khm(Y, L, \textbf{p}, \textbf{j}-\textbf{e}^i).$$
\eprop
\bpf
This follows from exactly the same argument as in the proof of Corollary 5.11 in Li \cite{li2019direct}.
\epf

\bibliography{Index}
\end{document}